\documentclass{commat}
\UseRawInputEncoding

\usepackage[all]{xy}
\usepackage{amsthm,array,amssymb,amscd,amsfonts,latexsym}

\usepackage{enumerate}
\usepackage{mathrsfs}

\usepackage{lscape}

\usepackage{graphics}
\usepackage{graphicx}
\usepackage[all]{xy}
\usepackage{tabularx,mdwtab}
\usepackage{lscape,longtable}

\usepackage{enumerate}

\usepackage{hyperref}

\usepackage{url}

\newcommand\surtab{\small\centering\extrarowheight=1\p@ \def\arraystretch{1.4}\tabcolsep3.6\p@ \doublerulesep=0.09\p@%
 \vskip 3\p@\advance\belowcaptionskip3\p@ \mathsurround\z@ \hyphenpenalty\z@ \doublehyphendemerits\z@ }
\newcommand\tcaption[2][0pt]{%
 \def\@tempa{#2}%
 \addtocounter{table}{1}\surtab\edef\@currentlabel{\thetable}%
 \ifx\@tempa\@empty \hfill \textbf{\tablename}~\thetable\kern#1 \else\textbf{\tablename}~\thetable.\kern 1ex{\boldmath\bfseries #2}\fi
}

\newcommand{\cc}{\raise .4pt \hbox{{$\scriptstyle{\bullet}$}}}

\newcommand{\Tran}{{\rm{Tran}\,}}
\newcommand{\Lin}{{\rm {Lin}\,}}
\newcommand{\re}{{\rm {Re}\,}}
\newcommand{\Tr}{{\rm {Tr}\,}}
\newcommand{\GL}{{\rm {GL}}}

\newcommand{\tto}{\to}
\newcommand{\sgn}{{\rm {sgn}}\,}

\newcommand{\ka}{{\bf k}}
\newcommand{\Ree}{{\bf R}}
\newcommand{\Cee}{{\bf C}}
\newcommand{\Zee}{{\bf Z}}
\newcommand{\Qee}{{\bf Q}}
\newcommand{\ii}{{\rm id}_{V}}

\newcommand{\codim}{{\rm codim}\,}
\newcommand{\rk}{{\rm rk}\,}
\newcommand{\Ker}{{\rm Ker}\,}

\title{%
    Discrete complex reflection groups
    }

\author{%
    Vladimir L. Popov
    }

\authorinfo[%
    V. L. Popov]{
    Steklov Mathematical Institute,
Russian Academy of Sciences, Gub\-kina 8,
Moscow 119991, Russia}{%
    popovvl@mi-ras.ru
    }

\abstract{%
    Here are reproduced slightly edited notes of my lec\-tures on the classification of discrete groups generated by complex reflections of Hermitian affine spaces
delivered  in October of 1980
at the University of Utrecht (see \cite{15} and MR83g:20049).
    }



\VOLUME{30}
\YEAR{2022}
\NUMBER{3}
\firstpage{303}
\DOI{https://doi.org/10.46298/cm.11249}
\begin{document}




\setcounter{tocdepth}{3}
{\footnotesize \tableofcontents}

\section*{\bf Foreword of July 10, 2023}

Below are the notes of five lectures I delivered in October 1980 at the University of Utrecht.\;Their subject matter is my classification of discrete groups generated by complex reflections of Hermitian affine spaces ${\mathbb C}^n$ (named here $r$-{\it groups}).\;The goal of lectures was to present the results of classification, along with the main ideas, statements, and proofs that make it possible to obtain these results.\;Getting the final answer inevitably requires some calculations.\;Being
 quite voluminous, the latter
are not included in the present text:\;I only explained what and how to calculate and gave the appropriate examples.

The prehistory of obtaining this classification
dates back to 1965 when E.\;B.\;Vinberg formulated to me the problem of classifying crystallographic $r$-groups.\;In 1967, I handed over to him my  manuscript \cite{13} containing many of the results and ideas of the desired theory, in particular,  the classification of crystallographic $r$-groups with primitive linear parts and one of the key ingredients of the theory, the operator $S$ (see Theorem \ref{opS}), whose usage later led me to classifying also $r$-groups with imprimitive linear parts thereby
completing solution to the problem and, moreover, extending it to the case of noncrystallo\-graphic $r$-groups.\;He, in turn, passed it to the author of \cite{14} who launched
his research on $r$-groups.\;Later, in Utrecht, I was told that in fact the problem
was posed earlier  by A.\;Borel.

In 1982, the notes of my Utrech lectures were published in the 15th issue of {\it Communica\-tions of Mathema\-tical Institute Rejksuniversiteit Utrecht},  a rather obscure periodical \cite{15}. In those days,  {\sf arXiv} did not yet exist and {\sf TeX} was not used to write mathema\-tical texts that
were instead typed on typewriters.\;As a result, for quite a long time these notes remained known only to a narrow circle of mathematicians.\;However, over time, their popularity increased, the number of references to them grew, and the results of the classification were used again and again.\;Therefore,  making these notes more accessible now becomes justified
since with the advent of {\sf TeX} and {\sf arXiv} technically this became possible.\;D.\;Leites, who took the initiative, organized the conversion of the typewritten text into a .tex file, for which I am sincerely grateful to him.\;I made some minor changes to this file, correcting obvious typos, stylistic and evident mathematical inaccuracies (like lacking some edge weights in the diagrams in Subsection\;\ref{refff}), and also  added  Subsections \ref{exa1}--\ref{exa2}, where  the example of one-dimensional groups is discussed in more detail.

As for the calculations mentioned above and so far kept by me in handwritten form, in principle,
the usage of {\sf arXiv} opens up the possibility of publishing them in the form of appropriate appendixes to these notes.\;This would
simplify  checking the reliability of the classification results and
protect against possible doubts
about their accuracy\footnote{Perhaps due to them, it was suggested in V.\;Goryunov, S.\;H.\;Man, {\it The complex crystallographic groups and symmetries of $J_1$}, Adv.\;Stud.\;Pure\;Math.\;{\bf 43} (2006), 55--72,   a construction of a group, which allegedly is an $r$-group absent in my classification.\;After many years,
in a letter
of April 12, 2023 to me, the first author advised that this construction is erro\-ne\-ous.\;He writes:\,``I realised
later on that the group presentation I came across was just a different pre\-sen\-tation of one of your groups''.\,However, I am not aware of the existence of a publication where the authors of the cited paper made this refutation public.}.\;Therefore, I do not exclude that somewhat later
 will take the advantage of this possibility although the conversion of a large handwritten mathematical text into a .tex file is quite a laborious\;task.

Since 1980, various aspects of the subject of $r$-groups have been explored in a number of papers.\;A full review of these explorations would be a separate endeavor.\;I will not undertake it here and  only mention  two of them:\;obtaining the presentations of such groups by generators and relations
(see  \cite{10}), and exploration of the quotient ${\mathbb C}^n/\Gamma$, where $\Gamma$ is a crystallographic $r$-group.\;If $T$ is the subgroup of all translations contained in $\Gamma$, then ${\mathbb C}^n/\Gamma$ is the quotient of the complex torus $A={\mathbb C}^n/T$ by the natural action of a finite linear $r$-group $K$.\;From the explicit description of $T$ obtained in these lectures, it follows that $A$ is in fact an Abelian variety of a rather special kind (for an a priori proof of this, see \cite{11}); therefore,  ${\mathbb C}^n/\Gamma$ is an algebraic variety.\;A striking achievement of the most recent time is the proof, obtained in \cite{12}, of the long standing conjecture that if  $\Gamma$ is irreducible, then the algebraic variety ${\mathbb C}^n/\Gamma=
A/K$ is a weighted projective space. In particular, this algebraic variety is rational.\;The latter property is considered as a counterpart for $r$-groups of the classical Chevalley's theorem  about freeness of the invariant algebras of finite linear $r$-groups.


\section*{\bf Introduction}\label{S0}


These notes essentially record the content of a~course of five lectures given at the Mathematical Institute of the Rijksuniversiteit Utrecht in October 1980. The aim of these lectures was to develop the theory of discrete groups generated by affine unitary reflections, and in fact, to provide
a~classification of these groups. I made an attempt to give an exposition in such a~way that our results would be comparable with the classical results of the theory of discrete groups generated by affine reflections in a~real Euclidian space, which was developed in the former half of XXth century by Coxeter, Witt, Stiefel and others. Both theories have much in common. However, the classification of complex groups is more complicated (and less geometrical) than the classification over the reals.

The problem of developing the theory of discrete groups generated by affine unitary reflections is a~comparatively old one; I was informed that it was posed by A.~Borel about
1965.\;My general aim during the~lectures was to explain the main ideas and to give proofs only of the theorems that are not of a~strictly technical nature. Therefore, I restricted myself to examples or simply to formulating results in all technical cases (which are, however, not always trivial). A more detailed exposition will appear elsewhere.

\vskip 2mm

I wish to thank the Mathematical Institute of the Rijksuniversiteit Utrecht for its hospitality. I am grateful to Professor T.~A.~Springer on whose initiative these lectures were given and written up. I am also grateful to A.~M.~Cohen for the interesting discussions I had with him and for his help in preparing these notes. My thanks go also to the secretaries of the University of Utrecht, the Netherlands for the careful typing of the manuscript.

\section{\bf Notation and formulation of the problem}\label{S1}

We assume in this section  that the ground field $\ka$ is $\Ree$ or $\Cee $.

\subsection{Notation}\label{1.1}

Let $E$ be \textit{an affine space} over $\ka$, $\dim E=n$, and let $V$ be its \textit{space of translations}. If $v\in V$, we denote by $\gamma _v$ the corresponding translation of $E$, i.e.,
$$
\gamma _v(a)=a+v \quad \text{for each}\quad a\in E.
$$
Let $A(E)$ be the \textit{group of all affine transformations} of $E$ and let
$$
\Tran A(E)=\{\gamma _v \mid v\in V\}.
$$
We denote by $\GL(V) \ltimes V$ \textit{the natural semidirect product} of $\GL (V)$ and $V$. Its elements are pairs $(P,v)$, where $P\in \GL (V)$ and $v\in V$, and the group operations are given by formulas
\begin{align*}
(P,v)(Q,w)&=(PQ,Pw+v),\\
(P,v)^{-1}&=(P^{-1},-P^{-1}v).
\end{align*}
Let
$$
\Lin \colon A(E)\to \GL(V)
$$
be the standard homomorphism defined by formula
$$
\gamma (a+v)=\gamma (a)+(\Lin \gamma )v \quad \text{for any}\quad \gamma \in A(E),\quad a\in E,\quad v\in V.
$$
If we take a~point $a\in E$ as the origin, we obtain an isomorphism
$$
\kappa _a\colon A(E)\to \GL (V)\ltimes V
$$
given by the formula
$$
\kappa _a(\gamma )=(\Lin \gamma ,\gamma (a)-a).
$$
Identifying $A(E)$ and $\GL (V)\ltimes V$ by means of $\kappa _a$, we obtain \textit{the action of} $\GL (V)\ltimes V$ \textit{on} $E$ given by the formula
$$
(P,v)q=a+P(q-a)+v \quad \text{for each}\quad q\in E.
$$
The dependence on $a$ is given by the formula
$$
\kappa _b(\gamma )=\kappa _a(\gamma _{a-b}\gamma \gamma _{b-a})\quad \text{for each}\quad \gamma \in A(E),\quad b\in E.
$$
For every $\gamma \in A(E)$ and $P\in \GL (V)$, we use the notation
\begin{align*}
H_\gamma &=\{a\in E \mid \gamma (a)=a\},\\
H_P&=\{v\in V \mid P_v=v\}.
\end{align*}
These are subspaces of $E$ and $V$, respectively.

Let $\langle\, \cdot \mid \cdot\, \rangle $ be a~positive-definite inner product on $V$, i.e., $V$ is an Euclidian if $\ka=\Ree $ (resp. Hermitian if $\ka=\Cee $) linear space with respect to $\langle \,\cdot \mid \cdot\, \rangle $ (linear in the first
coordinate). Let also
$$
{\rm Iso}(V) : =\{P\in \GL (V)\mid P \text{ preserves }\langle \,\cdot \mid \cdot\, \rangle \};
$$
this is a~compact group. The space $E$ becomes a~Euclidean (resp. Hermitian) affine metric space with respect to the distance given by the formula
$$
\rho (a,b)=\sqrt {\langle a-b \mid a-b \rangle } \quad \text{for each}\quad a,b\in E.
$$

\subsection{Motions and reflections}\label{more}

We say that $\gamma \in A(E)$ is a~\textit{motion} of $E$ if $\gamma $ preserves the distance $\rho $. It is easy to see that $\gamma $ is a~motion if and only if $\Lin \gamma \in {\rm Iso}(V)$.


\textit{An affine reflection} $\gamma \in A(E)$ is an element with the following properties:

\eject

{
\begin{enumerate}[\hskip 7.2mm\rm 1)]\itemsep=-.1ex
\item
$\gamma $ is a~motion,
\item
$\gamma $ has finite order,
\item
$\codim H_\gamma =1$.
\end{enumerate}

}

\noindent \textit{A linear reflection} $R\in \GL (V)$ is an element with the following pro\-per\-ties:
{
\begin{enumerate}[\hskip 7.2mm\rm 1)]\itemsep=-.1ex
\item
$R\in {\rm Iso}(V)$,
\item
$R$ has finite order,
\item
$\codim H_R=1$.
\end{enumerate}

}


\noindent The subspaces $H_\gamma $ and $H_R$ are called \textit{the mirrors} of $\gamma $ and $R$, respec\-ti\-vely.

Sometimes, when it is clear what we are talking about, we shall simply say \textit{reflection}.


If $R$ is a~linear reflection, then the line
$$
\ell_R=\{v\in V \mid v \perp H_R\}
$$
is called \textit{the root line} of $R$. If $v\in \ell_R$ and $\langle v \mid v \rangle =1$, then $Rv=\theta v$, where $\theta \ne 1$ is a~primitive root of 1 (if $\ka=\Ree $, then $\theta =-1$, if $\ka=\Cee $, then $\theta $ may be arbitrary). The pair $(v,\theta )$ completely determines $R$ and every pair $(u,\eta )$ with $\langle u \mid u \rangle =1$ and $\eta \ne 1$ a~primitive root of unity (
$=-1$ if $\ka=\Ree $), can be obtained in such a~way from a~reflection. We write
$$
R=R_{v,\theta }.
$$
Some properties of the reflections are contained in the following
\subsubsection{Proposition}\label{prope}
{\it
Let $\gamma \in A(E)$, $a\in E$, and $\kappa _a(\gamma )=(R,v)$. Then
\begin{enumerate}[\hskip 7.2mm\rm 1)]\itemsep=-.1ex
\item
$\gamma $ is a~reflection if and only if $R$ is a~reflection and $v \perp H_R$.
\item
If $\gamma $ is a~reflection and $R=R_{e,\theta }$ then
$$
H_\gamma =a+H_R+(1-\theta )^{-1}v.
$$
\item
$R_{e,\theta }v=v-(1-\theta )\langle v \mid e \rangle e$.
\item
If $\gamma $ is a~reflection and $\delta $ is a~motion, then $\delta \gamma \delta ^{-1}$ is a~reflection.
\end{enumerate}
}
\vskip 2mm
\begin{proof} It is left to the reader.
\end{proof}


\subsection{Main problem}

We say that a~subgroup $W$ of $A(E)$ is an \textit{$r$-group} if it is discrete and generated by affine reflections.

If $E$ and $E'$ are two affine spaces and $W \subseteq A(E)$ and $W' \subseteq A(E')$ are two arbitrary subgroups, then we say that $W$ and $W'$ are \textit{equivalent} if there exists an affine bijection $\phi \colon E \to E'$ such that
$$
W'=\phi W \phi ^{-1}.
$$

This means that after identifying $E$ and $E'$ by means of an arbitrary fixed isomorphism, the groups $W$ and $W'$, as subgroups of $A(E)$, have to be \textit{conjugate in} $A(E)$.

We want to emphasize here that even when $E$ and $E'$ are affine \textit{metric} spaces, $\phi $ \textit{need not} to be distance preserving.

Our main goal in these lectures is \textit{to classify $r$-groups up to equi\-valence}.

We will show now that in solving this problem one can restrict
to consideration of irreducible groups.

\subsection{Irreducibility}\label{1.4}

Let $W$ be a~subgroup of $A(E)$. We say that $W$ is \textit{reducible} if there exist affine metric spaces $E_j$, where $j=1,\ldots ,m$ for $m \geqslant 2$, and subgroups $W_j$ of $A(E_j)$ such that $W$ is
equivalent to $W_1 \times \cdots \times W_m \subseteq A(E_1 \times \cdots \times E_m)$. Otherwise $W$ is called \textit{irreducible}. Clearly, every group is isomorphic to a~product of irreducible groups (but its decomposition need not be unique).

\subsubsection{Theorem}\label{irred}
{\it  Let $W \subseteq A(E)$ be a~nontrivial subgroup \textup(possibly nondiscrete\textup) generated by affine reflections. Then
\begin{enumerate}[\hskip 7.2mm\rm a)]\itemsep=-.1ex
\item
$W$ is equivalent to a product
$W_1 \times \cdots \times W_m$,
where every group $W_j$ is
irreducible
and
is either generated by affine reflections or trivial \textup(hence, $1$-dimensional\textup), but not all $W_j$'s are trivial.

\item
$W$ is irreducible if and only if $\Lin W$ is an irreducible linear group \textup(generated by linear reflections\textup).

\item
The groups $W_1,\ldots, W_m$
are uniquely defined up to equivalence and numbering.

\item
Every product of the type described in {\rm a)} is a~group generated by reflections.
\end{enumerate}
}

\vskip 2mm

\begin{proof} a) The statement follows from the equality
$$
H_{(\gamma _1,\ldots ,\gamma _m)}=H_{\gamma _1}\times \cdots \times H_{\gamma _m}
$$
(hence $(\gamma _1,\ldots ,\gamma _m)$ is a~reflection if and only if one and only one of
$\gamma _1,\ldots, \gamma_m$ is a~reflection and the others are equal to 1).

b) The \lq\lq if\rq\rq{} part is obvious. Let us prove the \lq\lq only if\rq\rq{} part.

As the group $W$ is generated by reflections, the group $\Lin W$ lies in ${\rm Iso}(V)$. Therefore, $\Lin W$ is a~completely reducible linear group. Let
$$
V=\bigoplus\limits _{j=1}^mV_j,
$$
where
$V_1,\ldots, V_m$ are irreducible $\Lin W$-modules. Consider the sub\-spa\-ces
$$
E_j=a+V_j \quad \text{for each}\quad 1 \leqslant j \leqslant m.
$$
where $a\in E$ is an origin, and let
$$
\pi _j\colon W\tto A(E_j) \quad \text{for}\quad 1 \leqslant j \leqslant m
$$
be the morphism given by the formula
$$
\pi _j(\gamma )=\kappa _a^{-1} \bigl (\Lin \gamma \bigr |_{V_j},\quad p_j(\gamma (a)-\alpha )\bigr )
$$
(here $p_j\colon V \tto V_j$ is the natural projection). Let $W_j=\pi _j(W)$. Then it is not difficult to check that the map
$$
\phi \colon E \tto E_1 \times \cdots \times E_m,
$$
given by the formula
$$
\phi (q)=(a+p_1(q-a),\ldots ,a+p_m(q-a)) \quad \text{for each}\quad q\in E
$$
defines an equivalence of $W$ and
$W_1 \times \cdots \times W_m$.

c) Suppose that $W \subseteq A(E)$ and $W' \subseteq A(E')$ are two
equivalent groups generated by affine reflections and let $\phi \colon E \tto E'$ establish an 
equivalence of these groups. Let $W=W_1 \times \cdots \times W_r$ and let $W'=W_1' \times \cdots \times W_s'$ be decompositions into products of irreducible groups and let
\begin{alignat*}{2}
E&=E_1 \times \cdots \times E_r, &\quad V&=V_1 \oplus \cdots \oplus V_r,\\
E'&=E_1' \times \cdots \times E_s', &\quad V'&=V_1' \oplus \cdots \oplus V_s',
\end{alignat*}
be the corresponding decompositions of the affine spaces and its spaces of translations. We consider $W_j$ and $W'_l$ for all $j$ and $l$ as subgroups of $W$ and $W'$, respectively.

It is clear that $\Lin \psi $ yields an
equivalence of the linear groups $\Lin W$
$\Lin W'$ (in the usual sense). Hence $(\Lin \psi )V_j$ is a~simple $(\Lin W')$-submodule of $V'$ for every $j$.

Let $\gamma \in W$ be a~reflection. Then  $\gamma \in W_p$ for some $p$, see a) above. But $\psi \gamma \psi ^{-1}=\gamma '$ is also a~reflection (inside $W'$). Therefore, $\gamma '\in W_q'$ for some $q$.

The root line of $\Lin \gamma '$ is contained in $V_q'$, and this line is $(\Lin \psi )\ell$, where $\ell$ is the~root line of $\Lin \gamma $. Clearly, $\ell \subseteq V_p$. Hence,
$$
V_q' \cap (\Lin \psi )V_p \ne 0.
$$
It follows now from the irreducibility that, in fact,
$$
V_q'=(\Lin \gamma )V_p.
$$
Therefore,
$$
\psi W_p \psi ^{-1}\subseteq W_q'
$$
(because $W_p$ is generated by all reflections $\tau $ for which $\Lin \tau $ has its root line in $V_p$). The same proof is valid for the inverse inclusion. Hence
$$
\psi W_p \psi ^{-1}=W_q'
$$
and we can proceed by induction.

d) This is clear (see the equality an a)).\end{proof}

Therefore, \textit{from now on, we consider only the case of irreducible  $r$-gro\-ups}.



\subsection{What is already known: the case \boldmath $\ka=\Ree $}\label{1.5}

Let $W$ be an irreducible $r$-group in $A(E)$. There are two possibilities: either $W$ is \textit{finite}, or $W$ is \textit{infinite}.

\textit{If $W$ is finite, then there exists a~point in $E$ which is fixed under $W$}. Indeed, let $a\in E$ be an arbitrary point. Then
$$
b=a+\frac{1}{|W|}\sum\limits _{\gamma \in W}(\gamma (a)-a)
$$
is fixed under $W$. Therefore, $\kappa _b$ provides us with an isomorphism of $W$ with $\Lin W$, see Section \ref{1.1}, i.e., $W$ \textit{is a~linear group generated by reflections} (we identify $E$ and $V$ by choosing the point $b$ as the origin).

\vskip 2mm

\textbf{The case \boldmath $\ka=\Ree $.}

A beautiful classical theory concerning this case was developed by Coxeter, Witt, Stiefel (see [1]),
the results of which we recapitulate
below.

\vskip 2mm

1) \textbf{\boldmath $W$ is finite.}

In this case, $W$ is either the~Weyl group of an irreducible root system, or a~dihedral group, or one of two exceptional groups: ${H}_3$ or ${H}_4$.

The description of all these groups is usually given by means of their \textit{Coxeter graphs}. This is done in the following way. It is known that $W$ is generated by the elements $R_j$, where $1 \leqslant j \leqslant n$, which are reflections in the faces of a~Weyl chamber $C$. There exists a~unique set of vectors $e_j$, where $1 \leqslant j \leqslant n$, of unit length with the property: $R_j=R_{e_j,\theta _j}$ (where, in fact, $\theta _j=-1$), and
$$
C=\{v\in V \mid \langle e_j \mid v \rangle >0 \quad \text{for}\quad 1 \leqslant j \leqslant n\}.
$$
The angle between the mirrors $H_{R_i}$ and $H_{R_j}$ is of the form
$$
\frac{\pi }{m_{ij}},\quad \text{where}\quad m_{ij} \in \Zee  \quad \text{and}\quad m_{ij}\geqslant 2.
$$
The nodes of the Coxeter graph of $W$ are in bijective correspondence with the reflections $R_j$, where $1 \leqslant j \leqslant n$. Two nodes $R_i$ and $R_j$ are connected by an edge if and only if $m_{ij}\geqslant 3$. The weight of this edge is equal to $m_{ij}$ (if $m_{ij}=3$, then the weight is usually omitted). The complete list of
Coxeter graphs of finite irreducible $r$-groups is given by the following table.
$$
\left.\begin{array}{llll}
\unitlength 1mm 
\linethickness{0.4pt}
\ifx\plotpoint\undefined\newsavebox{\plotpoint}\fi 
\begin{picture}(42.57,0)(0,26)
\put(22.494,25.647){\makebox(0,0)[cc]{$\ldots$}}
\put(3.989,26.962){\circle{2.828}}
\put(5.356,26.962){\line(1,0){7.253}}
\put(13.97,26.904){\circle{2.828}}
\put(15.341,26.968){\line(1,0){3.679}}
\put(41.156,27.013){\circle{2.828}}
\put(39.789,27.013){\line(-1,0){7.253}}
\put(31.175,27.071){\circle{2.828}}
\put(29.804,27.007){\line(-1,0){3.679}}
\end{picture}%
&{\sf A}_n&
\unitlength 1mm 
\linethickness{0.4pt}
\ifx\plotpoint\undefined\newsavebox{\plotpoint}\fi 
\begin{picture}(35.422,0)(0,26)
\put(3.989,26.962){\circle{2.828}}
\put(5.356,26.962){\line(1,0){7.253}}
\put(13.97,26.904){\circle{2.828}}
\put(34.008,27.013){\circle{2.828}}
\put(32.641,27.013){\line(-1,0){7.253}}
\put(24.013,26.956){\circle{2.828}}
\put(15.399,26.961){\line(1,0){7.253}}
\put(19.078,29.274){\makebox(0,0)[cc]{$4$}}
\end{picture}%
&{\sf F}_4\\[4pt]
\unitlength 1mm 
\linethickness{0.4pt}
\ifx\plotpoint\undefined\newsavebox{\plotpoint}\fi 
\begin{picture}(42.57,0)(0,26)
\put(22.494,25.647){\makebox(0,0)[cc]{$\ldots$}}
\put(3.989,26.962){\circle{2.828}}
\put(5.356,26.962){\line(1,0){7.253}}
\put(13.97,26.904){\circle{2.828}}
\put(15.341,26.968){\line(1,0){3.679}}
\put(41.156,27.013){\circle{2.828}}
\put(39.789,27.013){\line(-1,0){7.253}}
\put(31.175,27.071){\circle{2.828}}
\put(29.804,27.007){\line(-1,0){3.679}}
\put(35.948,29.431){\makebox(0,0)[cc]{$4$}}
\end{picture}%
&{\sf B}_n={\sf C}_n&
\unitlength 1mm 
\linethickness{0.4pt}
\ifx\plotpoint\undefined\newsavebox{\plotpoint}\fi 
\begin{picture}(25.427,0)(10,26)
\put(13.97,26.904){\circle{2.828}}
\put(24.013,26.956){\circle{2.828}}
\put(15.399,26.961){\line(1,0){7.253}}
\put(19.078,29.274){\makebox(0,0)[cc]{$6$}}
\end{picture}%
&{\sf G}_2\\[8pt]
\unitlength 1mm 
\linethickness{0.4pt}
\ifx\plotpoint\undefined\newsavebox{\plotpoint}\fi 
\begin{picture}(42.623,0)(0,26)
\put(22.494,25.647){\makebox(0,0)[cc]{$\ldots$}}
\put(3.989,26.962){\circle{2.828}}
\put(5.356,26.962){\line(1,0){7.253}}
\put(13.97,26.904){\circle{2.828}}
\put(15.341,26.968){\line(1,0){3.679}}
\put(41.209,29.798){\circle{2.828}}
\put(31.175,27.071){\circle{2.828}}
\put(29.804,27.007){\line(-1,0){3.679}}
\put(41.199,24.223){\circle{2.828}}
\multiput(32.532,27.014)(.08801511,.03355972){83}{\line(1,0){.08801511}}
\multiput(32.637,26.909)(.08824182,-.03373952){81}{\line(1,0){.08824182}}
\end{picture}%
&{\sf D}_n&
\unitlength 1mm 
\linethickness{0.4pt}
\ifx\plotpoint\undefined\newsavebox{\plotpoint}\fi 
\begin{picture}(35.422,0)(10,26)
\put(13.97,26.904){\circle{2.828}}
\put(34.008,27.013){\circle{2.828}}
\put(32.641,27.013){\line(-1,0){7.253}}
\put(24.013,26.956){\circle{2.828}}
\put(15.399,26.961){\line(1,0){7.253}}
\put(19.078,29.274){\makebox(0,0)[cc]{$5$}}
\end{picture}%
&{\sf H}_3\\[8pt]
\unitlength 1mm 
\linethickness{0.4pt}
\ifx\plotpoint\undefined\newsavebox{\plotpoint}\fi 
\begin{picture}(45.351,0)(0,26)
\put(3.989,26.962){\circle{2.828}}
\put(5.356,26.962){\line(1,0){7.253}}
\put(13.97,26.904){\circle{2.828}}
\put(34.008,27.013){\circle{2.828}}
\put(32.641,27.013){\line(-1,0){7.253}}
\put(24.013,26.956){\circle{2.828}}
\put(15.399,26.961){\line(1,0){7.253}}
\put(43.937,27.009){\circle{2.828}}
\put(42.57,27.009){\line(-1,0){7.253}}
\put(24.018,25.647){\line(0,-1){7.253}}
\put(24.066,17.023){\circle{2.828}}
\end{picture}%
&{\sf E}_6&
\unitlength 1mm 
\linethickness{0.4pt}
\ifx\plotpoint\undefined\newsavebox{\plotpoint}\fi 
\begin{picture}(35.422,0)(0,26)
\put(3.989,26.962){\circle{2.828}}
\put(5.356,26.962){\line(1,0){7.253}}
\put(13.97,26.904){\circle{2.828}}
\put(34.008,27.013){\circle{2.828}}
\put(32.641,27.013){\line(-1,0){7.253}}
\put(24.013,26.956){\circle{2.828}}
\put(15.399,26.961){\line(1,0){7.253}}
\put(9.25,29.274){\makebox(0,0)[cc]{$5$}}
\end{picture}%
&{\sf H}_4\\[26pt]
\unitlength 1mm 
\linethickness{0.4pt}
\ifx\plotpoint\undefined\newsavebox{\plotpoint}\fi 
\begin{picture}(55.337,0)(0,26)
\put(3.989,26.962){\circle{2.828}}
\put(5.356,26.962){\line(1,0){7.253}}
\put(13.97,26.904){\circle{2.828}}
\put(34.008,27.013){\circle{2.828}}
\put(32.641,27.013){\line(-1,0){7.253}}
\put(24.013,26.956){\circle{2.828}}
\put(15.399,26.961){\line(1,0){7.253}}
\put(43.937,27.009){\circle{2.828}}
\put(42.57,27.009){\line(-1,0){7.253}}
\put(24.018,25.647){\line(0,-1){7.253}}
\put(24.066,17.023){\circle{2.828}}
\put(53.923,27.009){\circle{2.828}}
\put(52.556,27.009){\line(-1,0){7.253}}
\end{picture}%
&{\sf E}_7&
\unitlength 1mm 
\linethickness{0.4pt}
\ifx\plotpoint\undefined\newsavebox{\plotpoint}\fi 
\begin{picture}(20.427,0)(10,26)
\put(13.97,26.904){\circle{2.828}}
\put(24.013,26.956){\circle{2.828}}
\put(15.399,26.961){\line(1,0){7.253}}
\put(19.078,29.274){\makebox(0,0)[cc]{$p$}}
\end{picture}%
&\hspace*{-50pt}I_2(p) \text{ for } p=5 \text{ or } p \geqslant 7\\[26pt]
\unitlength 1mm 
\linethickness{0.4pt}
\ifx\plotpoint\undefined\newsavebox{\plotpoint}\fi 
\begin{picture}(45.351,0)(0,26)
\put(3.989,26.962){\circle{2.828}}
\put(5.356,26.962){\line(1,0){7.253}}
\put(13.97,26.904){\circle{2.828}}
\put(34.008,27.013){\circle{2.828}}
\put(32.641,27.013){\line(-1,0){7.253}}
\put(24.013,26.956){\circle{2.828}}
\put(15.399,26.961){\line(1,0){7.253}}
\put(43.937,27.009){\circle{2.828}}
\put(42.57,27.009){\line(-1,0){7.253}}
\put(24.066,17.023){\circle{2.828}}
\put(24.071,25.49){\line(0,-1){7.883}}
\put(53.923,27.009){\circle{2.828}}
\put(52.556,27.009){\line(-1,0){7.253}}
\put(63.923,27.009){\circle{2.828}}
\put(62.556,27.009){\line(-1,0){7.253}}
\end{picture}%
&&{\hskip -8mm \sf E}_8\\[26pt]
\end{array}
\right\}
\eqno{\boxed{1}}
$$

Set
$$
c_{ij}:=(1-\theta _i)(1-\theta _j)\langle e_i \mid e_j \rangle \langle e_j \mid e_i \rangle =
4 {\mbox{\rm cos}}^2 \frac \pi {m_{ij}}.
$$
One can change the weight $m_{ij}$ to the number $c_{ij}$ for all of the edges of the Coxeter graph. In this manner another weighted graph results.
Clearly, one graph determines the other.
We shall show later that the newly obtained graphs can be generalized to the case of $\ka=\mathbb C$.
\vskip 2mm


2) \textbf{\boldmath$W$ is infinite}.

In this case, $W$ is the affine Weyl group of an irreducible root system.

This group is \textit{a semidirect product} of $\Lin W$, which is the~(finite) Weyl group of a~certain root system $R$, and the
lattice of rank $n$ generated by the dual root system $\check R$. The groups $\Lin W$ thus obtained are distinguished from the others
in the above list $\boxed{1}$
in the following way
(see \cite[Chap.\,VI]{1}):

\subsubsection{Theorem {\rm (Linear parts of infinite  real irreducible {\it r}-groups)}}\label{lpre}
 {\it
Let $K \subseteq \GL (V)$ be a finite irreducible real $r$-group. Then the following properties are equivalent:
\begin{enumerate}[\hskip 7.2mm\rm a)]\itemsep=-.1ex
\item
$K=\Lin W$, where $W$ is an  infinite real irreducible $r$-group.
\item
There exists a~$K$-invariant lattice in $V$ of rank $n$.
\item
$K$ is defined over $\Qee $.
\item
$K$ is the Weyl group of a~certain irreducible root system, i.e., a~group whose Coxeter graph is one of ${\sf A}_n$, ${\sf B}_n$, ${\sf D}_n$, ${\sf E}_6$, ${\sf E}_7$, ${\sf E}_8$, ${\sf F}_4$, ${\sf G}_2$.
\item
All the numbers $c_{ij}$ lie in $\Zee  $.
\item
The ring with unity, generated over $\Zee  $ by all of the numbers $c_{ij}$, coincides with $\Zee  $.
\end{enumerate}
}


\subsection{What is already known: the case $\ka=\Cee $}\label{1.6}

Let $\ka=\Cee $ and let $W \subseteq A(E)$ be an irreducible $r$-group.

\vskip 2mm

1) \textbf{\boldmath $W$ is finite}.

Shephard and Todd, see [2], gave the complete list of such groups. A modern and unified approach was presented by Cohen, see \cite{3}.

We describe this classification in a~form that is more convenient for us, i.e., by means of certain graphs (as it was done in the real case).

Let $R_j=R_{e_j,\theta _j}$, where $1 \leqslant j \leqslant s$, be a~generating system of reflections of $W$. We can assume that
$$
\theta _j=e^{2 \pi i/m_j}.
$$
Therefore, this system (and hence, $W$) is \textit{uniquely defined by the system of lines $\ell_{R_j}$ in $V$ with the multiplicities} $m_j$ for $1 \leqslant j \leqslant
s$.

It is well known that an arbitrary set of \textit{vectors} in $V$ is uniquely (up to isometry) defined by means of a~certain set of numbers (more precisely, by the corresponding Gram matrix). Let us show that the same is true for an arbitrary set of lines in $V$ (i.e.,
points of the corresponding projective space).

\subsubsection{Proposition {\rm (Isometric systems of lines)}}\label{cyp}
{\it
Let $\{\ell_j\}_{j\in J}$ and $\{\ell_j'\}_{j\in J}$ be two sets of lines in $V$, and let $e_j\in \ell_j$ and $e'_j\in \ell'_j$ be arbitrary vectors with $1=\langle e_j \mid e_j \rangle =\langle e'_j \mid e'_j \rangle $. For 
every finite set of indices $j_1,\ldots ,j_d\in J$, consider the numbers
$$
h_{j_1 \ldots j_d}:=\langle e_{j_1}\mid e_{j_2}\rangle \langle e_{j_2}\mid e_{j_3}\rangle \cdots \langle e_{j_{d-1}}\mid e_{j_d}\rangle \langle e_{j_d}\mid e_{j_1}\rangle
$$
and
$$
h'_{j_1 \ldots j_d}:=\langle e'_{j_1}\mid e'_{j_2}\rangle \langle e'_{j_2}\mid e'_{j_3}\rangle \cdots \langle e'_{j_{d-1}}\mid e'_{j_d}\rangle \langle e'_{j_d}\mid e'_{j_1}\rangle .
$$
Then $h_{j_1 \ldots j_d}$ \textup(resp.\;$h'_{j_1 \ldots j_d}$\textup) is independent of the choice of the vectors $e_j$ \textup(resp.\;$e'_j$\textup) for $j\in J$. Moreover, the systems $\{\ell_j\}_{j\in J}$ and $\{\ell'_j\}_{j\in J}$ are isometric \textup(i.e., $g\ell_j=\ell'_j$ for each $j\in J$ and a~certain $g\in {\rm Iso}(V)$\textup) if and only if
$$
h_{j_1 \ldots j_d}=h'_{j_1 \ldots j_d}
\eqno\boxed{2}
$$
for each $j_1,\ldots ,j_d\in J$.
}
\vskip 2mm

\noindent{\it Proof.} We need only prove that if $\,\boxed{2}\,$
is fulfilled, then for every $i\in J$, there exists a~number $\lambda _i\in \Cee $ such that
$$
\langle e_j' \mid e_l' \rangle =\langle \lambda _je_j \mid \lambda _le_l \rangle \quad \text{ for every }j,l\in J,
$$
i.e., the Gram matrices in bases $\{e_j\}_{j\in J}$ and $\{\lambda _je_j\}_{j\in J}$ are the same (all the other statements are evident).

Let us fix an index $t\in J$. It is not difficult to see that one can assume that the system $\{\ell_j\}_{j\in J}$ is \lq \lq connected\rq \rq , i.e., for every $i\in J$ there exists a~sequence $j_1,\ldots ,j_d\in J$ such that
$$
j_1=i,\quad j_d=t \quad \text{and}\quad \langle e_{j_l}\mid e_{j_{l+1}}\rangle \ne 0 \quad \text{for each}\quad l=1,\ldots ,d-1.
$$
 Now a straightforward computation shows
 that one can take
$$
\lambda _j:=\prod\limits _{l=1}^{d-1}\frac {\langle e_{j_{l}}' \mid e_{j_{l+1}}' \rangle }{\langle e_{j_{l}}\mid e_{j_{l+1}}\rangle }.
\eqno \text{\textsquare}
$$

\vskip 2mm

We have seen above that every $r$-group is defined by a~system $\{\ell_j,m_j\}_{j\in J}$ of lines $\ell_j\in V$ with \textit{multiplicities} $m_j\in \Zee  $. It follows from the proposition that \textit{such a~system is uniquely} (\textit{up to isometry}) \textit{defined by the~system of numbers}
$$
c_{j_1 \ldots j_d}:=h_{j_1 \ldots j_d}\prod\limits _{l=1}^d(1-e^{2 \pi i/m_{j_l}})
$$
(one can derive from these numbers the multiplicities because $c_j=1-e^{2 \pi i/m_j}$). It will become clear later on why $h_{j_1 \ldots j_d}$ is multiplied by $\prod
_{l=1}^d(1-e^{2 \pi i/m_{j_l}}$) and not, say, by $\prod
 _{l=1}^d(e^{2 \pi i/m_{j_l}})$; the numbers $c_{j_1 \ldots j_d}$ are of great importance in the whole theory. We call them \textit{cyclic products}\index{cyclic products}.

So the group $W$ (with a~fixed generating system of reflections) is uniquely (up to
equivalence) defined by the corresponding set of cyclic products. As a~matter of fact one only needs to know the so-called \textit{simple cyclic products} $c_{j_1 \ldots j_m}$, i.e., those with all indices $j_1,\ldots ,j_m$ distinct, because
$$
c_{l _1 \ldots l _{p-1}l l _{p+1}\ldots l _{q-1}l l _{q+1}\ldots \ell _r}=c_{l _1 \ldots  l_{p-1}l l _{q+1}\ldots l _r}\cdot c_{l l _{p+1}\ldots l _{q-1}}.
$$
\vskip\abovedisplayskip
\begin{center}
\includegraphics[width=70mm]{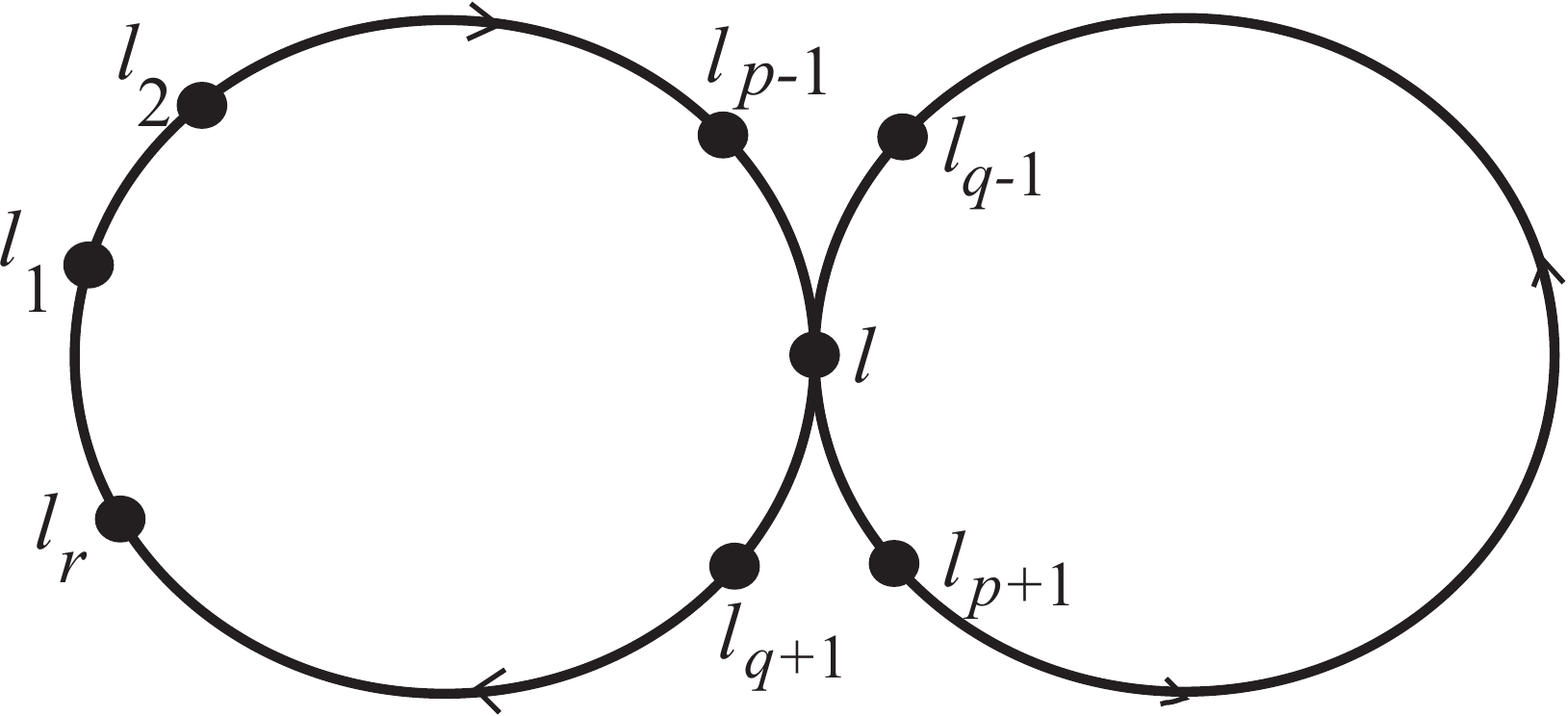}
\end{center}
\vskip\belowdisplayskip

We want to specify here \textit{several properties of the cyclic products}.

\vskip 2mm

a) If $j_1' \ldots j_d'$ is \textit{a cyclic permutation} of $j_1,\ldots ,j_d$, then
$$
c_{j_1 \ldots j_d}=c_{j_1' \ldots j_d'}.
$$
In other words, $c_{j_1 \ldots j_d}$ depends only on the~cycle $\sigma =(j_1,\ldots ,j_d)$. We use therefore \textit{the notation}
$$
c_\sigma :=c_{j_1 \ldots j_d}.
$$
In particular,
$$
c_{jl}=c_{lj}.
$$

b) If $\sigma =(j_1,\ldots ,j_d)$, then
$$
c_\sigma c_{\sigma^{-1}}=c_{j_1j_2}c_{j_2j_3}\cdots c_{j_{d-1}j_d}c_{j_dj_1}.
$$

c) One can reconstruct all the simple cyclic products (hence, all the cyclic products) only from the \lq \lq homologically independent\rq \rq {} ones. The following formula and drawing illustrates what we have in mind:
$$
c_{l _1 l _2 \ldots l _pj_dj_{d-1}\ldots j_1}\cdot c_{l _p l _{p+1}\ldots l _{q-1}l _q l_1j_1j_2 \ldots j_d}=c_{l _1 l _2 \ldots l _{p-1}l _p l _{p+1}\ldots l _{q-1}l _q}\cdot c_{l _1 j_1}\cdot c_{j_1j_2}\ldots c_{j_{d-1}j_d}\cdot c_{j_d l _p}.
$$
\vskip\abovedisplayskip
\begin{center}
\includegraphics[width=60mm]{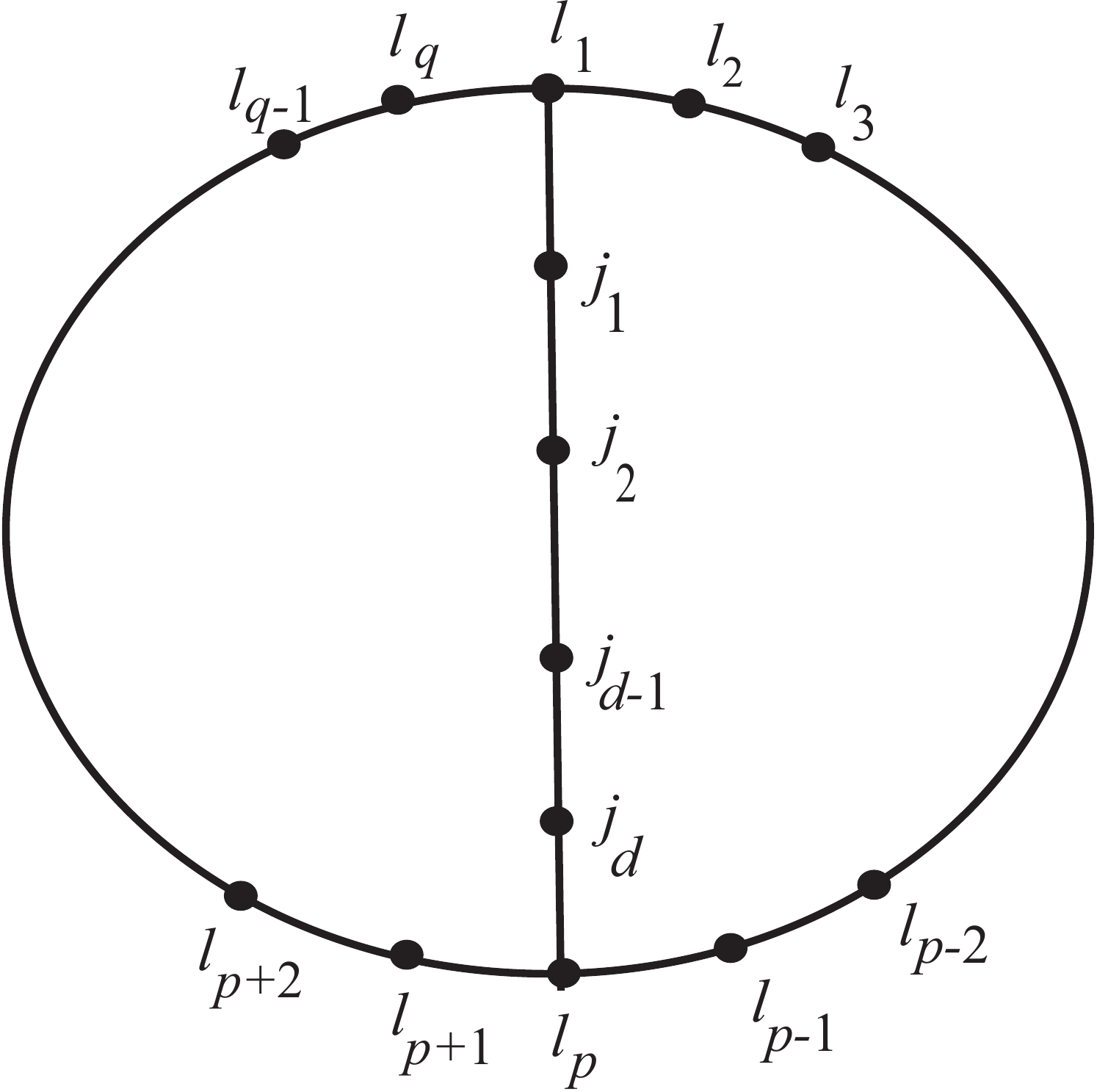}
\end{center}
\vskip\belowdisplayskip

A system of lines $\ell _j\in V$ with the multiplicities $m_j$ for $j\in J$, can be described by \textit{a graph} as follows.

The nodes of the graph are in bijective correspondence with the lines $\ell _j$ for $j\in J$. If a~node represents the line $\ell_j$, then \textit{this node has the weight} $m_j$. Two nodes $\ell _i$ and $\ell _j$ are connected by an edge if and only if $c_{ij}\ne 0$, and if they are connected, then \textit{the weight of the edge} is equal to $c_{ij}$.

\textit{It is convenient not to specify the weight of a node, respectively an edge, if it is equal to
$2$, respectively $1$. Below we follow this convention}. Moreover, every simple cycle of this graph is supplied with an arbitrary (but fixed) orientation and has weight equal to the corresponding cycle product:

\vskip\abovedisplayskip
\begin{center}
\includegraphics[width=53mm]{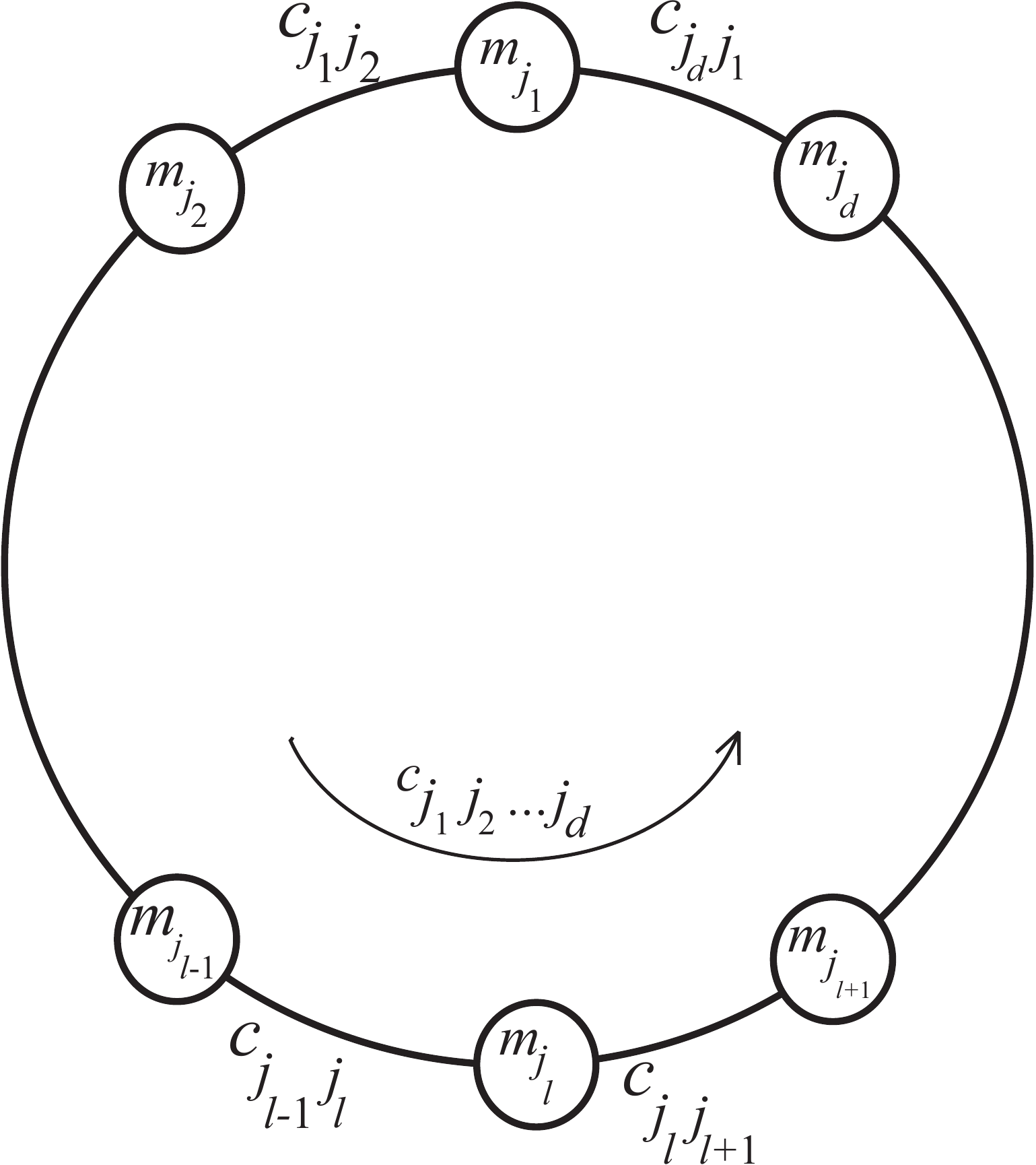}
\end{center}
\vskip -2mm

Therefore, we now have a~way \textit{to represent a~finite $r$-group $W \subseteq \GL (V)$ with a~fixed generating system of reflections by means of a~graph corresponding to the system of lines $\ell _{R_j}\subseteq V$ with the multiplicities $m_j$ for $1 \leqslant j \leqslant h$}. (Of course, using another system of generators one obtain \textit{another} graph which represents \textit{the same} group. This nonuniqueness in the representations of the group by means of its graph occurs because, contrary to the real case, there is no known canonical method for constructing a~generating system of reflections of a finite
complex
$r$-group. The problem\index{Open problem} of finding such a~method is still unsolved and seems to be very interesting.) It is easy to see that, $W$ being irreducible, the graph is \textit{connected}. This graph is called \textit{the graph of the group $W$}\index{group} (with respect to a~fixed generating system of reflections).

The classification of the finite complex irreducible $r$-groups $W$ was given in \cite{2,3} \textit{by means of generating systems of reflections}. It is now a~matter of more or less straightforward computation to reformulate the result \textit{by means of the graphs}. We need the following
notation to formulate the corresponding theorem:
$$
\omega =e^{2 \pi i/3},\quad \eta =e^{2 \pi i/5},\quad \varepsilon =e^{2 \pi i/8},\quad \zeta_m =e^{2 \pi i/m}.
$$

It appears \textit{a posteriori} that all the graphs under consideration are planar; they either have no simple cycles, or have only one such cycle (of length 3). We assume that this cycle is \textit{counter-clockwise oriented}.

We also \textit{fix a~numbering of the nodes} of the graph (in an arbitrary fashion). The number of the node is written \textit{beside} the node (but the weight of the node is written \textit{inside}).

In the table below, the ring with unity generated over $\Zee  $ by all cyclic products is also given. We need this ring later on; it plays an important role in the theory and does not depend on the choice or the generating system of reflections (and hence, on the graph that represents the group).

\subsubsection{Theorem {\rm (Classification of finite complex irreducible {\it r}-groups)}}\label{cfi}
{\it Up to equivalence, finite complex irredu\-cible
$r$-groups are precisely those
determined \textup(with respect to the~fixed generating systems of reflections\textup)
by the~graphs from Table\;{\rm 1}
below
\textup{(}our numbering of the groups coincides with that of
Shephard and Todd {\rm\cite{2}}; the notation of types is as in {\rm\cite{2,3,4}}{\textup)}.
}

\newpage

\begin{landscape}\addcontentsline{toc}{subsubsection}{Table 1. The finite complex irreducible
$r$-groups}
{\centering \footnotesize
\begin{longtable}{|r|c|c|c|c|}
\multicolumn{5}{l}{\textbf{Table 1.} The finite complex irreducible
$r$-groups}\\
\hline
No&Type&Graph&$\begin{array}{@{}c@{}}
\text{Ring generated}\\
\text{by cyclic products}
\end{array}
$&$\dim V$\\
\hline
\endfirsthead
\hline
{\rm No}&Type&Graph&$\begin{array}{@{}c@{}}
\text{Ring generated}\\
\text{by cyclic products}
\end{array}
$&$\dim V$\\
\hline
\endhead
1&${\sf A}_s$, $s \geqslant 1$&$\begin{array}{@{}c@{}}
{}\\[-12pt]
\text{\includegraphics{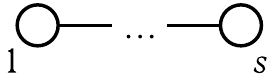}}
\end{array}
$&$\Zee  $&$s$\\
\hline
2&$\begin{array}{@{}c@{}}
G(m,1,s)\\
m \geqslant 2,\ s \geqslant 2,\\
\text{type }{\sf B}_s={\sf C}_s\\
\text{if }m=2
\end{array}
$&$\begin{array}{@{}ll@{}}
{}\\[-12pt]
\text{\includegraphics{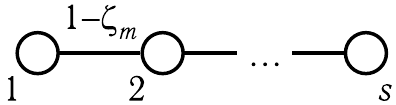}}\\
\text{\includegraphics{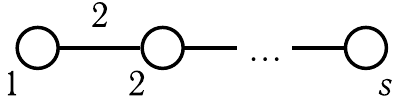}}&\begin{array}[b]{@{}c@{}}
\text{if} \\[-3pt]
m=2
\end{array}
\end{array}
$&$\begin{array}{@{}l@{}}
\Zee  [e^{2 \pi i/m}], \\
\Zee  [\omega ]\text{ if }m=3, 6,\\
\Zee  [i] \text{ if }m=2, \\
\Zee  \text{ if }m=2
\end{array}
$&$s$\\
\hline
2&$\begin{array}{@{}l@{}}
G(m,m,s) \\
m \geqslant 2,\ s \geqslant 3, \\
\text{type }{\sf D}_s\\
\text{if }m=2
\end{array}
$&$\begin{array}{@{}r@{}}
{}\\[-12pt]
\text{\includegraphics{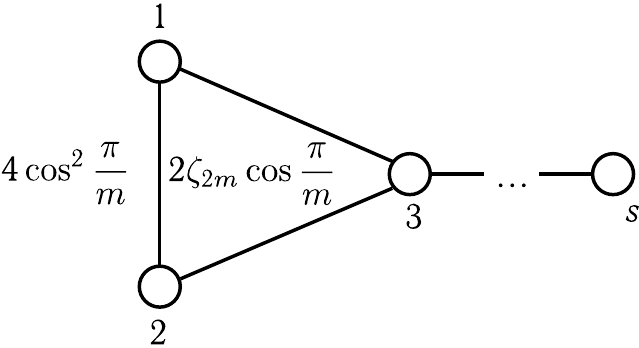}}\\
\text{\includegraphics{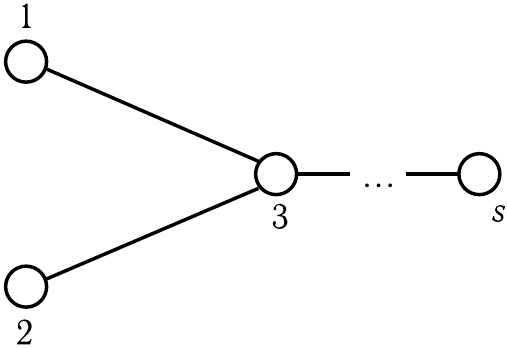}} \\[-36pt]
\text{if }m=2\\[16pt]
\end{array}
$&$\begin{array}{@{}l@{}}
\Zee  [e^{2 \pi i/m}], \\
\Zee  [\omega ]\text{ if }m=3,6, \\
\Zee  [i] \text{ if }m=4, \\[100pt]
\Zee  \text{ if }m=2
\end{array}
$&$s$\\
\hline
2&$\begin{array}{@{}l@{}}
G(m,m,2)={} \\
{}=I_2(m),\ m \geqslant 3, \\
{\sf A}_2 \text{ if } m=3, \\
{\sf B}_2 \text{ if } m=4, \\
{\sf G}_2 \text{ if } m=6
\end{array}
$&$\begin{array}{@{}ll@{}}
{}\\[-12pt]
\multicolumn{2}{@{}l@{}}{\text{\includegraphics{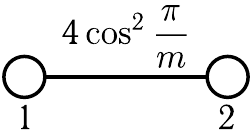}}} \\
\hspace*{-7pt}\text{\includegraphics{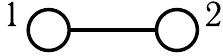}}&\begin{array}{@{}l@{}}
\text{if }m=3, \\[-5pt]
{}\end{array}
 \\
\hspace*{-7pt}\text{\includegraphics{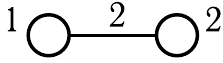}}&\begin{array}{@{}l@{}}
\text{if }m=4,\\[-5pt]
{}\end{array}
 \\
\hspace*{-5.5pt}\hskip .4mm\text{\includegraphics[width=23mm]{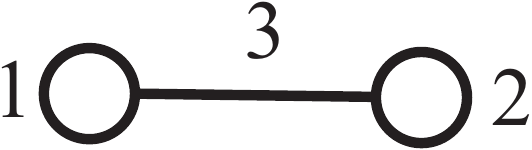}}&\begin{array}{@{}l@{}}
\text{if }m=6\\[-5pt]
{}\end{array}
\end{array}
$&$\begin{array}{@{}l@{}}
\Zee  \left [4 {\mbox{\rm cos}}^2 \frac{\pi}{m} \right ], \\
\Zee  \text{ if }m=3,4,6
\end{array}
$&2\\
\hline
2&$\begin{array}{@{}l@{}}
G(m,p,s-1) \\
m \geqslant 2, \\
s \geqslant 4 \\
p \mid m\\
p \ne 1,m
\end{array}
$&$\begin{array}{@{}l@{}}
{}\\[-12pt]
\text{\includegraphics{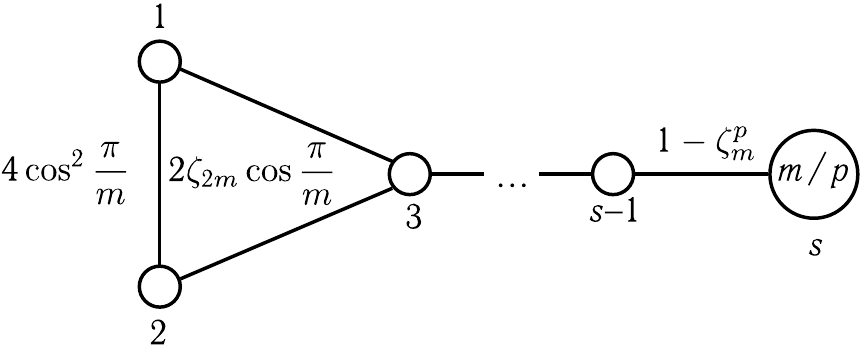}}
\end{array}
$&$\begin{array}{@{}l@{}}
\Zee  [e^{2 \pi i/m}], \\
\Zee  [i] \text{ if }m=4,p=2, \\
\Zee  [\omega ] \text{ if }m=6,p=2\\
\text{and }m=6,p=3
\end{array}
$&$s-1$\\
\hline
2&$\begin{array}{@{}l@{}}
G(m,p,2) \\
m \geqslant 2, \\
p \mid m, \\
p \ne 1,m
\end{array}
$&$\begin{array}{@{}l@{}}
{}\\[-12pt]
\text{\includegraphics[width=63mm]{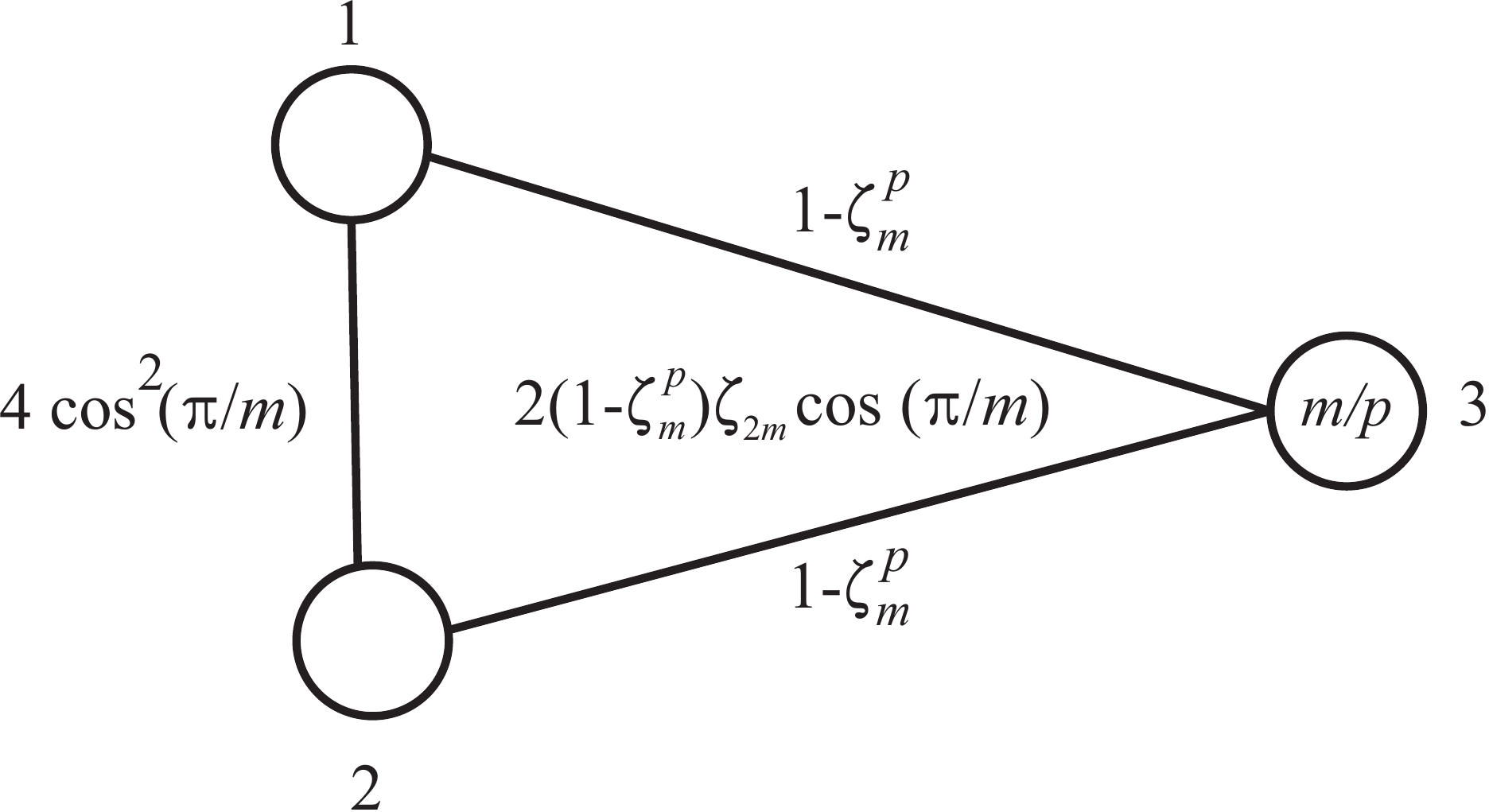}}
\end{array}
$&$\begin{array}{@{}l@{}}
\Zee  \left [\zeta _m,2 {\mbox{\rm cos}}\, \frac{\pi}{m},\zeta _m(1-\zeta _m^p)\right ], \\
\Zee  [2i]\text{ if }m=4,p=2, \\
\Zee  [\omega ] \text{ if }m=6,p=2, \\
\Zee  [2 \omega ]\text{ if }m=6,p=3
\end{array}
$&2\\
\hline
3&$[\ ]^m$&$\begin{array}{@{}l@{}}
{}\\[-12pt]
\text{\includegraphics{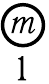}}
\end{array}
$&$\begin{array}{@{}l@{}}
\Zee  [e^{2 \pi i/m}], \\
\Zee  [\omega ]\text{ if }m=6,3, \\
\Zee  [i] \text{ if }m=4, \\
\Zee  \text{ if }m=2
\end{array}
$&1\\
\hline
4&$3[3]3$&$\begin{array}{@{}l@{}}
{}\\[-12pt]
\text{\includegraphics{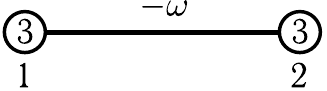}}
\end{array}
$&$\Zee  [\omega ]$&2\\
\hline
5&$3[4]3$&$\begin{array}{@{}l@{}}
{}\\[-12pt]
\text{\includegraphics{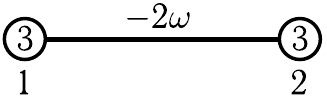}}
\end{array}
$&$\Zee  [\omega ]$&2\\
\hline
6&$3[6]2$&$\begin{array}{@{}l@{}}
{}\\[-12pt]
\text{\includegraphics{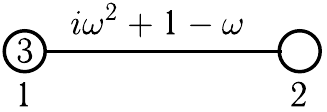}}
\end{array}
$&$\Zee  [e^{2 \pi i/12}]$&2\\
\hline
7&$\langle 3,3,2 \rangle _6$&$\begin{array}{@{}l@{}}
{}\\[-12pt]
\text{\includegraphics{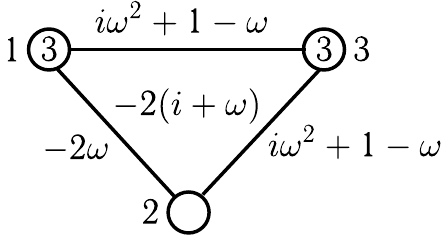}}
\end{array}$&$\Zee  [e^{2 \pi i/12}]$&2\\
\hline
8&$4[3]4$&$\begin{array}{@{}l@{}}
{}\\[-12pt]
\text{\includegraphics{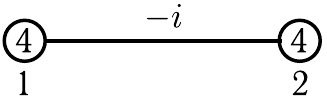}}
\end{array}
$&$\Zee  [i]$&2\\
\hline
9&$4[6]2$&$\begin{array}{@{}l@{}}
{}\\[-12pt]
\text{\includegraphics{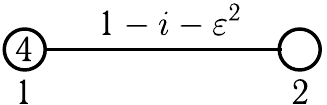}}
\end{array}
$&$\Zee  [\varepsilon]$&2\\
\hline
10&$4[4]3$&$\begin{array}{@{}l@{}}
{}\\[-12pt]
\text{\includegraphics[width=38mm]{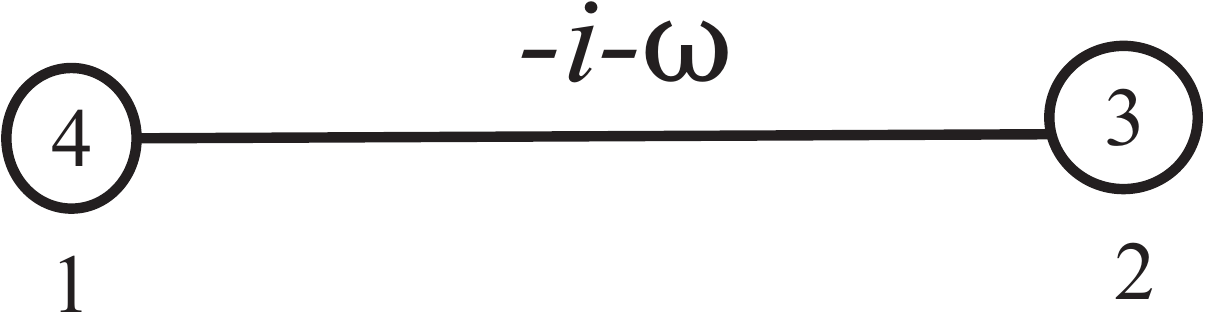}}
\end{array}
$&$\Zee  [e^{2 \pi i/12}]$&2\\
\hline
11&$\langle 4,3,2 \rangle _{12}$&$\begin{array}{@{}l@{}}
{}\\[-12pt]
\text{\includegraphics{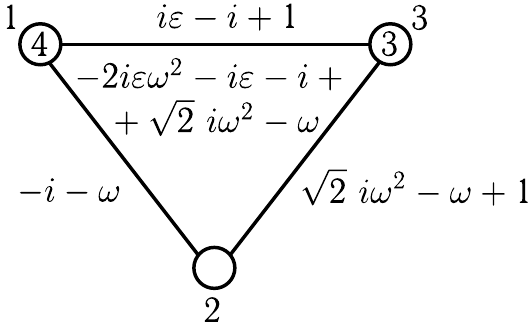}}
\end{array}
$&$\Zee  [e^{2 \pi i/24}]$&2\\
\hline
12&$\GL (2,3)$&$\begin{array}{@{}l@{}}
{}\\[-12pt]
\text{\includegraphics{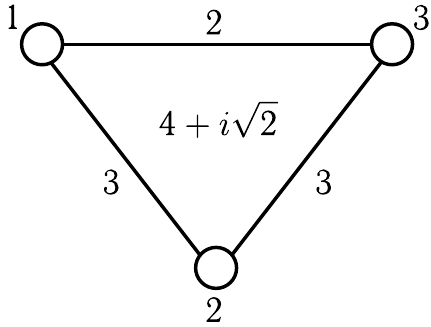}}
\end{array}
$&$\Zee  [i \sqrt 2]$&2\\
\hline
13&$\langle 4,3,2 \rangle _2$&$\begin{array}{@{}l@{}}
{}\\[-12pt]
\text{\includegraphics{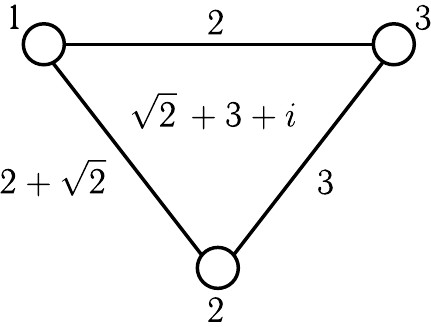}}
\end{array}
$&$\Zee  [i,\sqrt 2]$&2\\
\hline
14&$3[8]2$&$\begin{array}{@{}l@{}}
{}\\[-12pt]
\text{\includegraphics{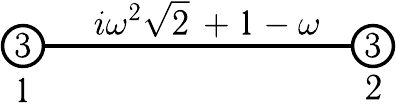}}
\end{array}
$&$\Zee  [\omega ,i \sqrt 2]$&2\\
\hline
15&$\langle 4,3,2 \rangle _6$&$\begin{array}{@{}l@{}}
{}\\[-12pt]
\text{\includegraphics{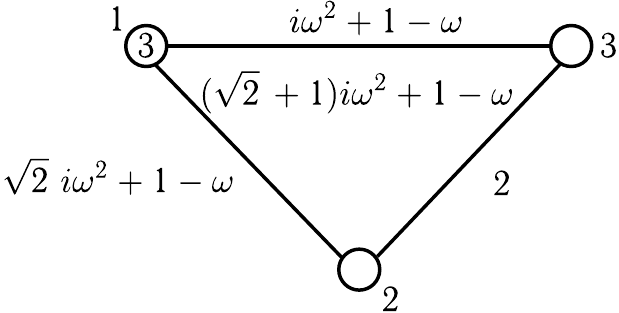}}
\end{array}
$&$\Zee  [i,\omega ,\sqrt 2]$&2\\
\hline
16&$5[3]5$&$\begin{array}{@{}l@{}}
{}\\[-12pt]
\text{\includegraphics[width=35mm]{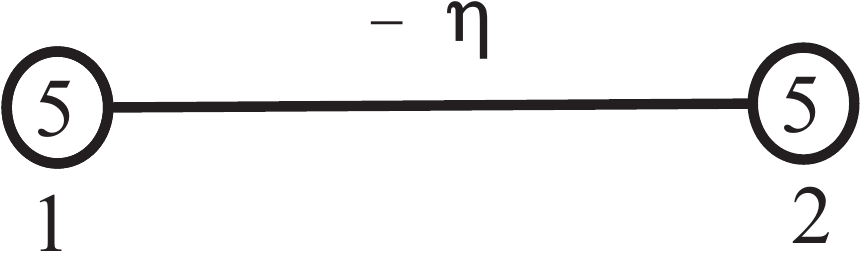}}
\end{array}
$&$\Zee  [\eta ]$&2\\
\hline
17&$5[6]2$&$\begin{array}{@{}l@{}}
{}\\[-12pt]
\text{\includegraphics[width=35mm]{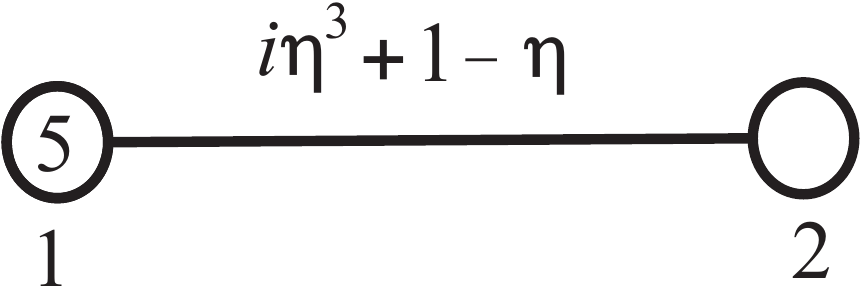}}
\end{array}
$&$\Zee  [e^{2 \pi i/20}]$&2\\
\hline
18&$5[4]3$&$\begin{array}{@{}l@{}}
{}\\[-12pt]
\text{\includegraphics[width=35mm]{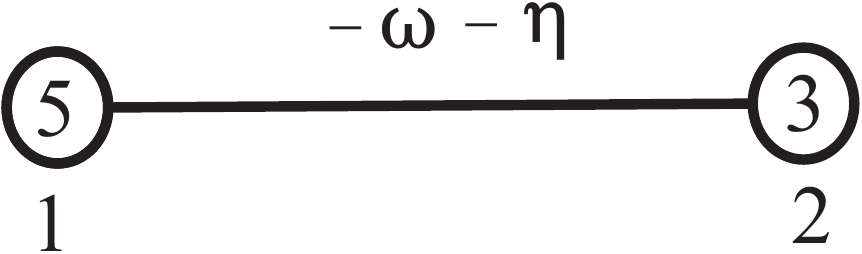}}
\end{array}
$&$\Zee  [e^{2 \pi i/15}]$&2\\
\hline
19&$\langle 5,3,2 \rangle _{30}$&$\begin{array}{@{}l@{}}
{}\\[-12pt]
\text{\includegraphics{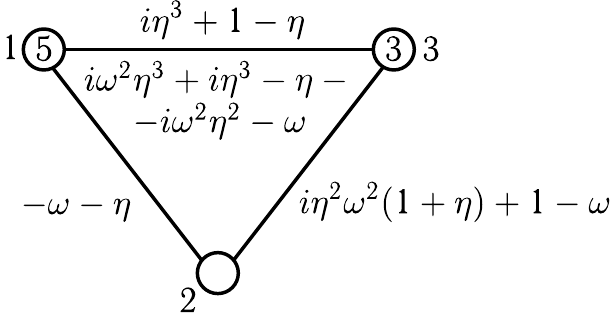}}
\end{array}
$&$\Zee  [e^{2 \pi i/60}]$&2\\
\hline
20&$3[5]3$&$\begin{array}{@{}l@{}}
{}\\[-12pt]
\text{\includegraphics{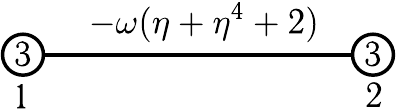}}
\end{array}
$&$\Zee  \left [\omega ,\frac{1+\sqrt 5}{2} \right ]$&2\\
\hline
21&$3[10]2$&$\begin{array}{@{}l@{}}
{}\\[-12pt]
\text{\includegraphics{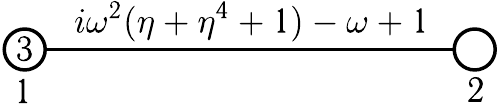}}
\end{array}
$&$\Zee  \left [\omega ,i\frac{1+\sqrt 5}{2} \right ]$&2\\
\hline
22&$\langle 5.3.2 \rangle _2$&$\begin{array}{@{}l@{}}
{}\\[-12pt]
\text{\includegraphics{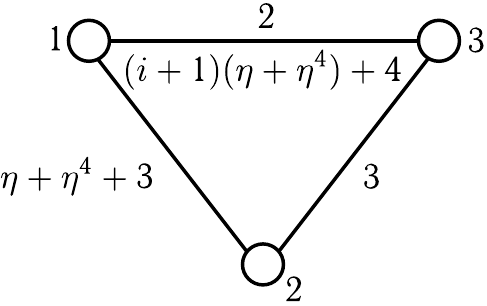}}
\end{array}
$&$\Zee  \left [\frac{\sqrt 5-1}{2},i \frac{\sqrt 5-1}{2} \right ]$&2\\
\hline
23&${\sf H}_3$&$\begin{array}{@{}l@{}}
{}\\[-12pt]
\text{\includegraphics[width=45mm]{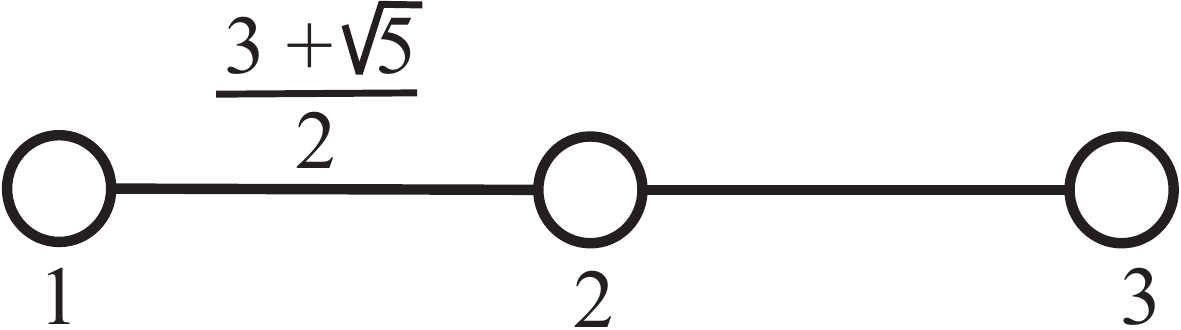}}
\end{array}
$&$\Zee  \left [\frac{1+\sqrt 5}{2} \right ]$&3\\
\hline
24&$J_3(4)$&$\begin{array}{@{}l@{}}
{}\\[-12pt]
\text{\includegraphics{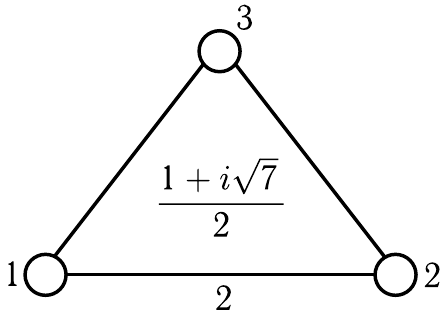}}
\end{array}
$&$\Zee  \left [\frac{1+i \sqrt 7}{2} \right ]$&3\\
\hline
25&$L_3$&$\begin{array}{@{}l@{}}
{}\\[-12pt]
\text{\includegraphics{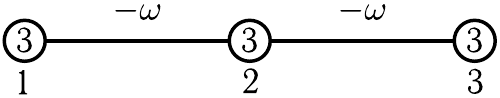}}
\end{array}
$&$\Zee  [\omega ]$&3\\
\hline
26&$M_3$&$\begin{array}{@{}l@{}}
{}\\[-12pt]
\text{\includegraphics{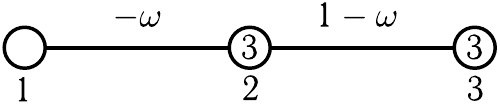}}
\end{array}
$&$\Zee  [\omega ]$&3\\
\hline
27&$J_3(5)$&$\begin{array}{@{}l@{}}
{}\\[-12pt]
\text{\includegraphics{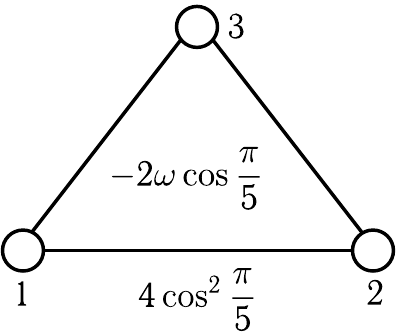}}
\end{array}
$&$\Zee  \left [\omega ,\frac{1+\sqrt 5}{2} \right ]$&3\\
\hline
28&${\sf F}_4$&$\begin{array}{@{}l@{}}
{}\\[-12pt]
\text{\includegraphics{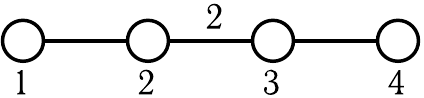}}
\end{array}
$&$\Zee  $&4\\
\hline
29&$\begin{array}{@{}l@{}}
[2\ 1;1]^4\\
=N_4
\end{array}
$&$\begin{array}{@{}l@{}}
{}\\[-12pt]
\text{\includegraphics{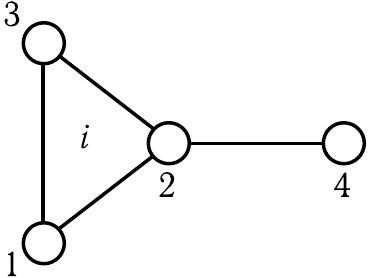}}
\end{array}
$&$\Zee  [i]$&4\\
\hline
30&$H_4$&$\begin{array}{@{}l@{}}
{}\\[-12pt]
\text{\includegraphics{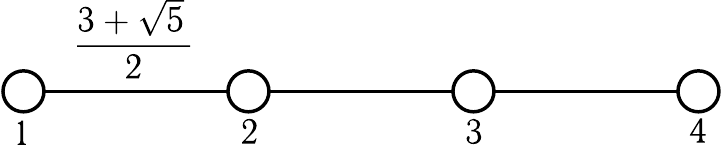}}
\end{array}
$&$\Zee  \left [\frac{1+\sqrt 5}{2} \right ]$&4\\
\hline
31&$\begin{array}{@{}l@{}}
\left [\left (\frac {1}{2}\,\gamma _3^4 \right )^{+1}\right ] \\
=EN_4
\end{array}
$&$\begin{array}{@{}l@{}}
{}\\[-12pt]
\text{\includegraphics[width=48mm]{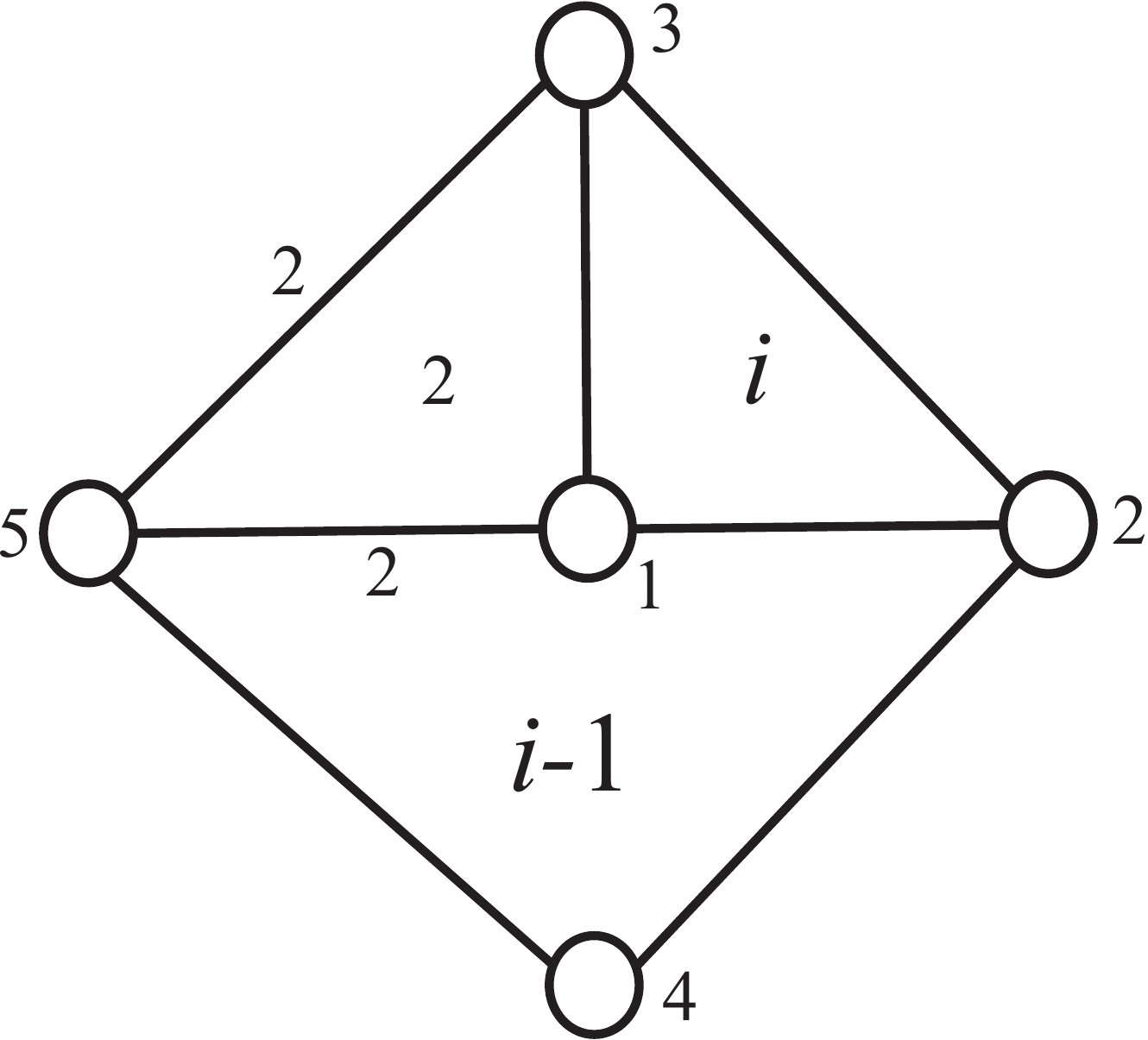}}
\end{array}
$&$\Zee  [i]$&4\\
\hline
32&$L_4$&$\begin{array}{@{}l@{}}
{}\\[-12pt]
\text{\includegraphics{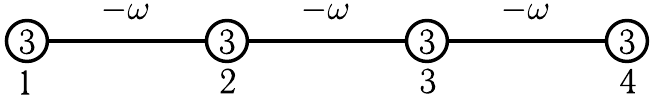}}
\end{array}
$&$\Zee  [\omega ]$&4\\
\hline
33&$\begin{array}{@{}l@{}}
[2\ 1;2]^3\\
=K_5
\end{array}
$&$\begin{array}{@{}l@{}}
{}\\[-12pt]
\text{\includegraphics{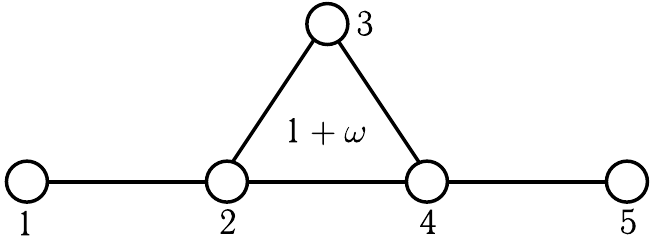}}
\end{array}
$&$\Zee  [\omega ]$&5\\
\hline
34&$\begin{array}{@{}l@{}}
[2\ 1;3]^3 \\
=K_6
\end{array}
$&$\begin{array}{@{}l@{}}
{}\\[-12pt]
\text{\includegraphics{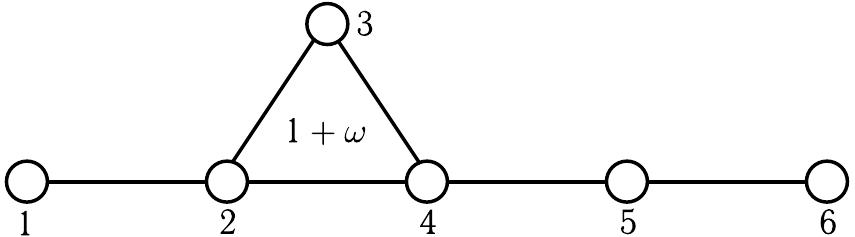}}
\end{array}
$&$\Zee  [\omega ]$&6\\
\hline
35&${\sf E}_6$&$\begin{array}{@{}l@{}}
{}\\[-12pt]
\text{\includegraphics{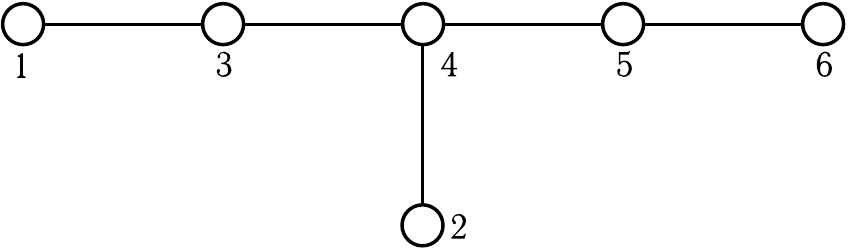}}
\end{array}
$&$\Zee  $&6\\
\hline
36&${\sf E}_7$&$\begin{array}{@{}l@{}}
{}\\[-12pt]
\text{\includegraphics{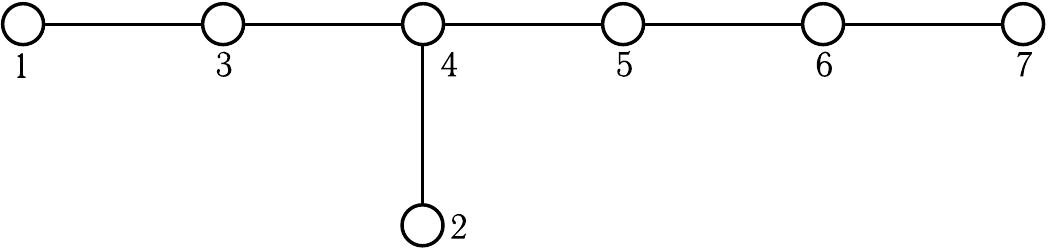}}
\end{array}
$&$\Zee  $&7\\
\hline
37&${\sf E}_8$&$\begin{array}{@{}l@{}}
{}\\[-12pt]
\text{\includegraphics{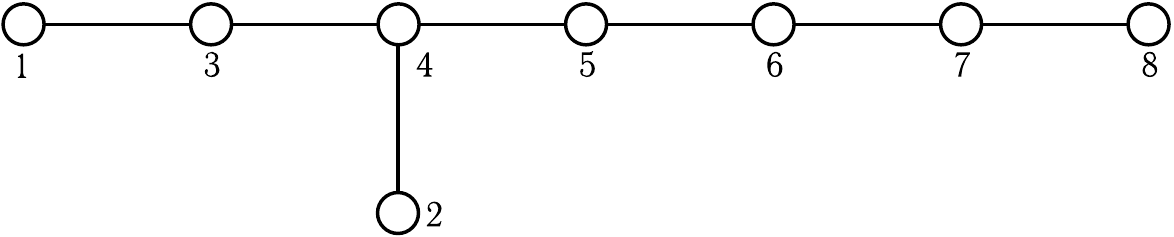}}
\end{array}
$&$\Zee  $&8\\
\hline         
\end{longtable}}
\end{landscape} 

\newpage

\subsubsection{Remark}
For
the groups from Table 1
that are of the form $\Lin W$, where $W$ is
an \textit{infinite} complex irreducible
$r$-group, one obtains from Table \ref{tab2} below an explicit description
of lines
$\ell _1,\ldots \ell_s$
(by means of specifying a~vector $e_j$ of unit length in each
$\ell _j$).\;For other groups, such explicit descriptions may be either
found in [3] or obtained directly from the graphs.

\vskip 2mm

2) \textbf{\boldmath $W$ is infinite}

This case was not investigated earlier and is our main concern in these lectures. The results are formulated in the next section. In his section, we only give several simple
 \textit{examples}, which show, first, that the groups under consideration do exist and, second, that we have here a~phenomenon which does not occur in the real case.

\subsection{Examples} We consider the case where $n=\dim E=1$. Let $a\in E$ be a~point. We identify $A(E)$ and $\GL (V)\ltimes V$ by means of $\kappa _a$, see Section  \ref{1.1}.

Let $\Gamma \ne 0$ be \textit{a lattice} in $V$.\;(\textrm{Hereafter a
lattice
means a~discrete subgroup of the additive group of a~vector space {\it not necessarily of maximal rank}}.) Let us consider a~subgroup
$$
W=\{(\pm 1,t)\mid \mathfrak t\in \Gamma \}
$$
of $\GL (V)\ltimes V$. This subgroup acts on $E$ (see Section  \ref{1.1}) and it is easy to check that it is an \textit{infinite complex irreducible
$r$-group}.

There are two possibilities: either $\rk \Gamma =1$ or $\rk \Gamma =2$.

If $\rk \Gamma =1$, i.e.,

 $$\Gamma=\mathbb Zv\;\;\mbox{where $v\in V, v\neq 0$}$$

\vskip\abovedisplayskip

\begin{center}
\includegraphics[width=
70mm]{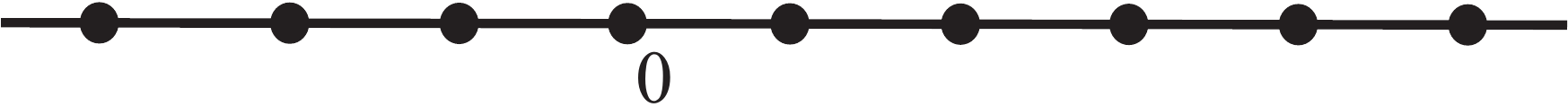}
\end{center}

\vskip\belowdisplayskip

\noindent
then $E/W$ is \textit{not compact}, i.e., \textit{by definition}, $W$ is \textit{a noncrystallographic group}. We denote this $W$ by $W(2, v)$. In this case, $W$ can be viewed as a~real $r$-group acting on $\Gamma \otimes _{\Zee  }\Ree $. As we know, it is impossible for an infinite real $r$-group to be noncrystallographic, see Section  \ref{1.5}.

If $\rk \Gamma =2$, i.e.,
$$
\Gamma=(\mathbb Z+\lambda\mathbb Z)v \;\;\mbox{where $v\in V, v\neq 0$, and $\lambda\in\mathbb C\setminus\mathbb R$}
$$

\vskip\abovedisplayskip

\begin{center}
\includegraphics[width=
70mm]{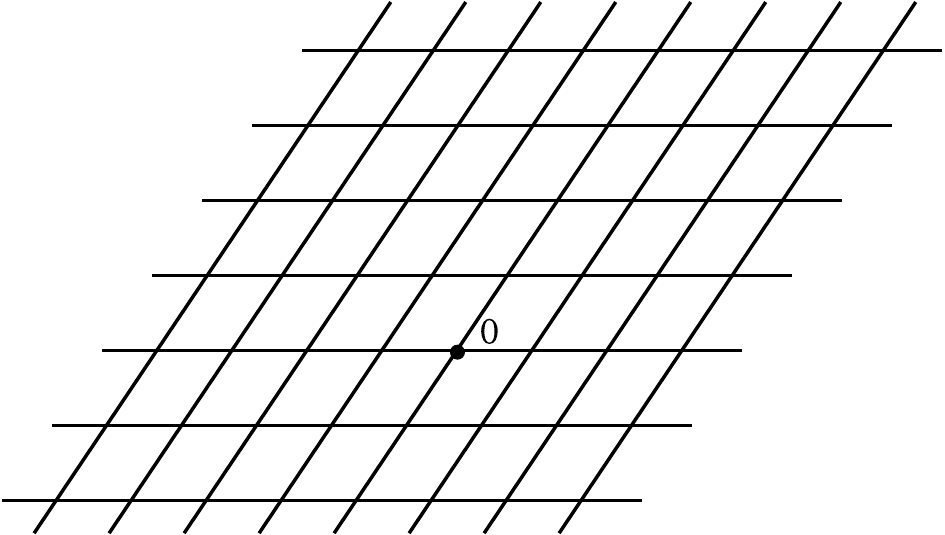}
\end{center}

\vskip\belowdisplayskip

\noindent
then $E/W$ is \textit{compact}, i.e.,
by definition,
$W$ is a~\textit{crystallographic group}. We denote this $W$ by $W(2,v,\lambda)$.



If
$\Gamma$ is a~lattice of equilateral triangles in $V$ ($=\Cee $), i.e.,

$$
\Gamma=
(\mathbb Z+\omega \mathbb Z)v\;\;\mbox{where $v\in V, v\neq 0$}
$$

\vskip\abovedisplayskip

\begin{center}
\includegraphics[width=40
mm]{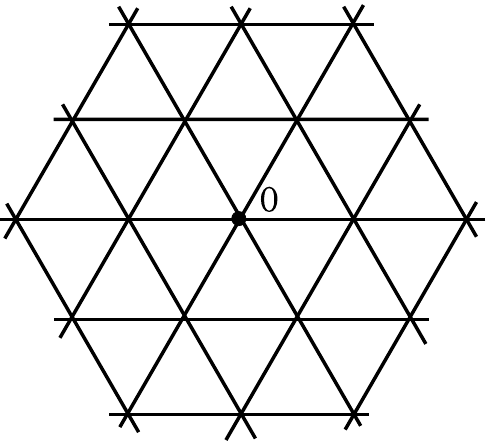}
\end{center}

\vskip\belowdisplayskip

\noindent
then  it is not difficult to check that
\begin{align*}
W(3, v)&:=\{(\omega ^l ,t)\mid t\in \Gamma,
l\in \Zee  \},\\
W(3, v)&:=\{(\pm \omega ^l ,t)\mid t\in \Gamma,
l\in \Zee  \}
\end{align*}
are infinite complex irreducible
crystallographic $r$-groups.

Analogously,
if $\Gamma$ is a~lattice of squares in $V$ ($= \Cee $), i.e.,


$$
\Gamma=
(\mathbb Z+i\mathbb Z)v\;\;\mbox{where $v\in V, v\neq 0$}
$$

\vskip\abovedisplayskip

\begin{center}
\includegraphics[width=70
mm]{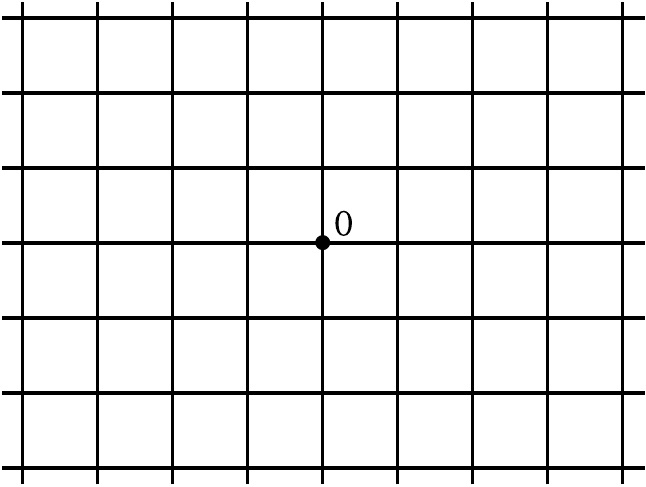}
\end{center}

\vskip\belowdisplayskip

\noindent
then
\begin{align*}
W(4,v)&:=\{(i^l ,t)\mid t\in \Gamma,
l \in \Zee  \}
\end{align*}
is an infinite complex irreducible
crystallographic $r$-group as well.

\subsubsection{Exercise}\label{exa1}
Prove that up to equivalence all infinite complex irreducible
1-dimensional $r$-groups are exhausted by groups of the following five types:
\begin{equation}
\label{ast}
(1)\; W(3, v),\;\;
(2)\; W(4, v),\;\;
(3)\; W(6,v),\;\;
(4)\; W(2,v, \lambda),\;\;
(5)\; W(2, v).
\tag{$\ast$}
\end{equation}

\subsubsection{Exercise}\label{exx}
Prove the following:

$\bullet$ In list  \eqref{ast}, the groups of different types
are not equivalent to each other.

$\bullet$ For every $d=2, 3, 4, 6$ and $v, u\in V$, the groups $W(d, v)$ and $W(d, u)$ are equivalent.

$\bullet$ Every $W(2,v,\lambda)$ is equivalent to a unique $W(2, v, \mu)$ where $\mu$ lies in the \lq\lq modular strip\rq\rq
\;$
\Omega =
\left \{z\in \Cee \mid {\rm Im}\,z >0;
-
1/2 \leqslant \re z<
1/2;
|z|\geqslant 1\text{ if }\re z \leqslant 0;
|z|>1 \text{ if }\re z>0 \right \}.
$

$\bullet$  For every $\lambda\in \Omega$ and  $v, u\in V$, the groups $W(2, v, \lambda)$ and $W(2, u, \lambda)$ are equivalent.

\vskip 4mm

Therefore, we see that there are continuous families of pairwise nonequivalent infinite $1$-dimensional complex irreducible $r$-groups depending on parameter $\lambda\in \Omega$.

\subsubsection{Exercise}
Prove the following:

$\bullet$ For every $d=2, 3, 4, 6$, the group  $W(d,v)$ is generated by two reflections.

$\bullet$ For every $v, \lambda$, the group $W(2,v, \lambda)$ is generated by three (but not by two) reflections.

\subsubsection{Exercise}\label{exa2}
For the groups from list \eqref{ast}, prove that the mirrors of reflections
and elements of $\Gamma$
are as depicted below.



\vskip 2mm

{\it Group} $W(3, v)$:

All reflections are of order $3$. Their mirrors are depicted as $\bullet$ and
$\odot\hskip -2.7mm\bullet$. The elements of $\Gamma$ are depicted as $\odot\hskip -2.7mm\bullet$.

\begin{figure}[htbp]
\begin{center}
\includegraphics[width=70
mm]{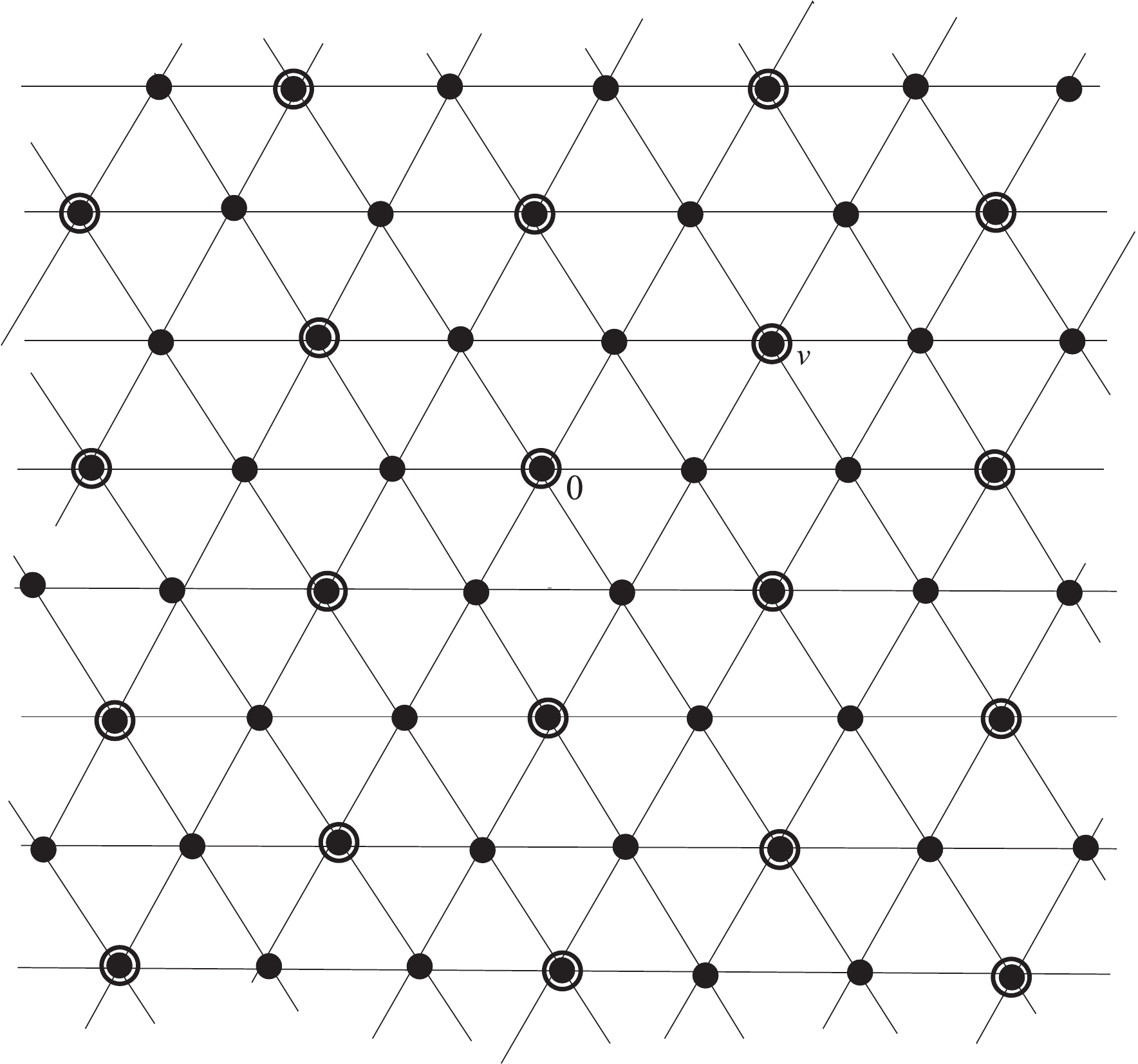}
\end{center}
\end{figure}


{\it Group}  $W(4, v)$:

 The mirrors of reflections of order $4$ are depicted as
 $\bullet$
 and $\odot\hskip -2.7mm\bullet$.
 The elements of $\Gamma$ are depicted as $\odot\hskip -2.7mm\bullet$.  The mirrors of reflections of order $2$
 are depicted as
$\circ$, $\bullet$,
 and $\odot\hskip -2.7mm\bullet$.

\begin{figure}[htbp]\begin{center}
\includegraphics[width=68
mm]{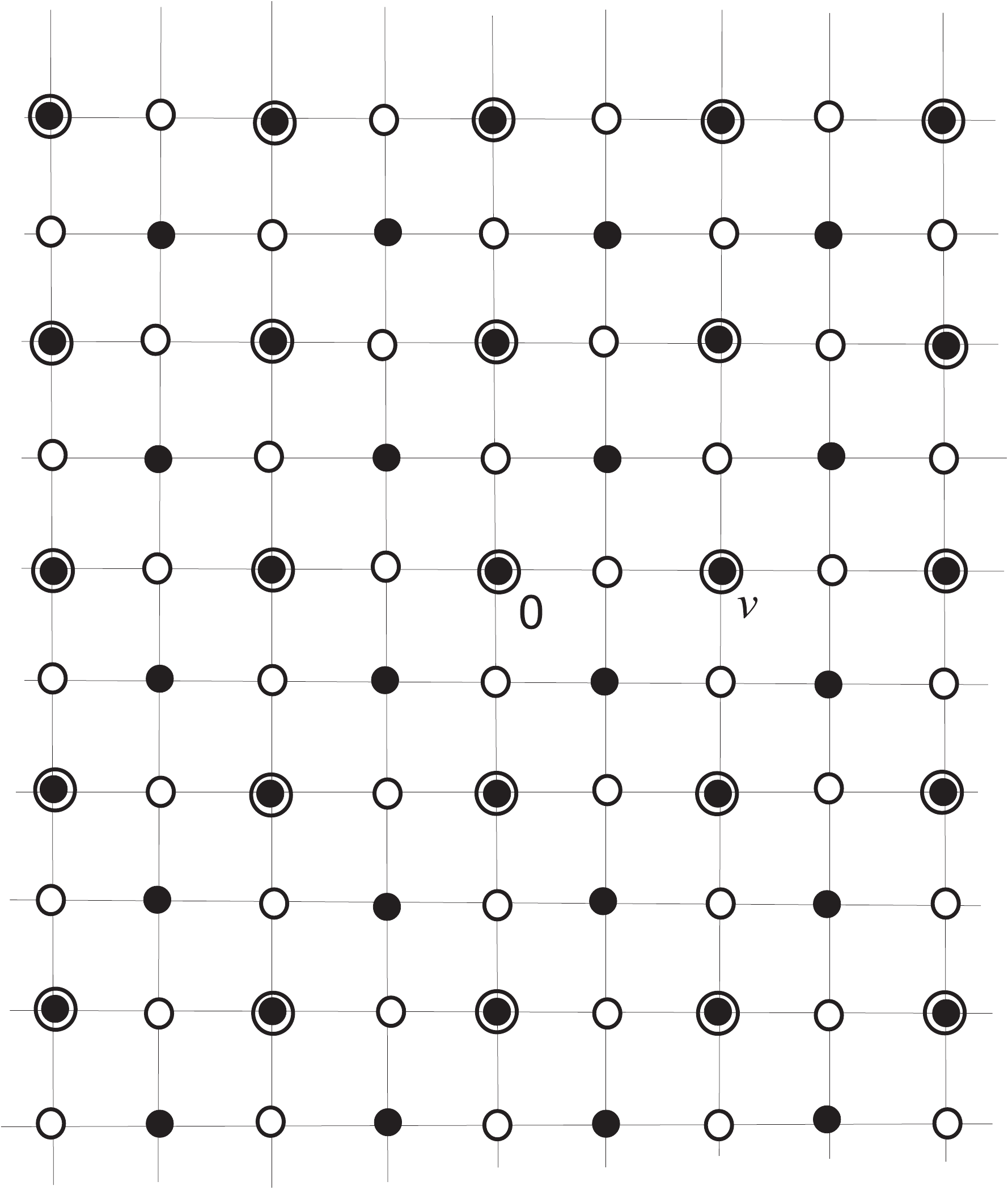}
\end{center}
\end{figure}


{\it Group} $W(6, v)$:

 The mirrors of reflections of order $6$ coincide with the elements of $\Gamma$ and are depicted as $\odot\hskip -2.7mm\bullet$. The mirrors of reflections of order $2$ are depicted as $\circ$ and $\odot\hskip -2.7mm\bullet$. The mirrors of reflections of order $3$ are depicted as
 $\bullet$ and $\odot\hskip -2.7mm\bullet$.


 \begin{figure}[htbp]\begin{center}
\includegraphics[width=75.5
mm]{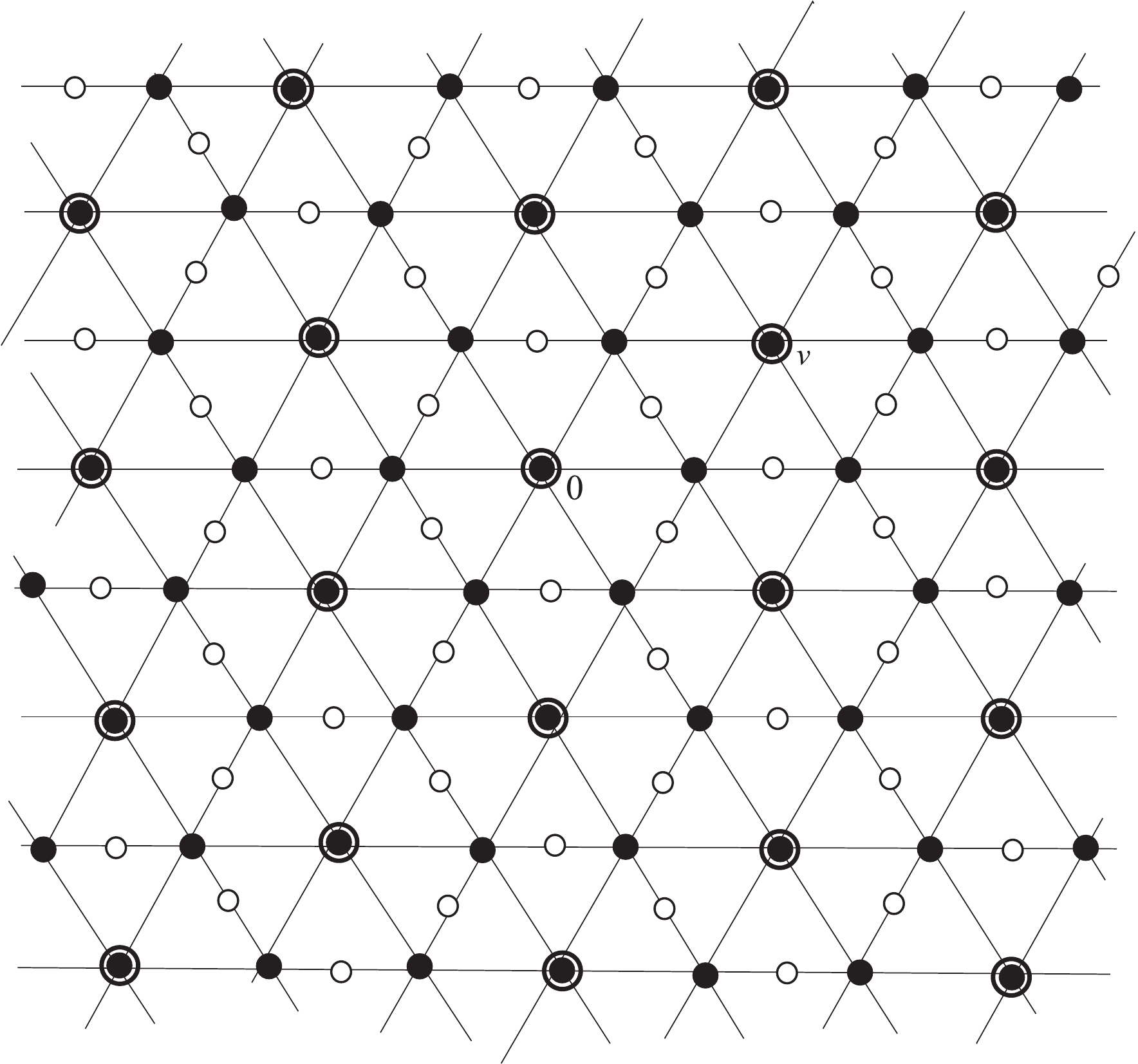}
\end{center}
\end{figure}


\eject

{\it Group} $W(2, v, \lambda)$:

  All reflections are of order $2$. Their mirrors are depicted as  $\bullet$ and $\odot\hskip -2.7mm\bullet$. The elements of $\Gamma$ are depicted as $\odot\hskip -2.7mm\bullet$.

\vskip 3mm

\begin{figure}[htbp]\begin{center}
\includegraphics[width=63
mm]{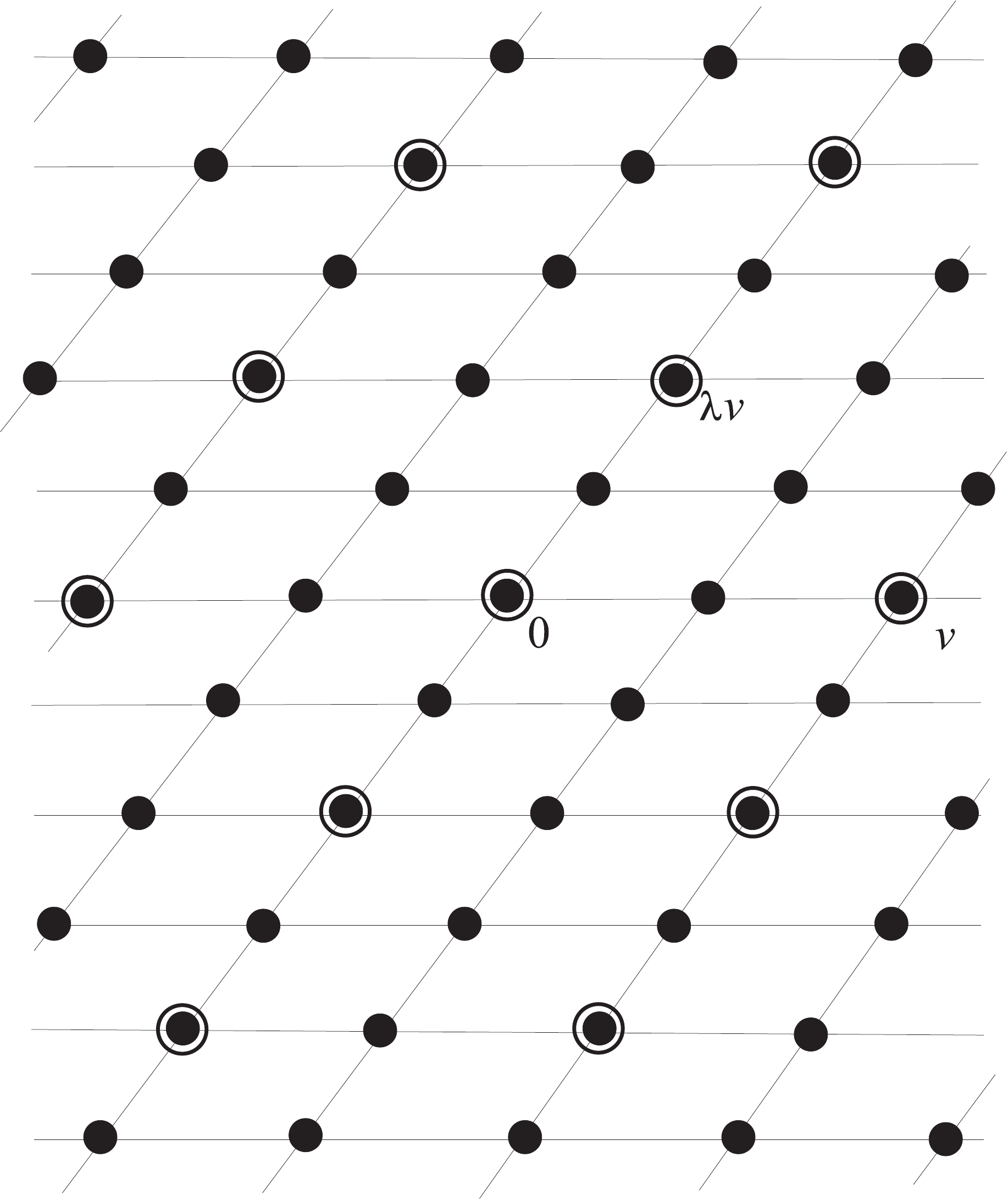}
\end{center}
\end{figure}

\vskip 2mm

{\it Group} $W(2, v)$:

  All reflections are of order $2$. Their mirrors are depicted as  $\bullet$ and $\odot\hskip -2.7mm\bullet$. The elements of $\Gamma$ are depicted as $\odot\hskip -2.7mm\bullet$.

\vskip 3mm

\begin{figure}[htbp]\begin{center}
\includegraphics[width=75
mm]{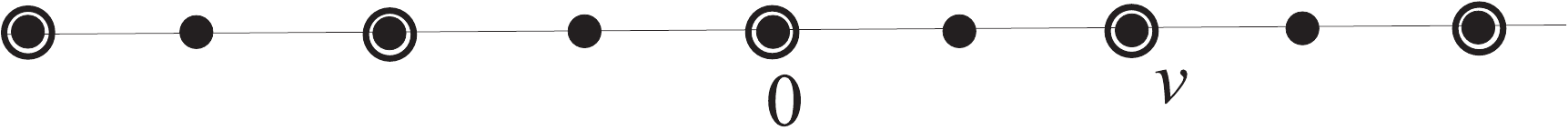}
\end{center}
\end{figure}


\section{\bf Formulation of the results}\label{Sfr}

In this section, we assume that $\ka =\Cee $.

\textit{Let $W$ be an irreducible infinite $r$-group}, $W \subseteq A(E)$. As we have seen in the example above, there are two possibilities: either $W$ is noncrystallographic (i.e., $E/W$ is not compact), or $W$ is crystallographic (i.e., $E/W$ is compact).

First, we shall describe the  structure of {\it noncrystallographic groups}. To do this, we need an auxiliary construction.

\subsection{Complexifications and real forms}

Let us consider $V$ as a~real vector space (of dimension $2n$). A linear subspace $V_{\Ree }$ of this real vector space is called \textit{a real form of} $V$ if

a) the natural map
$$
V_{\Ree }\otimes _{\Ree }\Cee \tto V
$$
is \textit{an isomorphism}, i.e., some (hence, any) $\Ree $-basis of $V_{\Ree }$ is a~$\Cee $-basis of $V$;

b) the restriction $\langle \,\cdot \mid \cdot\, \rangle \bigr |_{V_{\Ree }}$ of $\langle \,\cdot \mid \cdot\, \rangle $ to $V_{\Ree }$ is real-valued (hence, $V_{\Ree }$ is Euclidean with respect to $\langle \,\cdot \mid \cdot\, \rangle \bigr |_{V_{\Ree }}$).

If $V_{\Ree }$ is a~real form of $V$, then $V$ is \textit{the complexification} of $V_{\Ree }$.

Let $a\in E$ be a~point. We can consider $E$ as a~real affine space of dimension $2n$. The affine subspace
$$
E_{\Ree }=a+V_{\Ree }
$$
of this affine space is called \textit{a real form of} $E$, and $E$ is called \textit{the complexification} of $E_{\Ree }$.

It is clear that every real Euclidean linear (resp. affine) space is isomorphic to a~real form of a~certain complex Hermitian linear (resp. affine) space.

\subsubsection{Proposition}\label{jpre}
{\it The following properties hold:
\begin{enumerate}[\hskip 7.2mm\rm 1)]\itemsep=-.1ex
\item
${\rm Iso}(V)$ acts transitively on the set of real forms of $V$.
\item
The group of motions of $E$ acts transitively on the set of real forms of $E$.
\item
Every motion $\gamma $ of a~Euclidean affine space $E_{\Ree }$ can be extended in a~unique way to a~motion $\gamma _{\Cee }$ of $E$. This motion $\gamma _{\Cee }$ is called the complexification of $\gamma $ {\rm(}and $\gamma $ is called the real form of\;$\gamma _{\Cee }${\rm)}.
\item
$\dim _R H_\gamma = \dim _{\Cee }H_{\gamma _{\Cee }}$. Specifically, $\gamma $ is a~reflection if and only if $\gamma _{\Cee }$ is a~reflection.
\end{enumerate}
}

\vskip 2mm

\noindent {\it Proof} is left to the reader.\qed

\vskip 2mm

This proposition gives \textit{a method for constructing noncrystallographic infinite $r$-groups}. Indeed, let $G \subseteq A(E_{\Ree })$ be an infinite (real) $r$-group. Then  it is easy to see that
$$
G_{\Cee }=\{\gamma _{\Cee }\mid \gamma \in G\}\subseteq A(E)
$$
is an infinite complex noncrystallographic $r$-group (and $G_{\Cee }$ is irreducible if and only if $G$ is).

\subsection{Classification of infinite complex irreducible noncrystallographic $r$-groups: the result}

It appears that the construction above leads to every such group. More precisely, one has the following theorem (see also Section \ref{1.5},2).

\subsubsection{Theorem {\rm (Infinite complex irreducible noncrystallographic {\it r}-groups)}}\label{cncg}
{\it Let $W$ be an infinite complex irreducible
$r$-group. Then $W$ is non\-crystallographic if and only if it is equivalent to the complexification of an irreducible affine Weyl group.
}

\vskip 2mm

\noindent \textit{Proof} is given in Section  \ref{3.4}.\qed

\vskip 2mm

The description of {\it crystallographic groups} is much more complicated. In order to give this description we need some preparations and extra notation.

\subsection{Ingredients of the description}\label{2.3}

\textit{The subgroup of translations in $W$ will be denoted by} $\Tran W$.
$$
\Tran W = W \cap \Tran A(E),
$$
cf. Section  \ref{1.1}. It is clear that $\Tran W\lhd W$ and
$$
W/\Tran W \cong \Lin W.
$$

We usually identify $\Tran W$ with a~subgroup of the additive group of $V$ by means of the map $\gamma _v \mapsto v$. Clearly, this subgroup is a~$\Lin W$-\textit{invariant lattice} in $V$.

It will be proven in Section  \ref{3.1} that $\Tran W$ \textit{is a~full rank lattice} (i.e., of rank $2n$) and $\Lin W$ \textit{is a~finite group} (hence, $\Lin W$ is a~finite irreducible complex linear $r$-group, see Section  \ref{1.4}). Therefore, to describe $W$, one needs to point out a~group $\Lin W \subseteq \GL(V)$ from the Shephard and Todd list (i.e., from Theorem~\ref{cfi}),
a~$\Lin W$-invariant lattice $\Tran W \subseteq V$ of rank $2n$ and the way $\Lin W$ and $\Tran W$ are \lq\lq glued\rq\rq{} together. This is done below as follows:

\vskip 2mm

1) $\Lin W$ is given by its graph as in Sections\;\ref{1.6}--\ref{cfi}.

2) $\Tran W$ is given
as a set of some explicitly
described
linear combinations
of vectors $e_1,\ldots, e_s$.
Here, the set $R_j=R_{e_j,\theta _j}$ for $1 \leqslant j \leqslant s$ is a~fixed generating system of reflections of $\Lin W$ which is related to the graph of $\Lin W$ given in 1) as described in Sections~ \ref{1.6}--\ref{cfi}. To point out the vectors $e_1,\ldots, e_s$
explicitly, we assume that $V$ is a~subspace of a~standard Hermitian infinite-dimensional coordinate space $\Cee ^\infty $, i.e., the space, whose elements are the sequences ($a_1,a_2,\ldots $) of complex numbers with only a~finite number of nonzero elements $a_j$, and the~scalar product defined by the formula
$$
\langle (a_1,a_2,\ldots )\mid (b_1,b_2,\ldots )\rangle =\sum\limits _{j=1}^\infty a_j \overline b_j.
$$
The vectors
$e_1,\ldots, e_s$
are given by their coordinates in the~standard basis $\varepsilon _1,\varepsilon _2,\ldots $ of $\Cee ^\infty $, where $\varepsilon _j$ has a single nonzero entry $1$ in the position $j$:
$$
\varepsilon _j=(0,\ldots ,0,1,
0,\ldots ).
$$

3) The problem of
describing the \lq \lq glueing\rq \rq {} of $\Lin W$ and $\Tran W$ boils down to the determination of an extension of $\Tran W$ by $\Lin W$,
$$
0 \tto \Tran W \tto W \tto \Lin W \tto 1,
$$
that is described by means of cohomology. Let us show how it is done.

\subsection{Cohomology}\label{cohom}
Let $G$ be a~subgroup of $A(E)$ and write $T = \Tran G$, $K = \Lin G$. Choose a~point $a\in E$. Take $P\in K$ and let $\gamma \in G$ be such that $\Lin \gamma =P$. We have
$$
\kappa _a(\gamma )=(P,s(P)) \quad \text{where}\quad s(P) \in V.
$$
It is easy to see that the map
$$
\overline s\colon K \tto V/T,\qquad \overline s(P)=S(P)+T,
$$
is well-defined and is in fact a~1-cocycle, i.e.,
$$
\overline s(PQ)=\overline s(P)+P \overline s(Q) \quad \text{for all}\quad P,Q\in K
$$
(here, $K$ acts on $V/T$ in the natural way.)

Vice versa, if
$$
\overline r\colon K \tto V/T
$$
is an arbitrary 1-cocycle, let us consider an arbitrary map
$$
r\colon K \tto V
$$
such that
$$
\overline r(P)=r(P)+T \quad \text{for each}\quad P\in K.
$$
Then the set
$$
\{(P,r(P)+t )\mid t\in T,P\in K\}
$$
is a~subgroup $H$ of $A(E)$ with $\Lin H = K$ and $\Tran H = T$.

If we replace $a$ by an other point $b\in E$, then (see Section  \ref{1.1})
$$
\kappa _b(\gamma )=\kappa _a(\gamma _{a-b}\gamma \gamma _{b-a})=(P,s(P)+\underbrace {v-Pv}_{\text{1-coboundary}},\text{ where }v = a-b.
$$
\textit{Therefore, we have a~bijection between the set of $\Tran A(E)$-conjugacy classes of subgroups~$G$ of $A(E)$ with $\Lin G = K$, $\Tran G = T$ and the group} $H^1(K,V/T)$.

However, we have to consider subgroups of $A(E)$ \textit{up to equivalence}, i.e., up to $A(E)$-conjugation (and not just up to $\Tran A(E)$-conjugation!). This can be done as follows by means of an extra relation on $H^1(K,V/T)$. Let
$$
N(K,T) = \{Q\in \GL(V)\mid QKQ^{-1} = K,QT = T\}.
$$
If $Q\in N(K,T)$ and $\overline s\colon K \tto V/T$ is a~1-cocycle, resp. 1-coboundary, then it is easy to check that the map
$$
Q(\overline s) \colon K \tto V/T
$$
given by the formula
$$
Q(\overline s)(P)=Q \overline s(Q^{-1}PQ) \quad \text{for any}\quad P\in K
$$
is again a~1-cocycle, resp. 1-coboundary (here $Q$ acts on $V/T$ in the natural way). Therefore, we have an action of $N(K,T)$ on $H^1(K,V/T)$ (clearly, by means of automorphisms).

Let $\delta \in A(E)$ be such that $\Lin \delta G \delta ^{-1}=K$, $\Tran \delta G \delta ^{-1} = T$. We want to calculate the cocycle that corresponds to $\delta G \delta ^{-1}$. Changing $\delta $ to $\delta \gamma _v$, where $v = (\Lin \delta ^{-1} ) (a - \delta (a))$, we can assume that
$$
\kappa _a(\delta )=(Q,0) \quad \text{for}\quad Q\in N(K,T).
$$
Let $P \in K$ and $\lambda \in G$ be such that $\kappa _a(\lambda )=(Q^{-1}PQ,s(Q^{-1}PQ))$. Then
\[
\kappa _a(\delta \lambda \delta ^{-1})=(P,Qs(Q^{-1}PQ)).
\]
Therefore, the cocycle corresponding to $\delta G \delta ^{-1}$ is $Q(\overline s)$, where $\overline s$ is the cocycle corresponding to $G$.

\vskip 2mm

We see now that \textit{there is a~bijection between the set of classes of equivalent subgroups $ G \subseteq A(E)$ with $\Lin G = K$, $\Tran G = T$ and the set of $N(K,T)$-orbits in} $H^1(K,V/T)$.

\vskip 2mm

With all these facts in mind, \textit{we determine the extension $W$ {\rm(}of $\Tran W$ by $\Lin W${\rm)}
 by pointing out a~$1$-cocycle which represents the corresponding element of}
 $H^1 (\Lin W,V/\Tran W)$ (in fact, the whole $N(\Lin W,\Tran W)$-orbit in
 $H^1 (\Lin W,V/\Tran W)$). In order to do so, \textit{we need only to specify
 the values of this $1$-cocycle on the elements of a~generating system of reflections of} $\Lin W$. Technically it is more convenient to realize it as follows.

Let $\widetilde {\Lin W}$ be a~free group with free generating set
$r_1,\ldots, r_s$. We have the epimorphism
\[
\text{$\phi \colon \widetilde {\Lin W}\tto \Lin W$,\quad $\phi (r_j) = R_j$ for every $1 \leqslant j \leqslant s$.}
\]
 The kernel of $\phi $ is the subgroup of \lq\lq relations\rq\rq{} of $\Lin W$. This epimorphism determines
 in the natural way
 an action of $\widetilde {\Lin W}$ on $V$. A 1-cocycle $c$ of $\widetilde{\Lin W}$ with values 
 in $V$ is given by its values on the generators $r_j$, i.e., by the elements
 $$
 c(r_1),\ldots, c(r_s)\in V,
 $$
 and these values may be arbitrary (because $\widetilde {\Lin W}$ is free). It is easy to see that the formula
$$
R_j\mapsto c(r_j)+\Tran W,\quad \text{where}\quad 1 \leqslant j \leqslant s,
$$
defines a~1-cocycle of $\Lin W$ with values in $V/\Tran W$ if and only if $c(F) \in \Tran W$ for every $F \in \Ker \phi $. It is also clear that every 1-cocycle of $\Lin W$ with values in $V/\Tran W$ is obtained in such a~way.

\textit{We describe the extension $W$ {\rm(}of $\Tran W$ by $\Lin W${\rm)} by specifying
the vectors $c (r_1),\ldots\break\ldots, c (r_s )$.}

\vskip 2mm

We are now ready to formulate the results of the classification of infinite irreducible crystallogra\-phic $r$-groups.

\vskip 2mm

\textit{We denote by $K_d$ the finite linear irreducible $r$-group whose  number in the list of Shephard and Todd {\rm(}i.e., in the first column of Table~{\rm 1)}
is $d$}.
Note that this notations is in slight confusion with Cohen's notation $K_5,K_6$.



\subsection{Description of the group of linear parts: the result} First of all, there is an analogue of Theorem~\ref{1.5}.

\subsubsection{Theorem {\rm(Linear parts of infinite complex irreducible {\it r}-groups)}}\label{cirg2}
{\it
Let $K \subseteq \GL(V)$ be a finite complex irreducible $r$-group. Then the following properties are equivalent:
\begin{enumerate}[\hskip 7.2mm\rm a)]\itemsep=-.1ex
\item There exists a~nonzero $K$-invariant lattice in $V$.
\item There exists a~$K$-invariant lattice of rank $2n$ in $V$.
\item $K = \Lin W$, where $W$ is an infinite complex irreducible crystallographic $r$-group.
\item The ring with unity, generated over $\Zee  $ by all cyclic products of a~graph of $K$, lies in the ring of algebraic integers of a~purely imaginary quadratic extension of $\Qee $.
\item $K$ is defined over a~purely imaginary quadratic extension of $\Qee $.
\item $\begin{array}[t]{@{}l@{}}
K \text{is one of the groups:}\\
K_1;K_2\ (m=2,3,4,6);K_3\ (m=2,3,4,6);\\[-.3mm]
K_4;K_5;K_8;K_{12};K_{24};K_{25};K_{26};K_{28};K_{29};\\[-.3mm]
K_{31};K_{32};K_{33};K_{34};K_{35};K_{36};K_{37}.
\end{array}
$
\end{enumerate}
}

\vskip 2mm

\noindent{\it Proof}
is given in Section \ref{4.6}.

\vskip 2mm

Now we describe the crystallographic groups themselves.

\subsection{The list of  infinite complex irreducible crystallographic
\boldmath$r$-groups}

This list is given by Theorem \ref{cirg} below, in which we use the following notation:
$$
\Omega =\left \{z\in \Cee \mid {\rm Im}\,z >0;
-
1/2 \leqslant \re z<
1/2;
|z|\geqslant 1\text{ if }\re z \leqslant 0;
|z|>1 \text{ if }\re z>0 \right \}
$$
($\Omega$ is the \lq \lq modular strip\rq\rq; see Exercise \ref{exx});
$$[\alpha ,\beta ]=\{a \alpha +b \beta \mid a, b\in \Zee  \}\quad \mbox{for arbitrary $\alpha ,\beta \in \Cee $.}
$$

\subsubsection{Theorem {\rm (Infinite complex irreducible crystallographic {\it r}-groups)}}\label{cirg}
{\it Table $2$ below
gives the complete list of  infinite complex irreducible cry\-stal\-lographic
$r$-groups $W$ {\rm(}considered up to equivalence{\rm)}.
}

\vskip 2mm

\noindent{\it Proof}
(in which the relevant calculations are omitted)  is given in the subsequent sections.


\begin{landscape}\addcontentsline{toc}{subsubsection}{Table 2. The irreducible infinite complex crystallographic
$r$-groups}

\centering \footnotesize
\begin{longtable}{|l@{\extracolsep{.01mm}}|@{\extracolsep{.3mm}}
c@{\extracolsep{-.6mm}}|l|l|l@{\extracolsep{-.3mm}}|@{\extracolsep{-.3mm}}c@{\extracolsep{-.3mm}}|}
\multicolumn{6}{l}{\textbf{Table 2.} The irreducible infinite complex crystallographic 
$r$-groups}\label{tab2}\\[2mm]
\hline
\multicolumn{1}{|c|}{Notation of $W\vphantom {\displaystyle \frac{1}{2}} $}&$\,\begin{matrix}\dim W\\[-2.5mm] =n \end{matrix}\,$&\multicolumn{1}{c|}{$\Lin W$}&\multicolumn{1}{c|}{$\Tran W$}&\multicolumn{1}{c|}{$e_1,\ldots ,e_s$}&cocycle $c$\\
\hline
\endfirsthead
\hline
\multicolumn{1}{|c|}{Notation of $W\vphantom {\displaystyle \frac{1}{2}} $}&$\begin{matrix}\dim W\\[-2.5mm] =n \end{matrix}$&\multicolumn{1}{c|}{$\Lin W$}&\multicolumn{1}{c|}{$\Tran W$}&\multicolumn{1}{c|}{$e_1,\ldots ,e_s$}&cocycle $c$\\
\hline
\endhead
\hline
\endfoot
$\begin{array}{@{}l@{}}
{}\\[-12pt]
[{\sf A}_s]^\alpha \\
s \geqslant 1\\[-12pt]
{}
\end{array}
$&$s$&$\begin{array}{@{}l@{}}
{}\\[-12pt]
K_1, \\[-1mm]
\text{type }{\sf A}_s\\[-1mm]
s \geqslant 1\\[-12pt]
{}
\end{array}
$&$\begin{array}{@{}l@{}}
{}\\[-12pt]
[1,\alpha ]e_1+\cdots +[1,\alpha ]e_s, \\
\alpha \in \Omega \\[-12pt]
{}
\end{array}
$&$\begin{array}{@{}l@{}}
{}\\[-12pt]
e_j=(\varepsilon _j-\varepsilon _{j+1})/\sqrt 2\\
j=1,\ldots ,s\\[-12pt]
{}
\end{array}
$&\raisebox{-160pt}[0pt][0pt]{$c=0$}\\
\cline{1-5}
$\begin{array}{@{}l@{}}
[G(2,1,s)]_1^\alpha \\
s \geqslant 3\\[-12pt]
{}
\end{array}
$&&&$\begin{array}{@{}l@{}}
{}\\[-12pt]
[1,\alpha ]e_1+[1,\alpha ]\sqrt 2\,e_2+\cdots +[1,\alpha ]\sqrt s\,e_s, \\
\alpha \in \Omega \\[-12pt]
{}
\end{array}
$&&\\
\cline{1-1}\cline{4-4}
$\begin{array}{@{}l@{}}
{}\\[-12pt]
[G(2,1,s)]_2^\beta \\
s \geqslant 3\\[-12pt]
{}
\end{array}
$&&&
$\begin{array}{@{}l@{}}
{}\\[-12pt]
[1,\beta ]e_1+\left [1,\frac{1+\beta }{2} \right ]\!\!\sqrt 2\,e_2+;\delta +\left [1,\frac{1+\beta }{2} \right ]\!\!\sqrt 2\,e_s, \\
\beta \in \Omega \\[-12pt]
{}\\[-12pt]
{}
\end{array}
$&&\\
\cline{1-1}\cline{4-4}
$\begin{array}{@{}l@{}}
{}\\[-12pt]
[G(2,1,s)]_3^\gamma \\
s \geqslant 3\\[-12pt]
{}
\end{array}
$&$s$&\raisebox{0pt}[0pt][0pt]{$\begin{array}{@{}l@{}}
{}\\[-12pt]
K_2\\[-1mm]
\text{type} \\[-1mm]
G(2,1,s) \\[-1mm]
s \geqslant 3\\[-12pt]
{}
\end{array}
$}&$\begin{array}{@{}l@{}}
{}\\[-12pt]
[1,\gamma ]e_1+\left [\frac{1}{2},\gamma \right ]\sqrt 2\,e_2+\cdots +\left [\frac{1}{2},\gamma \right ]\sqrt 2\,e_s, \\
\gamma \in \Omega \\[-12pt]
{}
\end{array}
$&\raisebox{-67pt}[0pt][0pt]{$\begin{array}{@{}l@{}}
{}\\[-12pt]
e_1=\varepsilon _1, \\
e_j=(\varepsilon _{j-1}-\varepsilon _j)/\sqrt 2\\
j=2,\ldots ,s\\[-12pt]
{}
\end{array}
$}&\\
\cline{1-1}\cline{4-4}
$\begin{array}{@{}l@{}}
{}\\[-12pt]
[G(2,1,s)]_4^\delta \\
s \geqslant 3\\[-12pt]
{}
\end{array}
$&&&$\begin{array}{@{}l@{}}
{}\\[-12pt]
[1,\delta ]e_1+\left [1,\frac{\delta}{2} \right ]\sqrt 2\,e_2+\cdots +\left [1,\frac{\delta}{2} \right ]\sqrt 2\,e_s, \\
\delta \in \Omega \\[-12pt]
{}
\end{array}
$&&\\
\cline{1-1}\cline{4-4}
$\begin{array}{@{}l@{}}
{}\\[-12pt]
[G(2,1,s)]_5^\lambda \\
s \geqslant 3\\[-12pt]
{}
\end{array}
$&&&$\begin{array}{@{}l@{}}
{}\\[-12pt]
[1,\lambda ]e_1+\left [\frac{1}{2}, \frac{\lambda}{2} \right ]\sqrt 2\,e_2+\cdots +\left [\frac{1}{2},\frac{\lambda}{2} \right ]\sqrt 2\,e_s, \\
\lambda \in \Omega \\[-12pt]
{}
\end{array}
$&&\\
\cline{1-4}
$\begin{array}{@{}l@{}}
{}\\[-12pt]
[G(3,1,s)]_1\\
s \geqslant 2\\[-12pt]
{}
\end{array}
$&\raisebox{-15pt}[0pt][0pt]{$s$}&\raisebox{-15pt}[0pt][0pt]{$\begin{array}{@{}l@{}}
{}\\[-12pt]
K_2\\[-1mm]
\text{type} \\[-1mm]
G(3,1,s) \\[-1mm]
s \geqslant 2\\[-12pt]
{}
\end{array}
$}&$[1,\omega ]e_1+[1,\omega ]\sqrt 2\,e_2+\cdots +[1,\omega ]\sqrt 2\,e_s$&&\\
\cline{1-1}\cline{4-4}
$\begin{array}{@{}l@{}}
{}\\[-12pt]
[G(3,1,s)]_2 \\
s \geqslant 2\\[-12pt]
{}
\end{array}
$&&&$[1,\omega ]e_1+[1,\omega ]i \sqrt {\frac{2}{3}}\,e_2+\cdots +[1,\omega ]i \sqrt {\frac{2}{3}}\,e_s$&&\\
\cline{1-4}
$\begin{array}{@{}l@{}}
{}\\[-12pt]
[G(4,1,s)]_1\\
s \geqslant 2\\[-12pt]
{}
\end{array}
$&\raisebox{-15pt}[0pt][0pt]{$s$}&\raisebox{-15pt}[0pt][0pt]{$\begin{array}{@{}l@{}}
{}\\[-12pt]
K_2\\[-1mm]
\text{type} \\[-1mm]
G(4,1,s) \\[-1mm]
s \geqslant 2\\[-12pt]
{}
\end{array}
$}&$[1,i]e_1+[1,i]\sqrt 2\,e_2+\cdots +[1,i]\sqrt 2\,e_s$&
&
\\
\cline{1-1}\cline{4-4}
$\begin{array}{@{}l@{}}
{}\\[-12pt]
[G(4,1,s)]_2 \\
s \geqslant 2\\[-12pt]
{}
\end{array}
$&&&$[1,i]e_1+[1,i]\varepsilon e_2+\cdots +[1,i]\varepsilon e_s$&&\\[12mm]
\cline{1-4}
$\begin{array}{@{}l@{}}
{}\\[-12pt]
[G(6,1,s)] \\
s \geqslant 2\\[-12pt]
{}
\end{array}
$&$s$&$\begin{array}{@{}l@{}}
{}\\[-12pt]
K_2\\[-1mm]
\text{type}\\[-1mm]
G(6,1,s) \\[-1mm]
s \geqslant 2\\[-12pt]
{}
\end{array}
$&$[1,\omega ]e_1+[1,\omega ]\sqrt 2\,e_2+\cdots +[1,\omega ]\sqrt 2\,e_s$&
$\begin{array}{@{}l@{}}
{}\\[-12pt]
e_1=\varepsilon _1, \\
e_j=(\varepsilon _{j-1}-\varepsilon _j)/\sqrt 2\\
j=2,\ldots ,s\\[-12pt]
{}
\end{array}
$&\\
\cline{1-5}
$\begin{array}{@{}l@{}}
{}\\[-12pt]
[G(2,2,s)]^\alpha \\
s \geqslant 3\\[-12pt]
{}
\end{array}
$&$s$&$\begin{array}{@{}l@{}}
{}\\[-12pt]
K_2\\[-1mm]
\text{type}\\[-1mm]
G(2,2,s) \\[-1mm]
s \geqslant 3\\[-12pt]
{}
\end{array}
$&$\begin{array}{@{}l@{}}
{}\\[-12pt]
[1,\alpha ]e_1+\cdots +[1,\alpha ]e_s, \\
\alpha \in \Omega \\[-12pt]
{}
\end{array}
$&$\begin{array}{@{}l@{}}
{}\\[-12pt]
e_1=-(\varepsilon _1+\varepsilon _2)/\sqrt 2, \\
e_j=(\varepsilon _{j-1}-\varepsilon _j)/\sqrt 2, \\
j=2,\ldots ,s\\[-12pt]
{}
\end{array}
$&\\
\cline{1-5}
$\begin{array}{@{}l@{}}
{}\\[-12pt]
[G(3,3,s)] \\
s \geqslant 3\\[-12pt]
{}
\end{array}
$&$s$&$\begin{array}{@{}l@{}}
{}\\[-12pt]
K_2, \\[-1mm]
\text{type} \\[-1mm]
G(3,3,s) \\[-1mm]
s \geqslant 3\\[-12pt]
{}
\end{array}
$&$[1,\omega ]e_1+\cdots +[1,\omega ]e_s$&$\begin{array}{@{}l@{}}
{}\\[-12pt]
e_1=(\omega \varepsilon _1-\varepsilon _2)\sqrt{2}, \\
e_j=(\varepsilon _{j-1}-\varepsilon _j)/\sqrt 2, \\
j=2,\ldots ,s\\[-12pt]
{}
\end{array}
$&\\
\cline{1-5}
$\begin{array}{@{}l@{}}
{}\\[-12pt]
[G(4,4,s)] \\
s \geqslant 3\\[-12pt]
{}
\end{array}
$&$s$&$\begin{array}{@{}l@{}}
{}\\[-12pt]
K_2\\[-1mm]
\text{type} \\[-1mm]
G(4,4,s) \\[-1mm]
s \geqslant 3\\[-12pt]
{}
\end{array}
$&$[1,i]e_1+\cdots +[1,i]e_s$&$\begin{array}{@{}l@{}}
{}\\[-12pt]
e_1=(i \varepsilon _1-\varepsilon _2)/\sqrt 2, \\
e_j=(\varepsilon _{j-1}-\varepsilon _j)/\sqrt 2, \\
j=2,\ldots ,s\\[-12pt]
{}
\end{array}
$&\\
\cline{1-5}
$\begin{array}{@{}l@{}}
{}\\[-12pt]
[G(6,6,s)] \\
s \geqslant 3\\[-12pt]
{}
\end{array}
$&$s$&$\begin{array}{@{}l@{}}
{}\\[-12pt]
K_2, \\[-1mm]
\text{type} \\[-1mm]
G(6,6,s) \\[-1mm]
s \geqslant 3\\[-12pt]
{}
\end{array}
$&$[1,\omega ]e_1+\cdots +[1,\omega ]e_s$&$\begin{array}{@{}l@{}}
{}\\[-12pt]
e_1=((1+\omega )\varepsilon _1-\varepsilon _2)/\sqrt 2\\
e_j=(\varepsilon _{j-1}-\varepsilon _j)/\sqrt 2\\
j=2,\ldots ,s\\[-12pt]
{}
\end{array}
$&\raisebox{52pt}[0pt][0pt]{$c=0$}
\\
\cline{1-5}
$[G(2,1,2)]_1^\alpha $&&\raisebox{-35pt}[0pt][0pt]{$\begin{array}{@{}l@{}}
{}\\[-12pt]
K_2, \\[-1mm]
\text{type} \\[-1mm]
G(2,1,2) \\[-1mm]
=\text{type} \\[-1mm]
G(4,4,2)\\[-12pt]
{}
\end{array}
$}&$\begin{array}{@{}l@{}}
{}\\[-12pt]
\lbrack 1,\alpha ]e_1+[1,\alpha ]\sqrt 2\,e_2, \\
\alpha \in \Omega \\[-12pt]
{}
\end{array}
$&&\\
\cline{1-1}\cline{4-4}
$[G(2,1,2)]_2^\beta $&2&&$\begin{array}{@{}l@{}}
{}\\[-12pt]
[1,\beta ]e_1+\left [1,\frac{\beta}{2} \right ]\sqrt 2\,e_2, \\
\beta \in \Omega \\[-12pt]
{}
\end{array}
$&$\begin{array}{@{}l@{}}
{}\\[-12pt]
e_1=\varepsilon _1\\
e_2=(\varepsilon _1-\varepsilon _2)/\sqrt 2\\[-12pt]
{}
\end{array}
$&\\
\cline{1-1}\cline{4-4}
$[G(2,1,2)]_3^\gamma $&&&$\begin{array}{@{}l@{}}
{}\\[-12pt]
[1,\gamma ]e_1+\left [1,\frac{1+\gamma }{2} \right ]\sqrt 2\,e_2, \\
\gamma \in \Omega \\[10pt]
{}
\end{array}
$&&\\[20mm]
\cline{1-5}
$[G(6,6,2)]_1^\alpha $&&\raisebox{-55pt}[-0pt][0pt]{$\begin{array}{@{}l@{}}
{}\\[-12pt]
K_2, \\[-1mm]
\text{type} \\[-1mm]
G(6,6,2)\\[-12pt]
{}
\end{array}
$}&$\begin{array}{@{}l@{}}
{}\\[-12pt]
[1,\alpha ]e_1+[1,\alpha ](2+\omega )e_2, \\
\alpha \in \Omega \\[-12pt]
{}
\end{array}
$&&\\
\cline{1-1}\cline{4-4}
$[G(6,6,2)]_2^\beta $&\raisebox{-18pt}[-0pt][0pt]{2}&&$\begin{array}{@{}l@{}}
{}\\[-12pt]
[1,\beta ]e_1+\left [1,\frac{\beta}{3} \right ](2+\omega )e_2, \\
\beta \in \Omega \\[-12pt]
{}
\end{array}
$&\raisebox{-18pt}[-0pt][0pt]{$\begin{array}{@{}l@{}}
{}\\[-12pt]
e_1=((1+\omega )\varepsilon _1-\varepsilon _2)/\sqrt 2 \\
e_2=(\varepsilon _1-\varepsilon _2)/\sqrt 2\\[-12pt]
{}
\end{array}
$}&\\
\cline{1-1}\cline{4-4}
$[G(6,6,2)]_3^\gamma $&&&$\begin{array}{@{}l@{}}
{}\\[-12pt]
[1,\gamma ]e_1+\left [1,\frac{1+\gamma }{3} \right ](2+\omega )e_2, \\
\gamma \in \Omega \\[-12pt]
{}
\end{array}
$&&\\
\cline{1-1}\cline{4-4}
$[G(6,6,2)]_4^\delta $&&
&$\begin{array}{@{}l@{}}
{}\\[-12pt]
[1,\delta ]e_1+\left [1,\frac{2+\delta }{3} \right ](2+\omega )e_2, \\
\delta \in \Omega \\[-12pt]
{}
\end{array}
$&
&
\\
\cline{1-5}
$\begin{array}{@{}l@{}}
{}\\[-12pt]
[G(4,2,s-1)]_1\\
s \geqslant 3\\[-12pt]
{}
\end{array}
$&&\raisebox{-47pt}[0pt][0pt]{$\begin{array}{@{}l@{}}
{}\\[-12pt]
K_2, \\[-1mm]
\text{type} \\[-1mm]
G(4,2,s-1) \\[-1mm]
s \geqslant 3\\[-12pt]
{}
\end{array}
$}&\raisebox{-23pt}[0pt][0pt]{$T=[1,i]e_1+\cdots +[1,i]e_{s-1}$}&\raisebox{-77pt}[0pt][0pt]{$\begin{array}{@{}l@{}}
{}\\[-12pt]
e_1=(i \varepsilon _1-\varepsilon _2)/\sqrt 2\\
e_j=(\varepsilon _{j-1}-\varepsilon _j)/\sqrt 2\\
j=2,\ldots ,s-1\\
e_s=\varepsilon _{s-1}\\[-12pt]
{}
\end{array}
$}&
\raisebox{75pt}[0pt][0pt]{$c=0$}
\\
\cline{1-1}\cline{6-6}
$\begin{array}{@{}l@{}}
{}\\[-12pt]
[G(4,2,s-1)]_1^* \\
s \geqslant 3\\[-12pt]
{}
\end{array}
$&$s-1$&&&&$\begin{array}{@{}l@{}}
{}\\[-12pt]
c(r_j)=0, \\
j=1,\ldots ,s-1\\
c(r_s)=e_s/\sqrt 2\\[-12pt]
{}
\end{array}
$\\
\cline{1-1}\cline{4-4}\cline{6-6}
$\begin{array}{@{}l@{}}
{}\\[-12pt]
[G(4,2,s-1)]_2\\
s \geqslant 3\\[-12pt]
{}
\end{array}
$&&&$\begin{array}{@{}l@{}}
{}\\[-12pt]
T \cup \left (T+\frac{1+i}{2}\,(e_1+e_2)\right )\\
=[1,i]e_1+\cdots +[1,i]e_{s-1}+\frac{1}{\sqrt 2}\,[1,i]e_s\\[-12pt]
{}
\end{array}
$&&\\
\cline{1-4}
$[G(4,2,2)]_3$&2&$\begin{array}{@{}l@{}}
{}\\[-12pt]
K_2, \\[-1mm]
\text{type} \\[-1mm]
G(4,2,2)\\[-12pt]
{}
\end{array}
$&$[1,i]e_1+[1,i](1+i)e_2$&&\\
\cline{1-5}
$\begin{array}{@{}l@{}}
{}\\[-12pt]
[G(6,2,s-1)]_1\\
s \geqslant 3\\[-12pt]
{}
\end{array}
$&$s-1$&$\begin{array}{@{}l@{}}
{}\\[-12pt]
K_2 \\[-1mm]
\text{type} \\[-1mm]
G(6,2,s-1) \\[-1mm]
s \geqslant 3\\[-12pt]
{}
\end{array}
$&$[1,\omega ]e_1+\cdots +[1,\omega ]e_{s-1}$&\raisebox{-20pt}[0pt][0pt]{$\begin{array}{@{}l@{}}
{}\\[-12pt]
e_1=((1+\omega )\varepsilon _1-\varepsilon _2)/\sqrt 2, \\
e_j=(\varepsilon _{j-1}-\varepsilon _j)/\sqrt 2, \\
j=2,\ldots ,s-1\\
e_s=\varepsilon _{s-1}\\[-12pt]
{}
\end{array}
$}&\raisebox{18pt}[0pt][0pt]{$c=0$}\\
\cline{1-4}
$[G(6,2,2)]_2$&2&$\begin{array}{@{}l@{}}
{}\\[-12pt]
K_2, \\[-1mm]
\text{type} \\[-1mm]
G(6,2,2)\\[-12pt]
{}
\end{array}
$&$[1,\omega ]e_1+[1,\omega ](2+\omega )e_2$&&\\
\cline{1-4}
$\begin{array}{@{}l@{}}
{}\\[-12pt]
[G(6,3,s-1)]_1 \\
s \geqslant 3\\[-12pt]
{}
\end{array}
$&$s-1$&$\begin{array}{@{}l@{}}
{}\\[-12pt]
K_2, \\[-1mm]
\text{type} \\[-1mm]
G(6,3,s-1) \\[-1mm]
s \geqslant 3\\[-12pt]
{}
\end{array}
$&$[1,\omega ]e_1+\cdots +[1,\omega ]e_{s-1}$&
\raisebox{-20pt}[0pt][0pt]{$\begin{array}{@{}l@{}}
{}\\[-12pt]
e_1=((1+\omega )\varepsilon _1-\varepsilon _2)/\sqrt 2, \\
e_j=(\varepsilon _{j-1}-\varepsilon _j)/\sqrt 2, \\
j=2,\ldots ,s-1, \\
e_s=\varepsilon _{s-1}\\[-12pt]
{}
\end{array}
$}
&
\\
\cline{1-4}
$[G(6,3,2)]_2$&2&$\begin{array}{@{}l@{}}
{}\\[-12pt]
K_2, \\[-1mm]
\text{type} \\[-1mm]
G(6,3,2)\\[-12pt]
{}
\end{array}
$&$[1,2 \omega ]e_1+[2,\omega ]e_2$&&\\
\cline{1-5}
$[K_3(3)]$&&$\begin{array}{@{}l@{}}
{}\\[-12pt]
K_3\\[-1mm]
m=3\\[-12pt]
{}
\end{array}
$&$[1,\omega ]e_1$&&\\
\cline{1-1}\cline{3-4}
$[K_3(4)]$&1&$\begin{array}{@{}l@{}}
{}\\[-12pt]
K_3\\[-1mm]
m=4\\[-12pt]
{}
\end{array}
$&$[1,i]e_1$&$e_1=\varepsilon _1$&\\
\cline{1-1}\cline{3-4}
$[K_3(6)]$&&$\begin{array}{@{}l@{}}
{}\\[-12pt]
K_3, \\[-1mm]
m=6\\[-12pt]
{}
\end{array}
$&$[1,\omega ]e_1$&&\\
\cline{1-5}
$[K_4]$&
&$K_4$&$[1,\omega ]e_1+[1,\omega ]e_2$&$\begin{array}{@{}l@{}}
{}\\[-12pt]
e_1=\varepsilon _1,\;\; 
e_2=\frac{1-\omega }{3}\,(\varepsilon _1+\varepsilon _2+\varepsilon _3)\\[-12pt]
{}
\end{array}
$&\\
\cline{1-1}\cline{3-5}
$[K_5]$&&$K_5$&$[1,\omega ]e_1+[1,\omega ]\sqrt 2\,e_2$&$\begin{array}{@{}l@{}}
{}\\[-12pt]
e_1=\varepsilon _1, \;\;
e_2=\frac{1-\omega }{3}\,(\sqrt 2\,\varepsilon _1+\varepsilon _2)\\[-12pt]
{}
\end{array}
$&\\
\cline{1-1}\cline{3-5}
$[K_8]$&&$K_8$&$[1,i]e_1+[1,i]e_2$&$\begin{array}{@{}l@{}}
{}\\[-12pt]
e_1=\varepsilon _1,\;\; 
e_2=\frac{1-i}{2}\,(\varepsilon _1-\varepsilon _2)\\[-12pt]
{}
\end{array}
$&
\\
\cline{1-1}\cline{3-5}
$[K_{12}] \vphantom {\displaystyle \frac{1}{2}} $&\raisebox{14pt}[0pt][0pt]{2}&\raisebox{-16pt}[0pt][0pt]{$K_{12}$}&\raisebox{-16pt}[0pt][0pt]{$[1,i \sqrt 2]e_1+[1,i \sqrt 2]e_2$}&
$\begin{array}{@{}l@{}}
{}\\[2pt]
e_1=\frac{1}{\sqrt 2}\,\varepsilon _1+\frac{1+i}{2}\,\varepsilon _2, \\
e_2=\frac{\sqrt 2+(\sqrt 2-2)i}{4}\varepsilon _1+\frac{2+\sqrt 2-\sqrt 2i}{4}\varepsilon _2, \\
e_3=\frac{1}{\sqrt 2}\,\varepsilon _1+\frac{1-i}{2}\,\varepsilon _2\\[-30pt]
{}
\end{array}
$
&\raisebox{126pt}[0pt][0pt]{$c=0$}\\
\cline{1-1}\cline{6-6}
$[K_{12}]^*$&&&&&$\begin{array}{@{}l@{}}
{}\\[-12pt]
c(r_1)=c(r_2)=0\\
c(r_3)=\frac{1+i}{2}\,e_3\\[-12pt]
{}
\end{array}
$
\\
\cline{1-6}
$[K_{24}]$&
3
&$K_{24}$&$\begin{array}{@{}l@{}}
{}\\[-12pt]
\left [1,\frac{1+i \sqrt 7}{2} \right ]e_1+\left [1,\frac{1+i \sqrt 7}{2} \right ]e_2
+\left [1,\frac{1+i \sqrt 7}{2} \right ]e_3\\[-12pt]
{}
\end{array}
$&$\begin{array}{@{}l@{}}
{}\\[-12pt]
e_1=\varepsilon _2, \\
e_2=(1-i \sqrt 7)(\varepsilon _2+\varepsilon _3)/4, \\
e_3=\left (-\varepsilon _1-\varepsilon _2+\frac{1+i \sqrt 7}{2}\,\varepsilon _3 \right )/2\\[-12pt]
{}
\end{array}
$&\raisebox{-5pt}[0pt][0pt]{$c=0$}
\\
\cline{1-1}\cline{3-5}
$[K_{25}]$&\raisebox{-20pt}[0pt][0pt]{3}&$K_{25}$&$[1,\omega ]e_1+[1,\omega ]e_2+[1,\omega ]e_3$&$\begin{array}{@{}l@{}}
{}\\[-12pt]
e_1=\varepsilon _3, \\
e_2=\frac{1-\omega }{3}\,(\varepsilon _1+\varepsilon _2+\varepsilon _3), \\
e_3=-\omega \varepsilon _2\\[10pt]
{}
\end{array}
$&\\
\cline{1-1}\cline{3-5}
$[K_{26}]_1\vphantom {\displaystyle \frac{\frac{1}{2}}{\frac{1}{2}}} $&&\raisebox{-14pt}[0pt][0pt]{$K_{26}$}&$[1,\omega ]e_1+[1,\omega ]e_2+[1,\omega ]\sqrt 2\,e_3$&\raisebox{-14pt}[0pt][0pt]{$\begin{array}{@{}l@{}}
{}\\[-12pt]
e_1=\frac{1-\omega ^2}{3}\,(\varepsilon _1+\varepsilon _2+\varepsilon _3), \\
e_2=\varepsilon _3, \\
e_3=\frac{1}{\sqrt 2}\,(\varepsilon _2-\varepsilon _3)\\[-12pt]
{}
\end{array}
$}&\\
\cline{1-1}\cline{4-4}
$[K_{26}]_2\vphantom {\displaystyle \frac {\frac {1}{2}}{\frac {1}{2}}} $&&&$[1,\omega ]e_1+p1,\omega ]e_2+[1,\omega ]i \sqrt {\frac {2}{3}}\,e_3$&&\\
\cline{1-5}
$[F_4]_1^\alpha $&\raisebox{-86pt}[0pt][0pt]{4}&$
$&$\begin{array}{@{}l@{}}
{}\\[-12pt]
[1,\alpha ]e_1+[1,\alpha ]e_2+[1,\alpha ]\sqrt 2\,e_3+[1,\alpha ]\sqrt 2\,e_4, \\
\alpha \in \Omega \\[-12pt]
{}
\end{array}
$&\raisebox{-34pt}[0pt][0pt]{$\begin{array}{@{}l@{}}
{}\\[-12pt]
e_j=(\varepsilon _{j+1}-\varepsilon _{j+2})/\sqrt 2, \\
j=1,2, \\
e_3=\varepsilon _4, \\
e_4=(\varepsilon _1-\varepsilon _2-\varepsilon _3-\varepsilon _4)/2\\[-12pt]
{}
\end{array}
$}&
\\
\cline{1-1}\cline{4-4}
$[F_4]_2^\beta $&&
$\begin{array}{@{}l@{}}
{}\\[-12pt]
K_{28},\\
\text{type}\;{\sf F}_4
\\[-12pt]
{}
\end{array}
$
&$\begin{array}{@{}l@{}}
{}\\[-12pt]
[1,\beta ]e_1+\left [1,\frac{\beta}{2} \right ]\!e_2+[1,\beta ]\sqrt 2e_3+
\left [1,\frac{\beta}{2} \right ]\!\sqrt 2e_4, \\
\beta \in \Omega \\[-12pt]
{}
\end{array}
$&&\\
\cline{1-1}\cline{4-4}
$[F_4]_3^\gamma $&&&$\begin{array}{@{}l@{}}
{}\\[-12pt]
[1,\gamma ]e_1+[1,\gamma ]e_2+\left [1,\frac{1+\gamma }{2} \right ]\sqrt 2\,e_3+ \\
+\left [1,\frac{1+\gamma }{2} \right ]\sqrt 2\,e_4,\quad \gamma \in \Omega \\[-12pt]
{}
\end{array}
$&&\\
\cline{1-1}\cline{3-5}
$[K_{29}]$&&$K_{29}$&$[1,i]e_1+[1,i]e_2+[1,i]e_3+[1,i]e_4$&$\begin{array}{@{}l@{}}
{}\\[-12pt]
e_1=\frac{1}{\sqrt 2}\,(\varepsilon _2-\varepsilon _4),\;\; 
e_2=\frac{1}{\sqrt 2}\,(-i \varepsilon _2+\varepsilon _3) \\
e_3=\frac{1}{\sqrt 2}\,(-\varepsilon _3+\varepsilon _4) \\
e_4=\frac{-1+i}{2 \sqrt 2}\,(\varepsilon _1+\varepsilon _2+\varepsilon _3+\varepsilon _4)\\[-12pt]
{}
\end{array}
$&\raisebox{105pt}[0pt][0pt]{$c=0$}\\
\cline{1-1}\cline{3-5}
$[K_{31}] \vphantom {\displaystyle \frac{1}{2}}$&&\raisebox{-30pt}[0pt][0pt]{$K_{31}$}&\raisebox{-30pt}[0pt][0pt]{$[1,i]e_1+[1,i]e_2+[1,i]e_3+[1,i]e_4$}&\raisebox{-31pt}[0pt][0pt]{$\begin{array}{@{}l@{}}
{}\\[2pt]
e_1=\frac{1}{\sqrt 2}\,(\varepsilon _2-\varepsilon _4),\;\; 
e_2=\frac{1}{\sqrt 2}\,(-i \varepsilon _2+\varepsilon _3) \\
e_3=\frac{1}{\sqrt 2}\,(-\varepsilon _3+\varepsilon _4) \\
e_4=\frac{-1+i}{2 \sqrt 2}\,(\varepsilon _1+\varepsilon _2+\varepsilon _3+\varepsilon _4) \\
e_5=\frac{1-i}{\sqrt 2}\,\varepsilon _4\\[-2pt]
{}
\end{array}
$}&\\
\cline{1-1}\cline{6-6}
$[K_{31}]^*$&
&&&&$\begin{array}{@{}l@{}}
{}\\[-12pt]
c(r_j)=0\\
j=1,2,3,4, \\
c(r_5)\\= 
\frac{1+i}{2}\,e_5\\[-12pt]
{}
\end{array}
$\\
\cline{1-1}\cline{3-6}
$[K_{32}]$&\raisebox{0pt}[0pt][0pt]{4}&$K_{32}$&$[1,\omega ]e_1+[1,\omega ]e_2+[1,\omega ]e_3+[1,\omega ]e_4$&$\begin{array}{@{}l@{}}
{}\\[-12pt]
e_1=\varepsilon _3, \\
e_2=\frac{1-\omega }{3}\,(\varepsilon _1+\varepsilon _2+\varepsilon _3), \\
e_3=-\omega \varepsilon _2, \\
e_4=\frac{\omega ^2-\omega }{3}\,(-\varepsilon _1+\varepsilon _2+\varepsilon _3)\\[12pt]
{}
\end{array}
$&\raisebox{220pt}[0pt][0pt]{$c=0$}\\
\cline{1-5}
$[K_{33}]\vphantom {\begin{array}{@{}l@{}}
{}\\[12pt]
1\\
2\\[-12pt]
\end{array}
}$&5&$K_{33}$&\raisebox{-22pt}[0pt][0pt]{$[1,\omega ]e_1+\cdots +[1,\omega ]e_n$}&\raisebox{-22pt}[0pt][0pt]{$\begin{array}{@{}l@{}}
{}\\[-12pt]
e_1=\frac{\omega}{\sqrt 2}\,(\varepsilon _5+\varepsilon _6), \\
e_2=-\frac{\omega}{2 \sqrt 2}\,(-\varepsilon _1+(1+2 \omega )\varepsilon _2 \\
\hskip 26.5mm+ \varepsilon _3+\varepsilon _4+\varepsilon _5+\varepsilon _6), \\
e_j=\frac{1}{\sqrt 2}\,(\varepsilon _{j-2}-\varepsilon _{j-1}), \\
j=3,4,\ldots ,n\\[-12pt]
{}
\end{array}
$}&\raisebox{-30pt}[0pt][0pt]{$c=0$}\\
\cline{1-3}
$[K_{34}]\vphantom {\begin{array}{@{}l@{}}
{}\\[12pt]
1\\
2\\[-12pt]
\end{array}
}$&\raisebox{-16pt}[0pt][0pt]{6}&$K_{34}$&&&\\
\cline{1-1}\cline{3-5}
$[{\sf E}_6]^\alpha $&&$\begin{array}{@{}l@{}}
{}\\[-12pt]
K_{35},\\
\text{type}\;{\sf E}_6
\\[-12pt]
{}
\end{array}
$&\raisebox{-30pt}[0pt][0pt]{$\begin{array}{@{}l@{}}
{}\\[-12pt]
[1,\alpha ]e_1+\cdots +[1,\alpha ]e_n, \\
\alpha \in \Omega \\[-12pt]
{}
\end{array}
$}&\raisebox{-33pt}[0pt][0pt]{$\begin{array}{@{}l@{}}
{}\\[-15pt]
e_1=(\varepsilon _1-\varepsilon _2-\varepsilon _3-\varepsilon _4-\varepsilon _5 \\
\hskip 13.5mm -\varepsilon _6-\varepsilon _7+\varepsilon _8)/2 \sqrt 2, \\
e_2=(\varepsilon _1+\varepsilon _2)/\sqrt 2, \\
e_j=(-\varepsilon _{j-2}+\varepsilon _{j-1})/\sqrt 2, \\
j=3,\ldots ,n.\\[-12pt]
{}
\end{array}
$}&\\
\cline{1-3}
$[{\sf E}_7]^\alpha $&7&$\begin{array}{@{}l@{}}
{}\\[-12pt]
K_{36}, \\
\text{type}\;{\sf
E}_7\\[-12pt]
{}
\end{array}
$&&&\\
\cline{1-3}
$[{\sf E}_8]^\alpha $&8&$\begin{array}{@{}l@{}}
{}\\[-12pt]
K_{37}, \\
\text{type}\;{\sf
E}_8\\[-12pt]
{}
\end{array}
$&&&\\
\end{longtable}
\end{landscape}



\subsection{Equivalence}

\vskip 2mm

\subsubsection{Theorem {\rm(Equivalent groups from Table~\ref{tab2})}}\label{T??}
{\it
Table
$3$
below gives the complete list of equivalent groups $W$ and $W'$ from Table~{\rm\ref{tab2}} such that $W \ne W'$.
}


\begin{table}[h!]\centering
\addcontentsline{toc}{subsubsection}{Table 3.
Pairs of equivalent irreducible infinite complex crystallographic
$r$-groups}\label{ttt3}

\vskip 2mm

\begin{tabular}{|l|l|c|}
\multicolumn{3}{l}{\textbf{Table 3.}
Pairs of equivalent irreducible infinite complex crystallographic
$r$-groups}\label{tab3}\\[3mm]
\hline
\multicolumn{3}{|c|}{}\\[-14pt]
\multicolumn{1}{|c|}{$W$}&\multicolumn{1}{c|}{$W'$}&condition\\
\hline
&&\\[-1pt]
$[G(2,1,s)]_2^{1+\omega },\quad s \geqslant 3$&$[G(2,1,s)]_3^{1+\omega },\quad s \geqslant 3$&---\\[11pt]
$[G(2,1,s)]_2^{1+\omega },\quad s \geqslant 3$&$[G(2,1,s)]_4^{1+\omega },\quad s \geqslant 3$&---\\[11pt]
$[G(2,1,s)]_3^i,\quad s \geqslant 3$&$[G(2,1,s)]_4^i,\quad s \geqslant 3$&---\\[11pt]
$[G(2,1,2)]_2^\beta $&$[G(2,1,2)]_2^{-2/\beta }$&$-2/\beta \in \Omega $\\[11pt]
$[G(2,1,2)]_2^{1+\omega }$&$[G(2,1,2)]_3^{1+\omega }$&---\\[11pt]
$[G(2,1,2)]_2^\beta $&$[G(2,1,2)]_3^{1-2/\beta }$&$1-2/\beta \in \Omega $\\[11pt]
$[G(2,1,2)]_2^\gamma $&$[G(2,1,2)]_3^{(\gamma -1)/(\gamma +1)}$&$(\gamma -1)/(\gamma +1) \in \Omega $\\[11pt]
$[G(2,1,2)]_2^\beta $&$[G(2,1,2)]_3^{-1-2/\beta}$&$-1-2/\beta \in \Omega $\\[11pt]
$[G(6,6,2)]_2^\beta $&$[G(6,6,2)]_2^{-3/\beta }$&$-3/\beta \in \Omega $\\[11pt]
$[G(6,6,2)]_3^\gamma $&$[G(6,6,2)]^{(2 \gamma -1)/(\gamma +1)}$&$(2 \gamma -1)/(\gamma +1) \in \Omega $\\[11pt]
$[G(6,6,2)]_2^\beta $&$[G(6,6,2)]^{-1+3/\beta }$&$-1+3/\beta \in \Omega $\\[11pt]
$[G(6,6,2)]_2^\beta $&$[G(6,6,2)]_3^{2-3/\beta }$&$2-3/\beta \in \Omega $\\[11pt]
$[{\sf F}_4]_2^\beta $&$[{\sf F}_4]_2^{-2/\beta }$&$-2/\beta \in \Omega $\\[11pt]
$[{\sf F}_4]_2^{1+\omega }$&$[{\sf F}_4]_3^{1+\omega }$&---\\[11pt]
$[{\sf F}_4]_2^\beta $&$[{\sf F}_4]_3^{1-2/\beta }$&$1-2/\beta \in \Omega $\\[11pt]
$[{\sf F}_4]_3^\gamma $&$[{\sf F}_4]_3^{(\gamma -1)/(\gamma +1)}$&$(\gamma -1)/(\gamma +1) \in \Omega $\\[11pt]
$[{\sf F}_4]_2^\beta $&$[{\sf F}_4]_3^{-1-2/\beta }$&$-1-2/\beta \in \Omega $\\[13pt]
\hline
\end{tabular}
\end{table}

\eject

The rather technical proof of this theorem will not be given here.

\subsection{The structure of an extension of $\Tran W$ by $\Lin W$}

As we have seen in Section
\ref{1.5}, if $\ka =\Ree $, then the structure of an infinite irreducible $r$-group $W$ as an extension of $\Tran W$ by $\Lin W$ is very simple: it is always a~semidirect product. The situation is more complicated when $\ka = \Cee $, because there exist infinite complex irreducible
crystallographic $r$-groups $W$ which are \textit{not} semidirect products of $\Tran W$ and $\Lin W$.

\subsubsection{Theorem}\label{th1}
{\it
The groups $W$ from Table~$\ref{tab2}$ which are not semidirect products of $\Tran W$ and $\Lin W$ are
$$
[G(4,2,s)]_1^*,\quad [K_{12}]^* \quad \text{and}\quad [K_{31}]^*.
$$
}
\subsubsection{Theorem}\label{th2}
{\it
Let $K \subseteq \GL(V)$ be a~finite ir\-re\-ducible $r$-group and let $T \subseteq V$ be a~$K$-invariant lattice. Assume that there exists a~cry\-stal\-lographic $r$-group $W$ with $\Lin W = K$, $\Tran W = T$. Then the set of
elements of $H^1 (K,V/T)$
corresponding to such subgroups $W$ is, in fact, a~subgroup of $H^1 (K,V/T)$ and the order of this subgroup is $\leqslant 2$.
}

\subsection{The rings and fields of definition of \boldmath$\Lin W$}

As we have seen in Section \ref{1.5}, if $\ka = \Ree $, then the group $\Lin W$ for an infinite irreducible $r$-group $W$ is defined over $\Qee $. If $\ka =\Cee $, then $\Lin W $ for an infinite irreducible crystallographic $r$-group $W$ is defined over a~certain purely imaginary quadratic extension of $\Qee $, see Theorem\;\ref{cirg2}.
We describe this extension explicitly in the following Theorem \ref{cirg3}.

\subsubsection{Theorem}\label{cirg3}
{\it
Let $K \subseteq \GL(V)$ be a~finite complex irreducible
$r$-group.\;Then the ring with unity generated over $\Zee  $ by the set of all cyclic products related to an arbitrary generating system of reflections of $K$ coincides with the ring $\Zee  [\Tr K]$ generated over $\Zee  $ by the set of traces of all elements of $K$. The ring $\Zee  [\Tr K]$ is the minimal ring of definition of~$K$. This ring is isomorphic to $\Zee  $ if and only if $K$ is the complexification of the Weyl group of an irreducible root system.
}
\vskip 2mm

\noindent {\it Proof} is given in Sections~\ref{4.6} and \ref{fis}.

\vskip 2mm

It is easily seen from Table~1
and Theorem \ref{cirg3} that for the groups $K = \Lin W$, where $W$ is an infinite irreducible crystallographic $r$-group, one has
Table $4$ below.

\begin{table}[!h]\centering
\footnotesize
\addcontentsline{toc}{subsubsection}{Table 4.
Linear parts of
infinite complex irreducible crystallographic
$r$-groups}
\begin{tabular}{|c||l|l|@{\extracolsep{.5mm}}l@{\extracolsep{.5mm}}|@{\extracolsep{.5mm}}c@{\extracolsep{.5mm}}|l@{\extracolsep{.5mm}}|@{\extracolsep{.5mm}}c@{\extracolsep{.5mm}}|@{\extracolsep{.5mm}}c@{\extracolsep{.5mm}}|}
\multicolumn{8}{l}{\textbf{Table 4.}
Linear parts of infinite complex irreducible
crystallographic
$r$-groups}\label{tab4}\\[3mm]
\hline
$\Zee  [\Tr K]$&\multicolumn{1}{c|}{$\Zee  $}&\multicolumn{1}{c|}{$\Zee  [i]$}&\multicolumn{1}{c|}{$\Zee  [2i]$}&$\Zee  [i \sqrt 2]$&\multicolumn{1}{c|}{$\Zee  [\omega ]$}&$\Zee  [2 \omega ]$&$\Zee  \left [\frac {1+i \sqrt 7}{2} \right ]_{\mathstrut }^{\mathstrut }$\\
\hline
&&&&&&&\\[-10pt]
\hline
&&&&&&&\\[-6pt]
$K$
&
$\begin{array}[t]{@{}l@{}}
K_1={\sf A}_s,\\[-1.5mm]
s \geqslant 1;
\end{array}
$
&
$\begin{array}[t]{@{}l@{}}
G(4,1,s),\\[-1.5mm]
s \geqslant 2;
\end{array}
$
&
$G(4,2,2)$
&
$K_{12}$&$\begin{array}[t]{@{}l@{}}
G(3,1,s),\\[-1.5mm]
s \geqslant 2;
\end{array}
$
&
$G(6,3,2)$
&
$K_{24}$\\[6.5mm]
&
$\begin{array}[t]{@{}l@{}}
G(2,1,s)={\sf B}_s,\\[-1.5mm]
s \geqslant 2;
\end{array}
$
&
$\begin{array}[t]{@{}l@{}}
G(4,4,s),\\[-1.5mm]
s \geqslant 3;
\end{array}
$
&
&
&
$\begin{array}[t]{@{}l@{}}
G(6,1,s),\\[-1.5mm]
s \geqslant 2;
\end{array}
$
&
&\\[6.5mm]
&
$\begin{array}[t]{@{}l@{}}
G(2,2,s)={\sf D}_s,\\[-1.5mm]
s \geqslant 3;
\end{array}
$
&
$\begin{array}[t]{@{}l@{}}
G(4,2,s),\\[-1.5mm]
s \geqslant 3;
\end{array}
$
&
&
&
$\begin{array}[t]{@{}l@{}}
G(3,3,s),\\[-1.5mm]
s \geqslant 3;
\end{array}
$
&
&\\[6.5mm]
&
$G(6,6,2)={\sf G}_2;$
&
$\begin{array}[t]{@{}l@{}}
K_3\;
\text{($m=4$);}
\end{array}
$
&
&
&
$\begin{array}[t]{@{}l@{}}
G(6,6,s),\\[-1.5mm]
s \geqslant 3;
\end{array}
$
&
&\\[6.5mm]
&
$K_{28}={\sf F}_4;$
&
$K_8$;
&
&
&
$\begin{array}[t]{@{}l@{}}
G(6,2,s),\\[-1.5mm]
s \geqslant 2;
\end{array}
$
&
&\\[6.5mm]
&
$K_{35}={\sf E}_6$;
&
$K_{29}$;
&
&
&
$\begin{array}[t]{@{}l@{}}
G(6,3,s),\\[-1.5mm]
s \geqslant 3;
\end{array}
$
&
&\\[6.5mm]
&
$K_{36}={\sf E}_7$;
&
$K_{31}$
&
&
&
$\begin{array}[t]{@{}l@{}}
K_3
\text{($m=3,6$);}
\end{array}
$
&
&\\[4mm]
&
$K_{37}={\sf E}_8$
&
&
&
&
$K_4$;
&
&\\[4mm]
&
&
&
&
&
$K_5$;
&
&\\[4mm]
&
&
&
&
&
$K_{25}$;
&
&\\[4mm]
&
&
&
&
&
$K_{26}$;
&
&\\[4mm]
&
&
&
&
&
$K_{32}$;
&
&\\[4mm]
&
&
&
&
&
$K_{33}$;
&
&\\[4mm]
&
&
&
&
&
$K_{34}$
&
&\\[4mm]
\hline
&
&
&
&
&
&
&\\[-9pt]
$\begin{array}{@{}l@{}}
\text{field of} \\[-1.5mm]
\text{fractions} \\[-1.5mm]
\text{of}\; \Zee  [\Tr K]
\end{array}$&\multicolumn{1}{c|}{$\Qee $}&$\Qee (\sqrt {-1})$&$\Qee (\sqrt {-1})$&$\Qee (\sqrt {-2})$&$\Qee (\sqrt {-3})$&$\Qee (\sqrt {-3})$&$\Qee (\sqrt {-7})$\\[3pt]
\hline
\end{tabular}
\end{table}

\subsection{Further remarks} $\phantom 1$\label{fre}

a) In contrast to the real case, there exist \textit{$1$-parameter families} of inequivalent irreducible complex infinite crystallographic $r$-groups $W$ with a~fixed linear part $\Lin W$ (i.e., the groups with a~fixed linear part may have \textit{moduli}). We will see below that \textit{an irreducible crystallographic $r$-groups $W$ with $\Lin W = K$ has moduli if and only if $\Zee  [ \Tr K] = \Zee  $, i.e., if and only if $K$ is the complexification of the Weyl group of an irreducible root system}.

\vskip 1mm

b) It follows from Table~4
(and from the~known result of algebraic number theory) that the ring $\Zee  [ \Tr \Lin W] $, where $W$ is an infinite irreducible crystallographic $r$-group, is always a
principal ideal
domain.
It would be interesting to have
an a~priori proof of this fact.

\vskip 1mm

c) If $k=\Ree $, then it is known (and was \textit{a priori} proved in 1948--51 by Chevalley and Harish--Chandra) that assigning to every complex simply connected semisimple 
algebraic group its affine Weyl group yields a bijection between the set of isomorphism classes
of
complex simply connected semisimple
algebraic groups and the set of
equivalence classes of
all real crystallographic
$r$-groups.


\vskip 2mm

\noindent \textbf{Question.} 
Is there a \lq \lq global object\rq \rq (a counterpart of complex simply connected
semi\-simple algebraic group) naturally assigned to every complex crystallographic $r$-group, so that this assign\-ment yields a bijection between all equivalence classes of such $r$-groups and all isomorphism classes of
\lq \lq global objects\rq \rq?
\vskip 3mm
$$
\begin{array}{ccc}
\cline{1-1}\cline{3-3}
\multicolumn{1}{|c|}{\begin{array}{c}
\text{real}\\[0mm]
\text{crystallographic}\\[0mm]
\text{$r$-group}
\end{array}
}&\longleftrightarrow &\multicolumn{1}{|c|}{\begin{array}{c}
\text{complex simply connected}\\[0mm]
\text{semisimple algebraic groups}
\end{array}
}\\
\cline{1-1}\cline{3-3}
{}\\[0pt]
\cline{1-1}\cline{3-3}
\multicolumn{1}{|c|}{\begin{array}{c}
\text{complex}\\[0mm]
\text{crystallographic}\\[0mm]
\text{$r$-group}
\end{array}
}&\longleftrightarrow &\multicolumn{1}{|c|}{\begin{array}{c}
\text{?}
\end{array}
}\\
\cline{1-1}\cline{3-3}
\end{array}
$$
\vskip 5mm
\noindent
I do not know whether such a \lq\lq global object\rq\rq exists or not.
Curiously
that we can calculate (see below Remark \ref{corgg})
the group which, by analogy with the case of $k=\bf R$, might be  considered as the \lq \lq center\rq \rq {} of this hypothetical object.

\section{\bf Several auxiliary results and the classification of infinite complex irreducible noncrystallo\-graphic
 $r$-groups}\label{S3}

Now we move on to proofs. The greater part of these proofs relies heavily on the fact that for an infinite $r$-group $W$ the subgroup $\Tran W$ is sufficiently \lq \lq massive\rq \rq .

\subsection{The subgroup of translations}\label{3.1}

In order to clarify the last statement, we need the following
basic fact.

\subsubsection{Theorem {\rm (Existence of nonzero translations)}}\label{nontr}
{\it Let $W \subseteq A(E)$ be an infinite $r$-group. Then $\Tran W \ne 0$.
}

\begin{proof} Assume that $\Tran W = 0$. Then $\Lin\colon W \tto {\rm Iso}(V)$ is a~monomorphism.

 From the discreteness of $W$ it follows that $(\overline {\Lin W})^0$ is a~torus (here bar denotes the closure and the superscript
$0$ singles out the connected component of unity), see \cite[Chap.\,III, \S4, Exer.\,13a]{5}.
Set
$$
S = \Lin W \cap (\overline {\Lin W})^0.
$$
Then  $\Lin W \triangleright S$ and $[\Lin W\colon S] <\infty $. Let also
$$
V_0 := \{v\in V \mid (\overline {\Lin W})^0v = v\}
$$
and let $V_1$ be the subspace of $V$ determined from the decomposition
$$
V=V_0 \oplus V_1,\quad \text{where}\quad v_0 \perp V_1.
$$
The subspaces $V_0$ and $V_1$ are $S$-invariant. Definitely, $V_1 \ne 0$ (indeed, if not, then $S = 1$, and hence $| \Lin W| <\infty $
contrary to $| \Lin W| = |W|$).

The set
$$
\left \{P\in (\overline {\Lin W})^0 \mid H_P = V_0 \right \}
$$
(see notation in Section \ref{1.1}) is dense in $(\overline {\Lin W})^0$.
Therefore,
$\{P\in S \mid H_P= V_0\}$ is dense in $(\overline {\Lin W})^0$.
As $S$ is dense in $(\overline {\Lin W}) ^0$, it follows that
$$
\{P \in S \mid H_P = V_0\} \ne 0,
$$

\textit{We claim that there exists a~point $a\in E$ with $\gamma (a)- a~\in V_0$ for every $\gamma \in W$}.

\vskip 1mm

To prove this, consider an arbitrary point $b \in E$ and an operator $P_0\in S$ such that $H_{P_0} = V_0$. Take $\gamma _0\in W$ with $\Lin \gamma _0 = P_0$. We have $\gamma _0(b )-b=t_0+t_1$, where $t_0\in V_0$, $t_1\in V_1$. But the restriction of $\ii-P_0$ to $V_1$ is nondegenerate. Hence, there exists a~vector $t\in V_1$ such that $(\ii -P_0)t = t_1$. So, we have
$$
\gamma _0(b+t)-(b+t) = (\gamma _0(b)+P_0t)-(b+t) = (\gamma _0(b)-b) + (P_0-\ii)t = t_0+t_1-t_1 = t_0\in V_0.
$$
Put $a = b+t$. We have $\gamma _0(a)-a\in V_0$. Let us prove that $a$ is a~point wanted.

Let us first check that
$$
\gamma (a) -a\in V_0 \quad \text{if}\quad \gamma \in W \quad \text{and}\quad \Lin \gamma \in S.
$$
Write
$$
P = \Lin \gamma ,\quad v_0 =\gamma _0 (a)-a\in V_0 \quad \text{and}\quad v =\gamma (a)-a.
$$
We need to prove that $v\in V_0$. But $S$ is commutative, so $PP_0 = P_0P$. Hence, $\gamma \gamma _0 =\gamma _0 \gamma $. Now we have
\begin{align*}
\kappa _a(\gamma _0)&=(P_0,v_0)=\kappa _a(\gamma \gamma _0 \gamma ^{-1})=(P,v)(P_0,v_0)(P^{-1},-P^{-1}v)\\
&=(PP_0P^{-1},-PP_0P^{-1}v+Pv_0+v)=(P_0,-P_0v+Pv_0+v).
\end{align*}
So, $-P_0v+Pv_0+v = v_0$, i.e., $(\ii-P_0)v = (\ii-P)v_0$. But $P\in S$, hence $(\ii-P)v_0 = 0$. From $\Ker(\ii-P_0) = V_0$ it follows that $v\in V_0$. This establishes the claim.



Now we can prove that $\gamma (a)-a\in V_0$ for arbitrary $\gamma \in W$. Indeed, write $P = \Lin \gamma $ and $v =\gamma (a)-a$, as before. We have:
$$
\kappa _a(\gamma \gamma _0 \gamma ^{-1})=(PP_0P^{-1},-PP_0P^{-1}v+Pv_0+v).
$$

From $S \triangleleft \Lin W$ it follows that $V_0$ is $P$-invariant. Therefore, we have $Pv_0\in V_0$. Moreover, if $\gamma '=\gamma \gamma _0 \gamma ^{-1}$, then $\Lin \gamma' \in S$ and, by the claim,
$$
\gamma '(a)-a=-PP_0P^{-1}v+Pv_0+v\in V_0.
$$
Therefore, $-PP_0P^{-1}v+v \in V_0$, and hence
$$
-P_0P^{-1}v+P^{-1}v\in P^{-1}V_0=V_0,
$$
i.e., $(\ii-P_0)P^{-1}v\in V_0$. But the image of $\ii-P_0$ is $V_1$, hence $(\ii-P_0)P^{-1}v=0$. Therefore, $P^{-1}v\in V_0$, i.e., $v\in PV_0 = V_0$ and we are done.

\vskip 2mm

\textit{Now, we consider two subgroups of $W${\rm:} let $W'$, resp. $W''$, be the subgroup generated by those reflections $\gamma $ for which $H_\gamma \ni a$,
resp. $H_\gamma \not \ni a$}.

\vskip 2mm

The subgroup $W'$ is \textit{finite}. Indeed, identifying $A(E)$ and $\GL(V)\ltimes V$ by means of $\kappa _a$, we have $\kappa _a(\gamma ) = (R,0) \in {\rm Iso}(V)$ for each reflection $\gamma \in W'$, and hence for each $\gamma \in W '$. So $\kappa _a(W ')$ is a~discrete subgroup of a~compact group ${\rm Iso}(V)$, hence finite.

\vskip 2mm

We claim now that $W''$ is \textit{infinite}. Before proving this,
we shall show how to complete
the proof of the theorem if the claim holds. Thus, assume $W''$ is an infinite $r$-group.

Note that $\Tran W'' = 0$. Using the above arguments and constructions we obtain a~torus $(\overline {\Lin W''} ) ^0$, a~subgroup $S'' = \Lin W'' \cap (\overline {\Lin W''})^0$, and a~decomposition $V = V_0'' \oplus V_1''$, where
$$
V_0''=\left \{v\in V \mid (\overline {\Lin W''})^0v = v \right \}\quad \text{and}\quad V_0'' \perp V_1''.
$$
We also have $V_1'' \ne 0$ and $\{P\in S'' \mid H_P=V_0''\}\ne \varnothing $. But $W'' \subseteq W$, therefore, $\Lin W'' \subseteq \Lin W$, whence $(\overline {\Lin W''})^0 \subseteq (\overline {\Lin W})^0$. So $S'' \subseteq S$ and $V_0'' \supset V_0$. It follows now that $\gamma (a)-a\in V_0''$ and $\gamma (a)-a \ne 0$ for every $\gamma \in W''$.

Let $\gamma \in W''$ be a~reflection. Then
$$
\kappa _a(\gamma )= (\Lin \gamma ,\gamma (a)-a),
$$
and in view of the fact that $\gamma (a)-a \ne 0 $, we have:
$$
(\text{root line of $\Lin \gamma $)} =
\ka (\gamma (a)-a)\subseteq V_0''.
$$
Therefore, $H_{\Lin \gamma }\supset V_1$ and $\Lin \gamma $ trivially acts on $V_1$. But $W''$ is generated by reflections; therefore $\Lin W''$ trivially acts on $V_1$ which contradicts the existence of $P\in S'' \subseteq \Lin W''$ such that $H_P=V_0''$ and $V_1'' \ne 0$.

It remains to show that $|W''| = \infty $. Let us assume that this is not the case. Then  $W''$ has a~fixed point $b\in E$. By construction, $b \ne a$.
We shall show that $W'$ and $W''$ commute. This then yields that $W = W'W''$ is a~finite group (as, clearly, $W'$ and $W''$
generate~$W$), which is a~contradiction.

\eject

Let $\gamma ' \in W'$ and $\gamma ''\in W''$ be reflections. We have $a\in H_{\gamma '}$, $b\in H_{\gamma ''}$.

\begin{figure}[h!]\centering
\includegraphics
{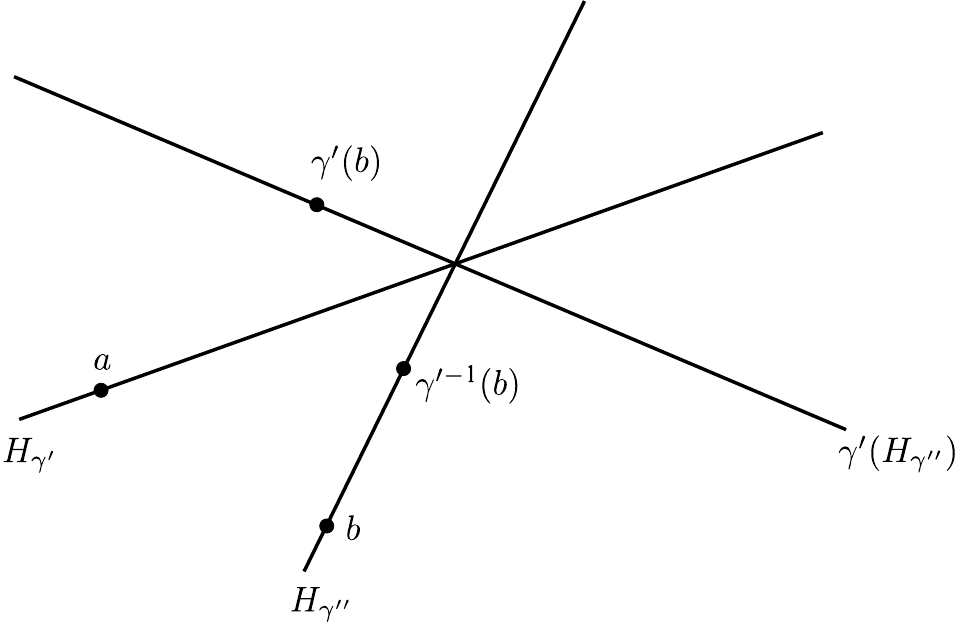}
\end{figure}

But $\gamma ' \gamma '' \gamma ^{\prime -1}$ is a~reflection with mirror $\gamma '(H_{\gamma ''})$. Hence, $\gamma ' \gamma '' \gamma ^{\prime -1}$ is either in $W'$ or in $W''$. If this element is in $W'$, then $\gamma '(H_{\gamma ''})\ni a$, i.e., $H_{\gamma ''}\ni \gamma ^{\prime -1}(a)=a$, which is a~contradiction. Therefore, $\gamma '(H_{\gamma ''})\ni b$ and $H_{\gamma ''}\ni \gamma ^{\prime -1}(b)$. We also have $b\in H_{\gamma ''}$ and $b \ne \gamma ^{\prime -1}(b)$ (for otherwise $b\in H_{\gamma '}$ which is absurd).

Consequently, there is a~unique line through $b$ and $\gamma ^{\prime -1}(b)$. This line lies in $H_{\gamma '}$ and is orthogonal to $H_{\gamma ^{\prime -1}}=H_{\gamma '}$. Hence, $H_{\gamma ''}\perp H_{\gamma '}$, i.e., $\gamma ' \gamma ''=\gamma '' \gamma '$.
\end{proof}

Theorem \ref{nontr} implies that $\Tran W$ is ``big enough'':

\subsubsection{Theorem}\label{n2n}
{\it Let $W$ be an infinite irreducible $r$-group. Then  $T = \Tran W$ is a~lattice of rank
$n = \dim\nolimits_{\ka} E$ if $\ka = \Ree $, and of rank $n$ or $2n$ if $\ka = \Cee $.
}

\begin{proof} Let $\ka = \Ree $. Then  $\Ree T$ is an invariant subspace of the $\Lin W$-module $V$. This subspace is nontrivial by
Theorem\;\ref{nontr}. Therefore, $\Ree T = V$ because of irreducibility (see Theorem~\ref{irred}),
and the assertion follows from the equality $\rk T = \dim \Ree T$.

Let $\ka = \Cee $. As above, we have $0 \ne \Cee T=\Ree T+i \Ree T=V$. Hence
$$
2n = \dim\limits_{\Ree }V \leqslant \dim \limits_{\Ree }\Ree T+\dim \limits_{\Ree }i \Ree T=2\rk T.
$$
Therefore, $\rk T \geqslant n$. But $\Ree T \cap i \Ree T$ is a $\Lin W$-invariant complex subspace of $V$, and hence either $\Ree T \cap i \Ree T=0$, or $\Ree T \cap i \Ree T=V$, i.e., $\Ree T=V$. If $\rk T = \dim \Ree T>n$, then $\Ree T \cap i \Ree T \ne 0$, hence $\Ree T=V$, i.e., $\rk T = 2n$.
\end{proof}

\subsubsection{Corollary}\label{coro}

{\it
\begin{enumerate}[\hskip 7.2mm\rm 1)]\itemsep=-.1ex
\item
$|\Lin W|<\infty $.
\item
If $\ka = \Cee $, then $W$ is a~crystallographic group if and only if $\rk T = 2n$.
\end{enumerate}
}

\begin{proof} 1) follows from the fact that $\Lin W$ is contained in a~compact group and
in $V$ there is a lattice of maximal rank and $\Lin W$-invariant.
\end{proof}

\subsection{Some auxiliary results}\label{auxi}

Let $K \subseteq \GL(V)$ be \textit{a~finite $r$-group and let $\mathcal H $ be the set of mirrors of all reflections from} $K$. Let $H\in \mathcal H $.

Then  it is easy to see that the subgroup of $K$ generated by all reflections $R\in K$ with $H_R = H$, \textit{is a~cyclic group}. Let $m(H)$ be \textit{the order of this group}.

\subsubsection{Theorem}\label{ctpgr}
{\it Let $\{R_j\}_{j\in J}$ be a~generating system of reflections of $K$ such that the order of $R_j$ is equal to $m(H_{R_j})$ for every $j\in J$. Then  each reflection $R\in K$ is conjugate in $K$ to $R_j^{l _j}$ for certain $j$ and $l _j$.
}

\begin{proof} Let $\mathcal O $ be the $K$-orbit of $H_R$ in $\mathcal H $ (since $PH_R=H_{PRP^{-1}}$, it follows that $K$ acts on $\mathcal H $).
The product of linear equations of all mirrors from $\mathcal O $ is a semi-invariant of $K$.
Let $\chi _{\mathcal O }$ be its character.
The latter is a homomorphism of $K$ into the multiplicative group of $\ka$. Therefore,
$\chi _{\mathcal O }(K)$ is a group
generated by $\chi _{\mathcal O }(R_j)$ for $j\in J$. We have also $\chi _{\mathcal O }(R)\ne 1$. Therefore, there exists a~$j\in J$ with $\chi _{\mathcal O }(R_j)\ne 1$. But $\chi _{\mathcal O } (R_j)\ne 1$ if and only if $\mathcal O $ is the orbit of $H_{R_j}$, see [6]. Therefore, $gH_R = H_{R_j}$ for a~certain $g\in K$ and the statement follows.
\end{proof}

\subsubsection{Theorem {\rm(Linear parts of reflections)}}\label{rR}
{\it Let $W$ be an $r$-group. Then  for any reflection $R\in \Lin W $, there exists a~ref\-lec\-tion $\gamma \in W$ with $\Lin \gamma = R$.
}

\begin{proof}
Let $\{\rho _j\}_{j\in J}$ be the set of all reflections of $W$. Then  $\{R_j= \Lin \rho _j\}_{j\in J}$ is a
system of reflections that generates $\Lin W$. Let $\mathcal O $ be the $\Lin W$-orbit of $H_R$ in $\mathcal H $. Then  there exists a~number $l $ with $H_{R_l } \in \mathcal O $, and an element $P_l \in \Lin W$ such that
$$
P_l H_{R_l }=H_{P_l R_l P_l ^{-1}}=H_R.
$$
Let $\pi _l \in W$ be such that $\Lin \pi _l =P_l $. We have:
$$
\Lin \pi _l \rho _l \pi _l ^{-1}=P_l R_l P_l ^{-1},
$$
and
$\pi _l \rho _l \pi _l ^{-1}$ is a~reflection, too. Therefore, $\chi _{\mathcal O }(\Lin W)$ is generated by $\chi _{\mathcal O }(R_j)$, where $j\in J'=\{j\in J \mid H_{R_j}=H_R\}$. We have
$$
\chi _{\mathcal O }(R)=\chi _{\mathcal O }(R_{j_1})\cdots \chi _{\mathcal O }(R_{j_s})=\chi _{\mathcal O }(R_{j_1}\cdots R_{j_s})
$$
for certain $j_1,\ldots ,j_s\in J'$. But $R_{j_1},\ldots ,R_{j_s}$ are the reflections whose mirrors are $H_{R}$. It follows now (see [6]) that $R = R_{j_1}\cdots R_{j_s}$. Let us consider the element $\rho =\rho _{j_1}\cdots \rho _{j_s}$. We have $\Lin \rho= R$ and it easily follows from the fact that the mirrors $H_{\rho _{j_1}},\ldots ,H_{\rho _{j_s}}$ are parallel to each other that $\rho $ is a~reflection.
\end{proof}

\subsection{Semidirect products}

Let $W$ be a~subgroup of $A(E)$. We have
$$
0 \tto \Tran W \hookrightarrow W \tto \Lin W \tto 1.
$$
When is $W$ a~semidirect product of $\Lin W$ and $\Tran W$? We have the following criterion.

\subsubsection{Theorem}\label{sdp}
{\it Let $|\Lin W|<\infty $. Then
\begin{enumerate}[\hskip 7.2mm\rm a)]\itemsep=-.1ex
\item  $W$ is a~semidirect product of $\Lin W$ and $\Tran W$ if and only if there exists a~point $a\in E$ such that $\Lin $ induces an isomorphism
between
the stabilizer $W_a$ of $a$ and the group
$\Lin W$.
\item  For every finite group $K \subseteq {\rm Iso}(V)$ and every $K$-invariant subgroup $T$ of $V$, there exists a~unique {\rm(}up to equivalence{\rm)} group $W \subseteq A(E)$ such that $\Lin W = K$, $\Tran W = T$ and $W$ is a~semidirect product of $\Lin W$ and $\Tran W$.
    \end{enumerate}
}

\textit{Proof} is left to the reader.\hfill \textsquare

\vskip 2mm

The point $a$ from part a) of  Theorem \ref{sdp}
is called a~{\it special point}\index{special point} of $W'$, see [1].

Using this theorem, we can clarify the structure of infinite $r$-groups in a~number of important cases:

\subsubsection{Theorem}\label{sdsd}
{\it Let $W \subseteq A(E)$ be a~group generated by reflections. Assume that ${| \Lin W|<\infty} $ and that $\Lin W$ is an essential group \textup{(}i.e., $\{v\in V \mid (\Lin W)v = v\} = \{0\})$ generated by $n = \dim \nolimits_{\ka} E$ reflections. Then  $W$ is a~semidirect product of $\Lin W$ and $\Tran W$.
}

\begin{proof}
Let $R_1,\ldots , R_n$ be a
system of reflections generating
$\Lin W$. Then  there exists a~reflection $\gamma _j\in W$ such that $\Lin \gamma _j=R_j$, where $j= l,\ldots , n$  (see the proof of Theorem \ref{rR}).
We have $\bigcap _{j=1}^nH_{R_j}=0$ because $\Lin W$ is essential. Therefore, $\bigcap_{j=1}^nH_{\gamma _j}\ne \varnothing $; more precisely, this intersection is a~single point $a\in E$.

We see now that $\Lin\colon W_a \tto \Lin W$ is a~surjective map, hence an isomorphism. Thus, we are done in view of Theorem~\ref{sdp}.
\end{proof}

The conditions of this theorem
always hold
if $\ka =\Ree $ and $W$ is an infinite irreducible $r$-group. If $\ka = \Cee $, in general this is not the case, though \lq \lq in most cases\rq \rq , it is.

\subsection{Classification of infinite complex irreducible
noncrystallographic $r$-groups}\label{3.4}

Let $\ka = \Cee $ and let $W$ be a~group as in the title of this section. Set $T = \Tran W$. Then  $\rk T = n = \dim \nolimits _{\Cee } E$, by Theorem \ref{n2n}
and Corollary \ref{coro}.


The $\Ree $-submodule $\Ree T$  of $V$ is $\Lin W$-invariant.\;Let
us consider the restriction of $\langle \,\cdot \mid \cdot\, \rangle $ to $\Ree T$. We claim that this restriction has only real values. Indeed, $\re \langle \,\cdot \mid \cdot\, \rangle $ defines a Euclidean structure on $\Ree T$ such that $\Lin W$ is orthogonal with respect to this structure. But $\Cee (\Ree T) = V$ because of irreducibility. Hence there exists a~canonical extension of $\re \langle \,\cdot \mid \cdot\, \rangle $
up to
a~Hermitian $\Lin W$-invariant scalar product, say $(\,\cdot \mid \cdot\, )$, on $V$.\;Thus, we have two $\Lin W$-invariant Hermitian structures on $V$, namely, $\langle \,\cdot \mid \cdot\, \rangle $ and $(\,\cdot \mid \cdot\, )$.
They are proportional because of irreducibility of $\Lin W$:
$$
\langle \,\cdot \mid \cdot\, \rangle =\lambda (\,\cdot \mid \cdot\, )\quad \text{for some}\quad \lambda \in \Cee .
$$
Considering their
restrictions to $\Ree T$, we obtain $\lambda = 1$. In other words, $\Ree T$ is a~real form of $V$ and the restriction of the action of $\Lin W$ to $\Ree T$ gives a finite real irreducible
$r$-group (and $\Lin W$ itself is the complexification of this group). It follows from the classification that this group is generated by $n$ reflections, see Section  \ref{1.5}. Therefore, $\Lin W$ is generated by $n$ reflections, too. From Theorem~\ref{sdsd}
it follows that $W$ is a~semidirect product. Let $a\in E$ be a~special point of $W$. Then  $a+\Ree T$ is a~real form of $E$. It is clear that $a+\Ree T$ is $W$-invariant. The restriction of $W$ to $a+\Ree T$ is a~real form of $W$. Therefore, this restriction is an affine Weyl group and $W$ is its complexification. This completes the proof of Theorem~\ref{cncg}.
\hfill \qed

\section{\bf Invariant lattices}\label{SinvLa}

As we have seen In Section  \ref{2.3}, one of the ingredients of the description of infinite complex crystallographic $r$-groups $W$ is the lattice $\Tran W$. In this section  we show how one can find all
full rank lattices in $V$ that are invariant with respect to
a given
finite $r$-group $K$ in $\GL(V)$.

We use the following \textit{notation}:

\vskip 1.5mm

$K \subseteq \GL(V)$ a~finite essential $r$-group, $n = \dim \nolimits _{\ka} V$,

$\Gamma \subseteq V$ a~$K$-invariant lattice,

$R_j:= R_{e_j,\mu _j}$ for $1 \leqslant j \leqslant s$, a
system of
reflections generating $K$,

$\mathcal L$ the set of root lines of all reflections in $K$ (see Section \ref{more}),

$\Gamma _\ell :=\ell \cap \Gamma $, where $\ell \in \mathcal L $,

$\Gamma _j:=\Gamma _{\ell _j}$ for $j=1,\ldots ,s$.

\subsection{Root lattices}\label{rl}

The set $$\Gamma ^0=\sum\limits _{\ell \in \mathcal L }\Gamma _\ell $$ \textit{is called the root lattice associated with} $\Gamma $.
\textit{If $\Gamma =\Gamma ^0$, then $\Gamma $ is called a~root lattice}.

\subsubsection{Theorem}\label{SSS}
{\it The set $\Gamma ^0$ is a~$K$-invariant lattice and $\rk \Gamma = \rk \Gamma ^0$.
}

\begin{proof}
It is clear that $\Gamma ^0$ is a $K$-invariant lattice, so let us prove the assertion about ranks.

As $K$ is essential, after a suitable renumbering, we can assume that
$$
V=\bigoplus\limits _{k=1}^n \ell _k.
$$
Put
$$
S=(\ii-R_{e_1,\mu _1})+\cdots +(\ii-R_{e_n,\mu _n}).
$$
If $v\in \Ker S$, then
$$
Sv=0=(1-\mu _1)\langle v\,|\,e_1 \rangle e_1+\cdots +(1-\mu _n)\langle v\,|\,e_n \rangle e_n,
$$
and therefore, $(1-\mu _j)\langle v\,|\,e_j \rangle e_j=0$, because $e_j\in \ell _j$.
Hence $\langle v\,|\,e_j \rangle =0$ for $1 \leqslant j \leqslant n$, i.e., $v=0$. So, $S$ is \textit{nonsingular}. But
$$
S \Gamma \subseteq \Gamma _{\ell _1}\oplus \ldots \oplus \Gamma_{\ell _n} \subseteq \Gamma ^0 \subseteq \Gamma
$$
 and $\rk S \Gamma = \rk \Gamma $. Therefore, $\rk \Gamma ^0 = \rk \Gamma $.
\end{proof}

\subsubsection{Corollary}\label{nnnn}

{\it
\begin{enumerate}[\hskip 7.2mm\rm 1)]\itemsep=-.1ex
\item If $\rk \Gamma = 2n$ {\rm(}hence $\ka = \Cee ${\rm)}, then $\rk \Gamma _\ell = 2$ for every $\ell \in \mathcal L $.
\item If $\rk \Gamma = n$, then $\rk \Gamma _\ell =1$ for every $\ell \in \mathcal L $.
\end{enumerate}
}

\begin{proof}
1) We may take $\ell _1=\ell $ in the previous proof. As this proof shows that $\rk(\Gamma _{\ell _1}\oplus \cdots \oplus \Gamma _{\ell _n} ) = \rk \Gamma _{\ell _1}+\cdots +\rk \Gamma _{\ell _n}=2n$, the first assertion follows from the evident inequality $\rk \Gamma _\ell \leqslant 2$.

2) The assertion can be proved similarly by means of reduction to the real form.
\end{proof}

It appears that one can reconstruct $\Gamma ^0$ from $\Gamma _1,\ldots, \Gamma _s$.

\subsubsection{Theorem}\label{T??}
 $$\Gamma ^0=\Gamma _1+\cdots +\Gamma _s.$$

\begin{proof}
Let $\ell \in \mathcal L $ and $u\in \Gamma _\ell $. Then  there exists $g\in K$ such that $gu\in \Gamma _j$ for a~certain $j$ (because every reflection in $K$ is conjugate to
a~power of some $R_j$, where $1 \leqslant j \leqslant s$, see Theorem \ref{ctpgr}).
Let
$$
\Gamma ':=\Gamma _1+\cdots +\Gamma _s.
$$
It is easy to see that $\Gamma '$ is $K$-invariant, hence $u\in \Gamma '$ and $\Gamma =\Gamma '$.
\end{proof}

The problem of finding all $K$-invariant lattices can be solved in two steps:
\begin{enumerate}[\hskip 7.2mm\rm 1)]\itemsep=-.1ex
\item
describing all $K$-invariant root lattices;
\item
describing all $K$-invariant lattices with a~fixed associated root lattice.
\end{enumerate}

First, we shall show how to solve problem 2).

\subsection{The lattices with a~fixed root lattice}

Set
$$
\Gamma ^*:=\{v\in V \mid (\ii -P)v\in \Gamma \;\; \mbox{for every $P\in K$}\}.
$$
Clearly, $\Gamma ^*$ is a~subgroup of $V$. It is more convenient to use another description of $\Gamma ^*$.

Let
$$
\pi :V \tto V/\Gamma
$$
be the~natural map and let $(V/\Gamma )^K$ be the set of fixed points of $K$ in $V/\Gamma $.
Then  $\Gamma ^*=\pi ^{-1}((V/\Gamma )^K)$. Therefore,
$$
\Gamma ^*=\{v\in V \mid (\ii-R_j)v\in \Gamma _j \;\;\mbox{for each $j=1,\ldots, s$}\}.
$$
It is clear that $\Gamma ^*$ is \textit{$K$-invariant} and $\Gamma \subseteq \Gamma ^*$.

\subsubsection{Theorem}\label{T??}
{\it
$\Gamma $ is a~lattice and $\rk \Gamma ^* = \rk \Gamma $.
}

\begin{proof}
First, let us show that $\Gamma ^*$ is a~lattice. If it is not a~lattice, then there exists a~vector $v\in \Gamma ^*$ such that
$$
\alpha v\in \Gamma ^*\;\; \mbox{for every $\alpha \in \Ree$}
$$
(because $\Gamma ^*$ is a~closed subgroup of $V$). Then
$$
(\ii-P)\alpha v\in \Gamma \;\; \mbox{for every $\alpha \in \Ree$ and $P\in K$.}
$$
Therefore, $(\ii-P) v = 0$, as $\Gamma $ is a~lattice, and hence $ v = 0$, because $K$ is an essential group. Therefore, $\Gamma ^*$ is a~lattice.

The space $V$ is Euclidean with respect to $\re \langle \,\cdot \mid \cdot\, \rangle $. Let $\Ree \Gamma ^\perp $ be the orthogonal complement of $\Ree \Gamma $ in this space. Let $v\in \Gamma ^*$ and $v = u + w$ for $u\in \Ree \Gamma $ and $w\in \Ree \Gamma ^\perp $. Then
$$
(\ii-P)v = (\ii-P)u + (\ii-P)w \subseteq \Gamma \subseteq \Ree \Gamma \;\; \mbox{for every
 $P\in K$}.
$$
Therefore, $(\ii-P )w = 0$, hence $w = 0$ (again because $K$ is essential). Thus, $\Gamma ^*\subseteq \Ree \Gamma $. This com\-p\-letes the proof.
\end{proof}

\subsubsection{Cohomological meaning of $\Gamma ^*$}

We have the exact sequence of groups
$$
0 \tto \Gamma \tto V \tto V/\Gamma \tto 0.
$$
It yields the exact cohomological sequence
$$
H^0(K,V) \tto H^0(K,V/\Gamma )\tto H^1(K,\Gamma ) \tto H^1(K,V)\tto \ldots .
$$
But $H^0(K,V)=0$ because $K$ is essential, $H^0(K,V/\Gamma )=(V/\Gamma )^K=\Gamma ^*/\Gamma $,
and 
 $H ^1(K,V) = 0$ because $V$ is divisible. Therefore,
$$
\Gamma ^*/\Gamma \simeq H^1(K,\Gamma ).
$$

Now we can explain how to find all $K$-invariant lattices with a~fixed $K$-invariant root lattice.

\subsubsection{Theorem}\label{latti}
{\it
Let $\Lambda $ be a~fixed $K$-invariant root lattice in $V$. Then  for every lattice $\Gamma $ in $V$, the following properties are equivalent:
\begin{enumerate}[\hskip 7.2mm\rm a)]\itemsep=-.1ex
\item
$\Gamma$ is a~$K$-invariant lattice and $\Gamma ^0=\Lambda $.
\item
$\Lambda \subseteq \Gamma \subseteq \Lambda ^*$ and $\Gamma _j=\Lambda _j$ for each $j=1,\ldots,  s$.
\end{enumerate}

\noindent
There are only finitely many
lattices $\Gamma $
with
properties {\rm a)} and {\rm b)}.
}

\begin{proof}
a) $\Rightarrow $ b). Let $v\in \Gamma $. Then
$$
(\ii-R_j)v=(1-\mu _j)\langle v\,|\,e_j \rangle e_j\in \Gamma \cap \ell _j=\Gamma _j \subseteq \Gamma ^0=\Lambda .
$$
Therefore, $v\in \Lambda ^*$ and $\Gamma \subseteq \Lambda ^*$.

b) $\Rightarrow $ a). Let $\Lambda \subseteq \Gamma \subseteq \Lambda ^*$. Then  $\Gamma =\pi ^{-1}\pi (\Gamma )$, where $\pi :V \tto V/\Lambda $ is the natural map, and $\pi (\Gamma )\subseteq \pi (\Lambda ^*)=(V/\Lambda )^K$. Hence $\pi (\Gamma )$ is $K$-invariant, and therefore, $\Gamma $ is also $K$-invariant. If $\Gamma _j=\Lambda _j$, then $\Gamma ^0=\Lambda $ (because $\Lambda =\Lambda ^0=\Lambda _1+\cdots +\Lambda _s$).
\end{proof}


\subsubsection{}
\hskip -2mm
We have already seen that for $K$ irreducible, one can take $s = n$ if $\ka = \Ree $ (see Section \ref{1.5}), and $s = n$ or $n + 1$ if $\ka = \Cee $ (see Theorem \ref{cfi}).
It appears that {\it if $s = n$, then there exists a~good
constructive
way to
find $\Lambda ^*$ by means of $\Lambda $.}

\subsubsection{Theorem}\label{opS}
{\it
Let $s = n$.
Consider the nonsingular \textup{(}see the proof of Theorem {\rm\ref{SSS}}\textup{)} linear operator
$$
S=(\ii-R_1)+\cdots +(\ii-R_n)
$$
Let $\Lambda $ be a~$K$-invariant lattice in $V$.  Then
\begin{enumerate}[\hskip 7.2mm\rm a)]\itemsep=-.1ex
\item
$\Lambda ^* \subseteq S^{-1}\Lambda $.
\item
If $\Lambda ^0=\Lambda $, then $\Lambda ^*=S^{-1}\Lambda $.
\end{enumerate}
}

\begin{proof}
a) Let $a\in \Lambda ^*$, then, by definition, $(\ii-R_j) a\in \Lambda $ for all $1 \leqslant j \leqslant n$, hence $Sa\in \Lambda $ and $a\in S^{-1}\Lambda $.

b) Let $\Lambda =\Lambda ^0$ and $u\in S^{-1}\Lambda $. Then
$$
Su=(\ii-R_1)u+\cdots +(\ii-R_n)u\in \Lambda =\Lambda _1 \oplus \cdots \oplus \Lambda _n.
$$
As $(\ii-R_j)u\in 
\ka e_j$ and $e_1,\ldots ,e_n$ are linearly independent, this yields
$(\ii-R_j)u\in \Lambda _j \subseteq \Lambda $.
Therefore, by the definition of $\Lambda ^*$, we have $u\in \Lambda ^*$.
\end{proof}

This theorem is useful in practice because one can explicitly describe the operator $S$: its matrix with respect to the basis $e_1,\ldots ,e_n$ is
$$
\left (\begin{array}{ccc}
(1-\mu _1)\langle e_1\,|\,e_1 \rangle &\ldots &(1-\mu _1)\langle e_n\,|\,e_1 \rangle \\
&\ldots \\
(1-\mu _n)\langle e_1\,|\,e_n \rangle &\ldots &(1-\mu _n)\langle e_n\,|\,e_n \rangle
\end{array}
\right ).
$$

Therefore, if $s = n$, then
the problem of
classifying
 $K$-invariant lattices with a fixed root lattice $\Lambda$ can be practically solved
as follows:

\vskip 1mm

{\it First, find $\Lambda^*$ using Theorem {\rm{\ref{opS}}b)}, and then look over all lattices intermediate between $\Lambda$ and $\Lambda^*$, applying Theorem {\rm\ref{latti}b)}.}

\vskip 1mm

\subsubsection{}\label{sn+1}

How can this problem
be solved 
if
$s=n + 1$?

\vskip 1.3mm

It is not difficult to see that one can take $R_1,\ldots ,R_{n+1}$ in such a~way that $R_1,\ldots ,R_n$ generate a~subgroup $K'$ of $K$, which is also \textit{irreducible} (the numbering in Table~1
has this property).

\subsubsection{Example}
$$
K=K_{31}
\quad \begin{array}{c}
\includegraphics[width=50mm]{t35new-eps-converted-to}
\end{array}
\qquad K'=K_{29} \begin{array}{c}
\includegraphics{t33-eps-converted-to}
\end{array}
$$
Given this, the solution to the specified problem when $s=n+1$ is the following:

\vskip 1mm
\textit{Following the above strategy for the case of $s=n$, find all $K'$-invariant lattices and then
select
among
them the ones invariant under} $R_{n+1}$.

\vskip 1mm

In Section \ref{4.8}, we will give rather detailed examples of how this can be practically done, after we have explained (in Section \ref{4.7}) how to find invariant root lattices.

By now we want to discuss several general properties of invariant lattices and explain the significance of root lattices in the theory of infinite $r$-groups.

\subsection{Further remarks on lattices}\label{frl}

\subsubsection{Theorem}\label{genK}
{\it
Let $F \subseteq \GL(V)$ be a~finite
irreducible group and $n = \dim\nolimits_{\ka} V$.\;If $\;\Gamma $ is a~nonzero $F$-invariant lattice in $V$, then  $\rk \Gamma=n$
for $\ka = \Ree $, and $\rk \Gamma = n$ or $2n$ for $\ka = \Cee $. Moreover, if $F$ is an $r$-group, $\ka = \Cee $, and $\rk \Gamma = n$, then $F$ is the complexification of the Weyl group of a~certain irreducible root system in a real form of $V$.
}

\vskip 1mm
\noindent{\it Proof}
as in Sections \ref{3.1} and \ref{3.4}.\hfill \qed

\vskip 2mm
\subsubsection{Corollary}\label{??}
{\it  Keeping the notation of Theorem {\rm\ref{genK}}, let $\ka = \Cee $ and let $F$ be an $r$-group. If $F$ has an invariant lattice of rank $n$ in $V$, then $F$ has an invariant lattice of rank $2n$ in $V$.
}

\begin{proof}
By Theorem \ref{genK},
the group $F$ is the complexification of the~Weyl group $W$ of a certain irreducible root system in  a real form of $V$.\;Let
$\Gamma$ be a $W$-invariant lattice of rank $n$ in this real form of $V$ 
(for example, such is the lattice generated by the specified root system).
Then  for every $z\in \Cee \setminus \Ree $, the~lattice $\Gamma + z \Gamma $ is $F$-invariant and has rank $2n$.
\end{proof}

It should be noted that if $\ka = \Ree $ and $K$ is the~Weyl group of a root system in $V$, then the~root lattices of rank $n$ in $V$ (as definded in Section \ref{rl}) are
the~lattices $\Lambda $ of \textit{radical weights}\index{radical weights} of the root systems in $V$ whose Weyl group is $K$, and $\Lambda ^*$ is the~lattice of weights of such a root system; see
\cite{1}.
For such $\Lambda$, the matrix of $S$ with respect to a~basis of $\Lambda$ formed by simple roots, is
the \textit{Cartan matrix}\index{Cartan matrix} of the corresponding root system.

\subsection{Properties of the operator \boldmath $S$}\label{4.4}

\subsubsection{Theorem}\label{proS}
{\it
Let $K$ be irreducible, let $s = n$, let $\Lambda $ be a~nonzero $K$-invariant lattice in $V$, and let $S = (\ii-R_1) + \cdots + (\ii-R _n)$. Then
\begin{enumerate}[\hskip 6.2mm\rm a)]\itemsep=-.1ex
\item
$\begin{array}[t]{@{}l@{}}
\text{$|{\mbox{\rm det}}\, S|$ and $\Tr S\in \Zee$ if
$\ka=\mathbb R$},\\[-1mm]
\text{$|{\mbox{\rm det}}\, S|^2$ and $2 \re \Tr S \in \Zee  $ if
$\ka=\mathbb C$.}
\end{array}
$
\item
$|{\mbox{\rm det}}\, S|$ depends only on $K$ but not on the
generating
system of reflections in\;$K$.
\item
$\begin{array}[t]{@{}l@{}}
\text{$|{\mbox{\rm det}}\, S|$ is divisible by $[\Lambda :\Lambda ^0]$ if $\rk \Lambda = n$,}\\[-1mm]
\text{$|{\mbox{\rm det}}\, S|^2$ is divisible by $[\Lambda :\Lambda ^0]$ if $\rk \Lambda = 2n$.}
\end{array}
$
\item
$\Lambda =\Lambda ^0 \Rightarrow [\Lambda ^*:\Lambda ]=\begin{cases}
|{\mbox{\rm det}}\, S| &\text{if }\rk \Lambda = n,\\[-1mm]
|{\mbox{\rm det}}\, S|^2 &\text{if }\rk \Lambda = 2n.
\end{cases}$
\item
If $\Lambda =\Lambda ^0$ and $d_1,\ldots ,d_r$ are the invariant factors of a
matrix of the endomorphism of $\Lambda $ defined by $S$ arranged so that $d_1|\ldots |d_r$ and $d_j>0$, then
$$
H^1(K,\Lambda )\simeq \Lambda ^*/\Lambda \simeq \Zee  /d_1 \Zee  \oplus \ldots \oplus \Zee  /d_r \Zee  .
$$
Specifically, $d_1,\ldots ,d_r$ do not depend on the generating system of reflections in $K$.
\item
$\Lambda =\Lambda ^0 \Rightarrow |H^1(K,\Lambda )|=\begin{cases}
|{\mbox{\rm det}}\, S|&\text{if }\rk \Lambda = n,\\[-1mm]
|{\mbox{\rm det}}\, S|^2 &\text{if }\rk \Lambda = 2n.
\end{cases}$
\end{enumerate}
}

\begin{proof}
If $\rk \Lambda = n$, then, by considering the corresponding real form, we reduce the problem to the case of $\ka = \Ree $.

If $\rk \Lambda = 2n$, then $\ka = \Cee $. It is known that, in this case, for every $P \in \GL(V)$ one has
$$
{\mbox{\rm det}}\, P_{\Cee \mid \Ree }=|{\mbox{\rm det}}\, P|^2 \quad \text{and}\quad \Tr P_{\Cee \mid \Ree }= 2 \re \Tr P
$$
(here, $P_{\Cee \mid \Ree }$ is $P$, considered as the~linear operator of the~$2n$-dimensional real vector space $V$).

Now, a) follows from the fact that a~basis of $\Lambda $ is an $\Ree $-basis of $V$. We have also:
$$
\Lambda ^0 \subseteq \Lambda \subseteq (\Lambda ^0)^*=S^{-1}\Lambda ^0,\quad \text{so}\quad [\Lambda :\Lambda ^0]\mid [S^{-1}\Lambda ^0:\Lambda ^0].
$$
But $S^{-1}\Lambda ^0/\Lambda ^0 \simeq \Lambda ^0/S \Lambda ^0$, whence
$$
[S^{-1}\Lambda ^0:\Lambda ^0]=\begin{cases}
[(\Lambda ^0)^*]:\Lambda ^0]=|{\mbox{\rm det}}\, S| &\text{if}\quad \ka=\Ree, \\[-1mm]
|{\mbox{\rm det}}\, S_{\Cee \mid \Ree }|&\text{if}\quad \ka=\Cee .
\end{cases}
$$
The assertions b), c), d), e), f) follow from these equations.
\end{proof}

\subsubsection{Corollary}\label{Cooo}
{\it \begin{enumerate}[\hskip 5mm\rm a)]\itemsep -.1ex

  \item If $|{\mbox{\rm det}}\, S|=1$, then $\Lambda =\Lambda ^0$, i.e., every $K$-invariant lattice in $V$ is a~root lattice.

\item  Let $\ka = \Cee $ and suppose $|{\mbox{\rm det}}\, S|^2$ is a~prime number. Then  $\Lambda =\Lambda ^0$ or $\Lambda =(\Lambda ^0)^*$.
    \end{enumerate}
}

\subsubsection{Remark}
One can calculate ${\mbox{\rm det}}\, S$ \textit{directly from the graph of} $K$. Namely,
$$
{\mbox{\rm det}}\, S=\sum\limits _\sigma \sgn \sigma \cdot a_\sigma,
$$
where $\sigma $ runs through the permutations of degree $n$, and $a_\sigma :=c_\alpha c_\beta \cdots c_\gamma $ for the cycle decomposition  $\sigma =\alpha \beta \cdots \gamma $ (i.e., expressing of $\sigma$
as a product of disjoint cycles).

\subsubsection{Examples}\label{exa}

a) $K=K_1$, type ${\sf A}_n$. 
We denote by $S({\sf A}_n)$ the operator $S$ for the graph
$$
\includegraphics{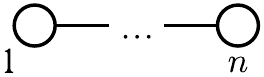}
$$
The nonzero simple cyclic products
are $c_{j,j+1}=1$ for $1 \leqslant j \leqslant n-1$, and $c_j=2$ for $1 \leqslant j \leqslant n$.

We have ${\mbox{\rm det}}\, S({\sf A}_1) = 2 $. Assume, by induction, that ${\mbox{\rm det}}\, S({\sf A}_k) = k + 1$ if $k < n$. Then
\[
{\mbox{\rm det}}\, S({\sf A}_n) = c_1{\mbox{\rm det}}\, S({\sf A}_{n-1}) - c_{12} \cdot {\mbox{\rm det}}\, S({\sf A}_{n-2}) = 2n - (n-1) = n + 1.
\]

b) $K = K_2$, type $G(m,m,n)$. Let us consider the generating system of reflections given by the graph
$$
\includegraphics{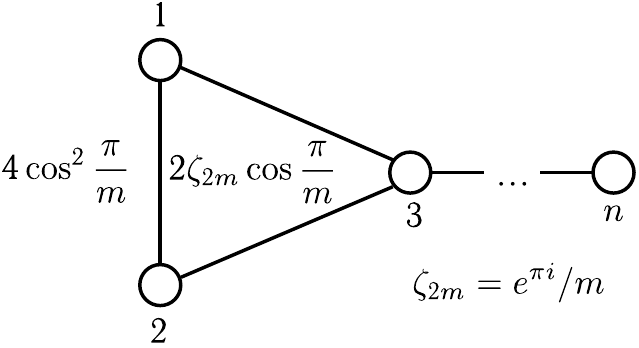}
$$

The list of all nonzero simple cyclic products is
\begin{align*}
c_j&=2\;\; \mbox{for $1 \leqslant j \leqslant n$;}\\
c_{13}&=c_{23}=c_{34}=c_{45}=\ldots =c_{n-1,n}=1;\\
c_{12}&=4 {\mbox{\rm cos}}^2 \frac{\pi}{m};\\
c_{123}&=2 {\mbox{\rm cos}}\,\frac{\pi}{m}\,e^{i/m};\\
c_{132}&=2 {\mbox{\rm cos}}\,\frac{\pi}{m}\,c^{-\pi i/m}.
\end{align*}

It is easy to see that
\begin{multline*}
{\mbox{\rm det}}\, S=c_1 {\mbox{\rm det}}\, S({\sf A}_{n-1})-c_{12}{\mbox{\rm det}}\, S({\sf A}_{n-2})-c_{13}\cdot c_2 {\mbox{\rm det}}\, S({\sf A}_{n-3})\\
+c_{123}{\mbox{\rm det}}\, S({\sf A}_{n-3})+c_{132}{\mbox{\rm det}}\, S({\sf A}_{n-3})
=2n-4 {\mbox{\rm cos}}^2 \frac{\pi}{m}\,(n-1)-1 \cdot 2 \cdot (n-2)\\
+2 {\mbox{\rm cos}}\,\frac{\pi}{m}\,e^{\pi i/m}(n-2)+2 {\mbox{\rm cos}}\,
\frac{\pi}{m}\,e^{-\pi i/m}(n-2)\\
=4-4 {\mbox{\rm cos}}^2 \frac{\pi}{m}\,(n-1)+4 {\mbox{\rm cos}}^2 \frac{\pi}{m}\,(n-2)=4 {\mbox{\rm sin}}^2 \frac{\pi}{m}.
\end{multline*}

\subsubsection{Remark}\label{corgg}
If $\ka = \Ree$ and $\Lambda $ is the lattice of radical weights of a root system ${\sf R}$ in $V$ with the Weyl group $K$, then $d_1,\ldots ,d_r$ are invariant factors of the~Cartan matrix of  ${\sf R}$
and $H^1 (K,\Lambda )$ is isomorphic to the center of a~complex simply connected semisimple algebraic group $G$ whose affine Weyl group is
$K\ltimes \Lambda$ (thus the root system of $G$ is ${\sf R}^{\vee}$, the dual of ${\sf R}$). 

\subsection{Root lattices and infinite \boldmath $r$-groups}\label{4.5}

We shall now show that if there exists a~nonzero $K$-invariant lattice $\Lambda $, then there also exists an infinite $r$-group $W$ with $\Lin W = K$ and $\Tran W =\Lambda $.

\subsubsection{Theorem}\label{ssp}
{\it
Let $K \subseteq \GL(V)$ be a~finite
$r$-group and let $\Lambda $ be a~nonzero $K$-invariant lattice in $V$. Then  the semidirect product of $K$ and $\Lambda $ is an $r$-group if and only if $\Lambda =\Lambda ^0$.
}

\begin{proof}
Let $W$ be the semidirect product of $K$ and $\Lambda $. Let $a$ be a~special point. We identify $W$ and $\kappa _a(W) = K \ltimes \Lambda \subseteq \GL(V)\ltimes V$. Then  the set of all reflections in $W$ is
$$
\{(R,v)\mid R \text{ is a~reflection of $K$ and $v\in \ell _R \cap \Lambda $}\}.
$$
Let $W^0$ be a~subgroup of $W$ generated by this set. Then  $K \subseteq W^0$. Also $\Lambda ^0 \subseteq W^0$, because $\Lambda ^0=\sum
 _{\ell \in \mathcal L }\Lambda _\ell $.
Moreover, if $v\in \Lambda _\ell $ and $\ell $ is the root line of a~reflection $R$, then $(R,v)(R^{-1},0)=(1,v)$. Therefore, we have $K \ltimes \Lambda ^0 \subseteq W^0$. But every reflection in $W$ lies in $K \ltimes \Lambda ^0$ (because it has the form $(R,v)$, where $R$ is a~reflection and $v\in \ell _R \cap \Lambda =\Lambda ^0$, see Proposition \ref{prope}).
Therefore, $W^0=K \ltimes \Lambda ^0$ is generated by reflections.
Hence $K \ltimes \Lambda = W = W^0 = K \ltimes \Lambda ^0$ if and only if $\Lambda =\Lambda ^0$.
\end{proof}

\subsubsection{Corollary}\label{sdprp}
{\it
Let $K$ be a finite $r$-group in $\GL(V)$ generated by $n=\dim_{\ka}V$ reflections. Then all infinite $r$-groups $W$ in $A(E)$ with $\Lin W=K$
are exactly all groups $ 
K \ltimes \Lambda $, where $\Lambda $ is a~nonzero
root lattice for $K$.
}

\subsection{Description of the group of linear parts; proof}\label{4.6}

Now we shall prove
Theorem \ref{cirg2}
and a part of Theorem \ref{cirg3}.


\begin{proof}
a) $\Rightarrow $ b) is already proved in Section \ref{frl}.

b) $\Rightarrow $ a) is trivial.

c) $\Rightarrow $ b) is also trivial: such a~lattice is $\Tran W$.

b) $\Rightarrow $ c) Let $\Gamma $ be an invariant lattice of rank $2n$. Then  the semidirect product of $K$ and $\Gamma ^0$ is a~crystallographic $r$-group because of
Theorem \ref{ssp}.

b) $\Rightarrow $ d) Let us first prove that $\Zee  [\Tr K]$ coincides with the ring with unity generated over $\Zee  $ by all cyclic products.

We have $\Zee  [\Tr K] = \Zee  [ \Tr \Zee  K ] $. Indeed, clearly $\Tr K \subseteq \Tr \Zee  K $, hence $\Zee  [\Tr K] \subseteq \Zee  [ \Tr \Zee  K ]$. The reverse inclusion follows from the fact that $\Tr$ is an additive function. But
$\ii -R_1,\ldots \ii-R_s$
generate the ring $\Zee  K$ (here $R_1,\ldots R_s$
is a~generating system of reflections of $K$). Therefore, the monomials $(\ii-R_{j_1})\cdots (\ii-R_{j_r})$ generate $\Zee K $ as a~$\Zee  $-module. We have
$$
(\ii-R_{j_1})(\ii-R_{j_r})(\ii-R_{j_{r-1}})\cdots (\ii-R_{j_2})e_{j_1}=c_{j_1 \ldots j_r}e_{j_1},
\eqno\boxed{3}
$$
and it easily follows from this equality that
$$
\Tr (\ii-R_{j_1})(\ii-R_{j_r})(\ii-R_{j_{r-1}})\cdots (\ii-R_{j_2})=c_{j_1 \ldots j_r}.
$$
Therefore, $\Zee  [\Tr \Zee  K]=\Zee  [\ldots ,c_{j_1 \ldots j_r},\ldots ]$ and we are done.

By now, let $\Gamma $ be a~$K$-invariant lattice of rank $2n$. It follows from equality $\,\boxed{3}\,$ that
$$
c_{j_1 \ldots j_r}\Gamma _{j_1}\subseteq \Gamma _{j_1} \quad \text{for every}\quad j_1,\ldots ,j_r.
$$
But $\rk \Gamma _j=2$, see Section \ref{nnnn}.
Therefore, $c_{j_1 \ldots j_r}$ is an integral element of a~certain purely imaginary quadratic extension of $\Qee $; denote this extension by $L_{j_1}$.

Let $L = L_1$; let us show that $c_{j_1 \ldots j_r} \in L$ for every $j_1,\ldots ,j_r$. The group $K$ being irreducible, there exists
$$
\alpha =c_{1 l _1 l _2 \ldots l _tj_1 l _t l _{t-1}\ldots l _1}\ne 0
$$
$$
\hskip 8mm \includegraphics[width=75mm]{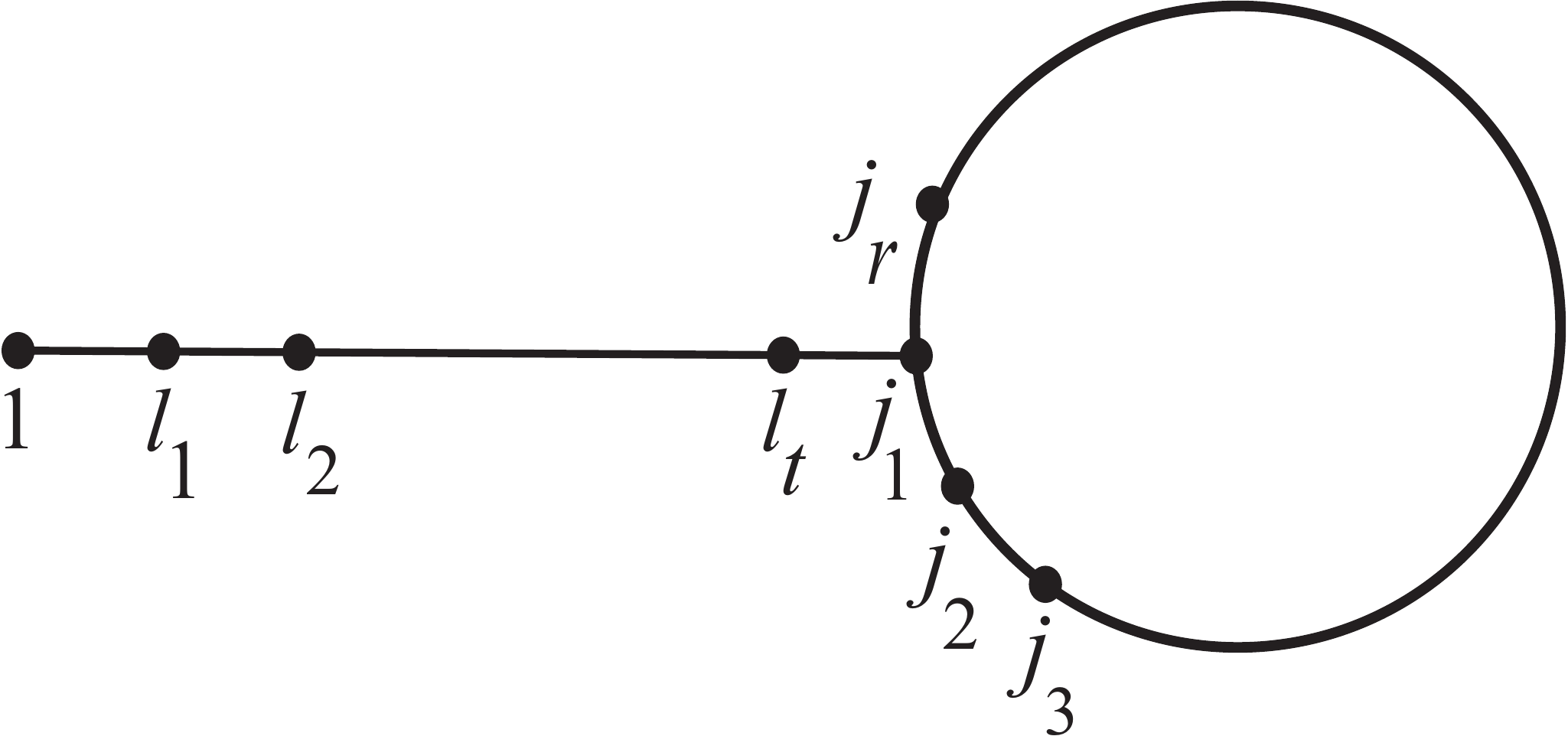}
$$
\noindent
for certain $l _1,\ldots ,l _t$. Let
\begin{align*}
\beta &=c_{1 l _1 l _2 \ldots l _tj_1j_2 \ldots j_rj_1 l _t l _{t-1}\ldots l _1},\\
\gamma &=c_{j_1j_2 \ldots j_r}.
\end{align*}
Then  $\alpha ,\beta \in L$ and $\beta =\alpha \gamma $. But $\alpha \ne 0$, hence $\gamma =\beta /\alpha \in L$. Therefore, the ring $\Zee  [\Tr K]$ lies in the maximal order of $L$.

d) $\Rightarrow $ e) This is proved in \cite[Lemma 1.2]{7}.

d) $\Rightarrow $ a) Let $\Zee  [ \Tr K] \subseteq D$, where $D$ is the maximal order of a~certain purely imaginary quadratic extension $L$ of $\Qee $. The group $K$ being irreducible and $D$ being integrally closed, one can again apply \cite[Lemma 1.2]{7} and see that $K$ is defined over $D$. Therefore, there exists a~$K$-invariant $D$-submodule $\Gamma $ of $V$ such that the natural map $\Gamma \bigotimes\nolimits _D \Cee \tto V$ is an isomorphism. But $D$ is a~Dedekind ring and $\Gamma $ is a~torsion-free $D$-module of rank $n$. Therefore, $\Gamma $ is isomorphic to a~direct sum of $n$ fractional ideals of the field $L$, see, e.g., \cite[Theorem 22.5]{8}. Let $J_1,\ldots ,J_n$ be these ideals. Then  there exists a~$\Cee $-basis $v_1,\ldots ,v_n$ of $V$ such that
$$
\Gamma =J_1v_1+\cdots +J_nv_n.
$$
But $D$ is a~lattice of rank 2 in $\Cee $ and, for every fractional ideal $J$ of $L$, there exists a~nonzero $d\in D$ such that $d \cdot J \subseteq D$. Hence, $J$ is
also
a~lattice of rank 2 in $\Cee $.

e) $\Rightarrow $ d) Let $K$ be defined over a~purely imaginary quadratic extension $L$ of $\Qee $. Then  $\Tr P\in L$ for every $P \in K$, hence $\Zee  [\Tr K] \subseteq L$. But $\Tr P$ is an integral algebraic number. Therefore, $\Zee [ \Tr K]$ lies in a~maximal order of $L$.

d) $\Rightarrow $ f) Using Table~1,
we can easily find those $K$ for which $\Zee  [ \Tr K]$ lies in the maximal order of a~certain purely imaginary quadratic extension of $\Qee $. When finding, it is convenient to use that each generator of such a  ring $\Zee  [ \Tr K]$
should necessarily be an integral element of a~certain purely imaginary quadratic extension of $\Qee $. This necessary condition is verified by means of the following criterion:



\vskip 2mm

\textit{A number $z\in \Cee $ is an integral element of a~certain purely imaginary quadratic extension of $\Qee $ if and only if $|z|^2 \in \Zee  $ and $2 \re z \in \Zee  $}.

\vskip 2mm
 \textit{A posteriori} it appears that this necessary condition is also sufficient.

The output of this finding are precisely the groups listed in f).

\subsubsection{Example}
Let $K = K_2$, type $G(m,m,s)$, $s \geqslant 3$. The graph of $K$ is
$$
\includegraphics{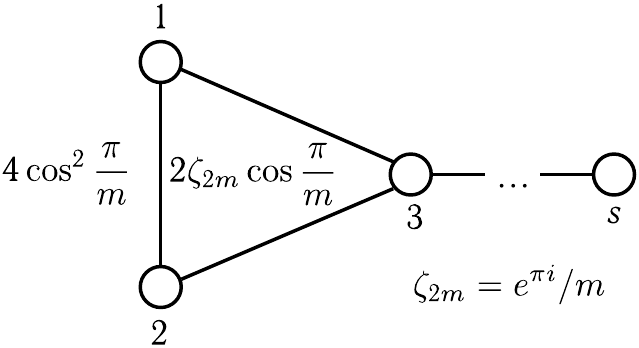}
$$
We have $\Zee  [ \Tr K] = \Zee  [e^{2\pi i/m}]$. The specified necessary condition for the generator $e^{2\pi i/m}$ gives: $2 {\mbox{\rm cos}}\,\frac{2 \pi }{m}\in \Zee  $.
This yields  $m = 2,3,4,6$. In these cases, $\Zee  [ \Tr K]$ lies in a~maximal order of a~certain purely imaginary extension of $\Qee $.
\end{proof}

\subsection{Description of root lattices}\label{4.7}
We assume now that $\ka = \Cee $. We first describe (up to similarity) all $K$-invariant root lattices of rank $2$n when $s = n$ (Section \ref{4.5} implies that only these lattices are of interest to us).

\subsubsection{Theorem}\label{2lat}
{\it
Let $K$ be a finite irreducible $r$-group in $\GL(V)$ generated by a system of reflections $R_1,\ldots, R_s$. For every $j=1, \ldots, s$, let $\Zee [\Tr K]_j$ be the subring of $\Zee [\Tr K]$ generated over $\Zee$ by
all cyclic products of the form $c_{j\ldots}$.
Let $\Lambda _j \subseteq \Cee e_j$ for $j=1,\ldots, s$ be a~set of lattices of rank $2$ and let $\Gamma =\Lambda _1+\cdots +\Lambda _s$. In order for $\Gamma $ to be a~$K$-invariant root lattice with $\Gamma _j=\Lambda _j$ for each $j=1,\ldots, s$, it is necessary, and if $s = n$, also sufficient that the following conditions hold:
\begin{enumerate}[\hskip 7.2mm\rm a)]\itemsep=-.1ex
\item
$\Zee  [\Tr K]_j\Lambda _j \subseteq \Lambda _j$ for each $j=1,\ldots,  s$.
\item
$(\ii -R_k)\Lambda _j \subseteq \Lambda _k$ and $(\ii-R_j)\Lambda _k \subseteq \Lambda _j$ for every $j$ and $k$.

Moreover, {\rm b)} is equivalent to the condition
\item
for every $j\neq k$ such that
$c_{kj}\ne 0$, one has
$$
(\ii-R_k)\Lambda _j \subseteq \Lambda _k \subseteq c_{kj}^{-1}(\ii-R_k)\Lambda _j.
$$
\end{enumerate}
}

\begin{proof} Let $\Gamma$ be $K$-invariant.

Property a) follows from the inclusion $c_{j_1 \ldots j_r}\Gamma _{j_1}\subseteq \Gamma _{j_1}$ 
because, by definition, $\Zee  [\Tr K]_{j_1}$ is generated
over $\Zee $
by all cyclic products of the form $c_{j_1\ldots}$.

The \lq \lq necessary part\rq \rq {} of b) follows from the $K$-invariance of the lattice.

Let us prove the \lq \lq sufficient part\rq \rq . If $s = n$, then $\Gamma $ is in fact a~direct sum of $\Lambda _j$ for $1 \leqslant j \leqslant n$; hence, $\Gamma $ is a~lattice. This lattice is invariant under $\ii-R_k$ for every $k$: if $k \ne j$, then the invariance follows from b); if $k= j$, then it follows from a). Hence, $\Gamma $ is $K$-invariant.

Now let us prove that b)\,$\Leftrightarrow$\,c).

b) $\Rightarrow $ c). One obtains the proof by applying the operator $\ii-R_k$ to both sides of the inclusion $ (\ii-R_j)\Lambda _k \subseteq \Lambda _j$.

c) $\Rightarrow $ b). 
Apply $\ii-R_j$ to both sides of the inclusion ${\Lambda _k \subseteq c_{kj}^{-1}(\ii-R_k)\Lambda _j}$.
\end{proof}

\subsubsection{Corollary}\label{Uniq}
{\it
Let $\Gamma $ be a~nonzero $K$-invariant lattice. For every $j<k$ such that  $|c_{kj}|=1$, the lattices $\Gamma _j$ and $\Gamma _k$ uniquely determine each other by the formulas
$$
\Gamma _k=(\ii-R_k)\Gamma _j \quad \text{and}\quad \Gamma _j=(\ii-R_j)\Gamma_k.
$$
}

\begin{proof}
By c) of Theorem \ref{2lat}, we have
$$
(\ii-R_k)\Gamma_j \subseteq \Gamma_k \subseteq c_{kj}^{-1}(\ii-R_k)\Gamma_k.
$$
As the index of the left lattice in the right
lattice
is $|c_{kj}|^2=1$,
the inclusions are in fact equalities. Applying $\ii-R_j$, we get
$$
c_{jk}\Gamma_j \subseteq (\ii-R_j)\Gamma_k \subseteq \Gamma_j.
$$
Again, because of the above reason, the inclusions are in fact equalities.
\end{proof}

If $s = n+1$, we also need to know $\Gamma^*$ for a~$K' $-invariant lattice $\Gamma $, and to select those $\Gamma \subseteq \Lambda \subseteq \Gamma ^*$ for which
\begin{enumerate}[\hskip 7.2mm\rm a)]\itemsep=-.1ex
\item
$\Lambda $ is $R_{n+1}$-invariant,
\item
$\Lambda _j=\Gamma _j$ for each $j=1,\ldots, s$.
\end{enumerate}
Therefore, we also need to know $(\Gamma ^*)^0$. We have the following description of this lattice:

\subsubsection{Theorem}\label{***}
{\it
\[
\Gamma _j^*=\bigcap\limits _{\substack{k \text{\normalfont{} such that }\\
c_{jk}\ne 0}}c_{jk}^{-1}(\ii-R_j)\Gamma _k\;\;\text{for each $j=1,\ldots,  s$}.
\eqno\boxed{4}
\]
In particular, $ \Gamma _j^*=\Gamma _j$ if there is a~number $k$ such that $|c_{jk}|=1$.
}

\vskip 2mm
\noindent{\it Proof.} For $\lambda\in\mathbb C$, we have
\[
\lambda e_j\in \Gamma _j^* \;\; \text{if and only if $(\ii-R_k)\lambda e_j\in \Gamma _k $ for each $k$
such that $c_{jk}\neq 0$},
\]
because $c_{jk}=0$ means that $(\ii-R_k)\lambda e_j=0$.

Assume that $c_{jk}\ne 0$, i.e., $\langle e_j \mid e_k \rangle \ne 0$. If $\lambda (\ii-R_k)e_j\in \Gamma _k$, then, applying $\ii-R_j$ to the both sides of this inclusion,
we obtain $\lambda c_{jk}e_j\in (\ii-R_j)\Gamma _k$, i.e., $\lambda e_j\in c_{jk}^{-1}(\ii-R_j)\Gamma _k$.

Vice versa, if for some $\lambda \in \mathbb C$ and every $k$ such that $c_{jk}\neq 0$, we have $\lambda e_j\in c_{jk}^{-1}(\ii-R_j)\Gamma _k$, then, applying $\ii-R_k$ to both sides of this inclusion, we obtain $(\ii-R_k)\lambda e_j\in c_{jk}^{-1}c_{jk}\Gamma _k=\Gamma _k$, i.e., $\lambda e_j\in \Gamma ^*$.

In order to prove the second assertion, let us apply $\ii-R_j$ to $\Gamma _k \supseteq (\ii-R_k)\Gamma _j$. We obtain
$$
\Gamma _j \supseteq (\ii-R_j)\Gamma _k \supseteq c_{jk}\Gamma _j;
$$
whence $c_{jk}\Gamma _j=\Gamma _j$ if $|c_{jk}|=1$.
Thus, the latter two inclusions are in fact the equalities; hence
\[
c_{jk}^{-1}(\ii-R_j)\Gamma _k=\Gamma _j, 
\]
which implies the claim in view of $\,\boxed{4}$.\qed


\vskip 2mm
\textit{
We assume now that $s = n$ unless otherwise stated}.

\vskip 1.5mm

\subsubsection{Algorithm for constructing \boldmath$K$-invariant root full rank lattices: Case 1}\label{Case1}

In this case, we will only consider the
groups $K$ from Theorem~\ref{1.6} with the following property:

\vskip 1mm

\textit{Every two nodes of the graph of $K$ {\rm(}see Table~{\rm 1)}
can be connected by a~path of edges such that the absolute value of the weight of each edge is equal to $1$}.

\vskip 1mm

These are the following groups:

$K_1$;

$K_2$, type $G(6,1,s)$ for $s \geqslant 2$, type $G(m,m,s)$ for $m = 2, 3, 4, 6$ and  $s \geqslant 3$;

$K_3$ for $m =3,4,6 $;

$K_4$; $K_8$; $K_{24}$; $K_{25}$; $K_{29}$; $K_{32}$; $K_{33}$; $K_{34}$; $K_{35}$; $K_{36}$; $K_{37}$.

\vskip 2mm

\noindent\textbf{Algorithm}

Take an arbitrary
lattice $\Delta$ of rank 2 in $\Cee $ such that
$\Zee  [\Tr K]\Delta\subseteq \Delta$.

Put $\Lambda _1=\Delta e_1$.

Given a vertex of the graph of $K$ with number $l \geqslant 2$,
consider an arbitrary path of edges whose end points are the vertices with numbers
1 and
$l $ and
the absolute value of the weight of each edge equals 1.
$$
\includegraphics[width=50mm]{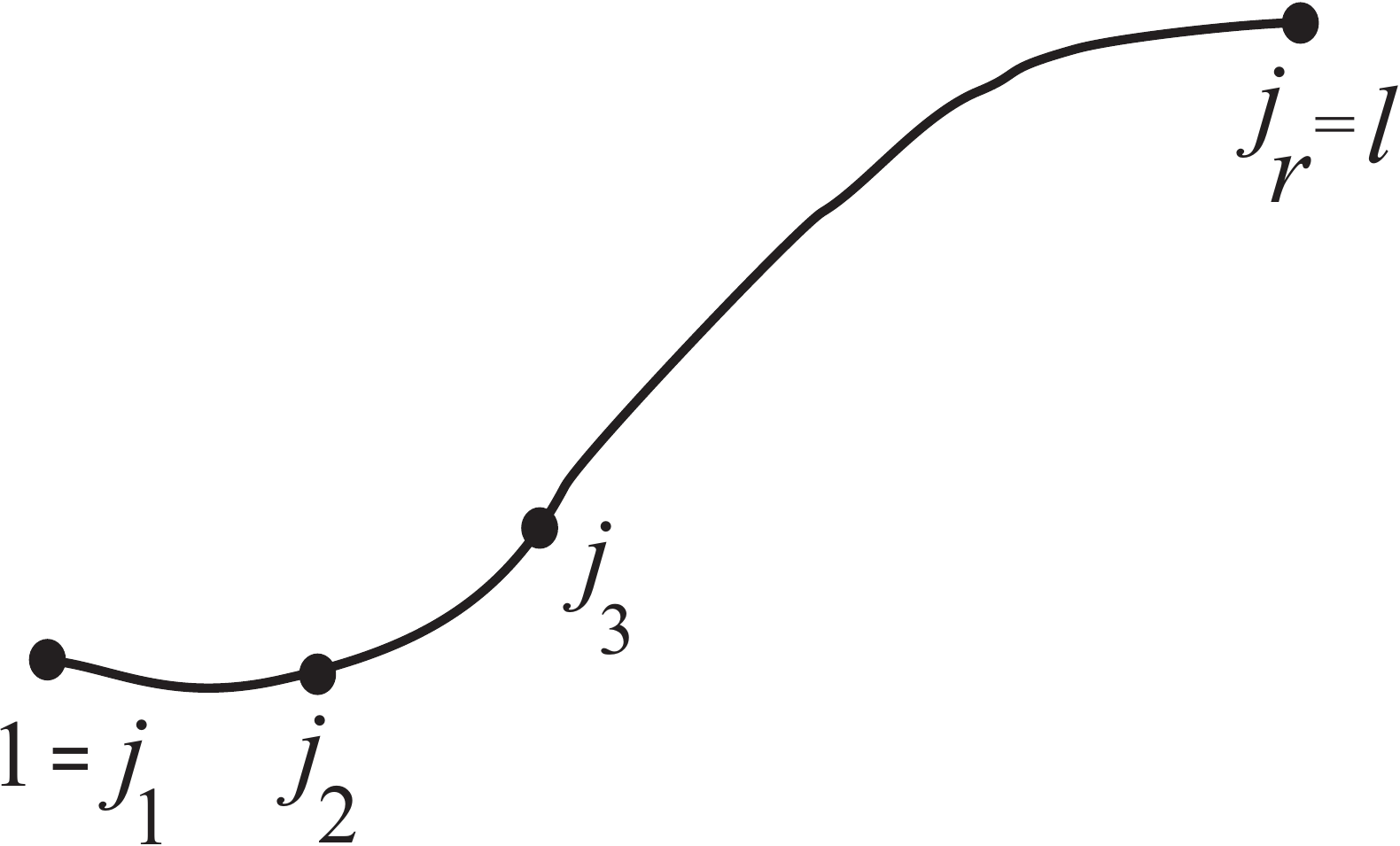}
$$
Put
$$
\Lambda _l =(\ii-R_{j_r})\cdots (\ii-R_{j_3})(\ii-R_{j_2})\Lambda _1
=\Delta \bigg(\prod\limits _{k=2}^r(1-\mu _{j_k})\langle e_{j_{k-1}}\,|\,e_{j_k}\rangle \!\bigg)e_{l}.
$$

\textit{We claim that}
\begin{enumerate}[\hskip 7.2mm\rm a)]\itemsep -.1ex
\item
$\Gamma =\Lambda _1+\cdots +\Lambda _n$ \textit{is a~$K$-invariant root lattice in $V$\;of
rank $2n$}.
\item
$\Gamma $ \textit{does not depend on the construction {\rm(}i.e., on the choice of the paths{\rm)}}.
\item
\textit{Each $K$-invariant root lattice in $V$\;of
rank $2n$ is obtained in this way}.
\end{enumerate}

\begin{proof}
Assertions b) and c) follow from Corollary \ref{Uniq}.
Let us prove a). We check conditions a) and b) of Theorem~\ref{2lat}.
Condition a) is clearly fulfilled, so we only need to check condition b).

We use the notation of this condition. By construction, we have
\begin{align*}
\Lambda_k&=(\ii-R_{l_p})\cdots (\ii-R_{l_2})\Lambda _1,\\
\Lambda _j&=(\ii-R_{r_q})\cdots (\ii-R_{r_2})\Lambda _1,
\end{align*}
for some sequences $1=l_1, l_2,\ldots, l_{p-1}, l_p=k$ and
$1=r_1, r_2,\ldots, r_{q-1}, r_q=j$
such that the $|c_{l_{d-1}l_d}|=|c_{r_{t-1}r_t}|=1$ for all $d$ and $t$:
$$
\includegraphics[width=80mm]{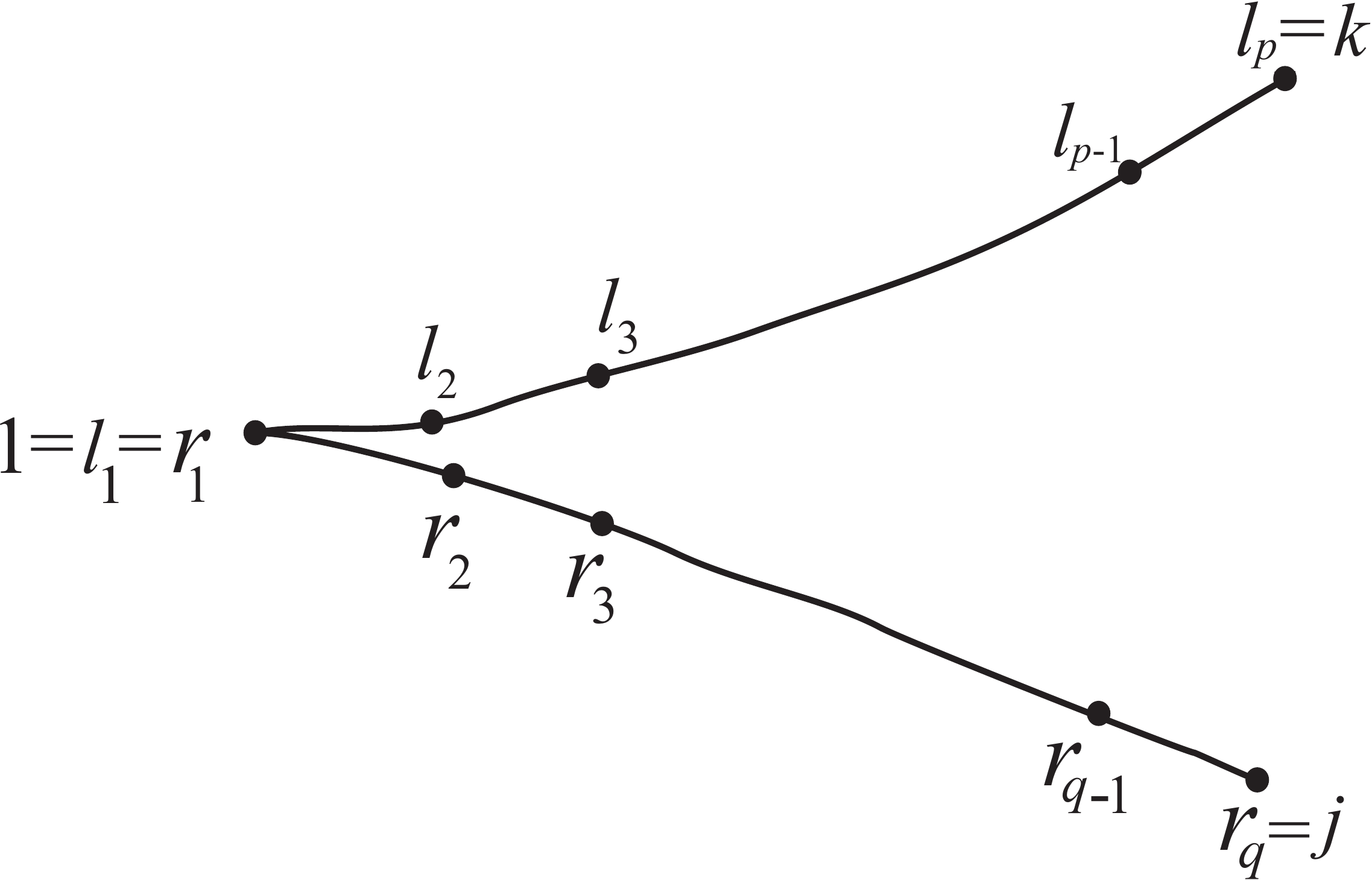}
$$
Let $P=(\ii-R_{l_1})\cdots (\ii-R_{l_{p-1}})$. Then
$$
P \Lambda _k=c_{l_1l_2 \ldots l_{p-1}l_pl_{p-1}\ldots l_2}\Lambda _1 \subseteq \Lambda _1
$$
(thanks to the construction of $\Lambda _1$). But
$$
c_{l_1l_2 \ldots l_{p-1}l_pl_{p-1}\ldots l_2}=c_{l_1l_2}c_{l_2l_3}\cdots c_{l_{p-1}l_p}.
$$
Hence $|c_{l_1l_2 \ldots l_{p-1}l_pl_{p-1}\ldots l_2}|=1$, and therefore, $$P \Lambda _k=\Lambda _1.$$

Let us consider $(\ii-R_k)\Lambda _j$. We have
\begin{align*}
P(\ii-R_k)\Lambda _j&=(\ii-R_{l_1})\cdots (\ii-R_{l_{p-1}})(\ii-R_{l_p})(\ii-R_{r_q})\cdots (\ii-R_{r_2})\Lambda _1\\
&=c_{r_1 r_2\ldots r_{q}l_pl_{p-1} \ldots l_2}
\Lambda _1 \subseteq \Lambda _1.
\end{align*}
Therefore,
$$
P(\ii-R_k)\Lambda _j \subseteq \Lambda _1=P \Lambda _k.
$$
But the restriction of $P$ to $\Cee e_k$ has trivial kernel, because $P \Lambda _k=\Lambda _1$. Therefore, $(\ii-R_k)\Lambda _j \subseteq \Lambda _k$.
The same arguments show that
$(\ii-R_j)\Lambda _k \subseteq \Lambda _j$.
\end{proof}

Therefore, we only need to find all lattices $\Delta$
in
$\Cee $ of rank 2 such that $\Zee  [ \Tr K]\Delta\subseteq \Delta$.\;If $\Zee  [ \Tr K]=\Zee$, then this condition holds for every lattice of rank 2 in $\mathbb C$, and such a lattice is similar to a unique lattice $[1, \tau]$ with $\tau\in \Omega$ (see, e.g., [9, Chap.\,II, \S7, Sect.\,7, Rem.]). By Theorem \ref{cirg2}, if  $\Zee  [ \Tr K]\neq\Zee$, then $\Zee  [\Tr K]$ is an order of a purely imaginary quadratic extension of $\mathbb Q$. It is not difficult to see that
in his case $\Delta$ is similar to an ideal in $\Zee  [ \Tr K]$.\;As Table 4 shows that $\Zee  [ \Tr K]$ is a principal ideal domain (see b) in Section \ref{fre}),
this implies that, up to similarity,
$$
\Delta =\Zee  [ \Tr K] \quad \text{if}\quad \Zee  [\Tr K]\ne \Zee  .
$$

\subsubsection{Example}
$K=K_2$, type $G(4,4,s)$, $s \geqslant 3$. The graph is
$$
\includegraphics{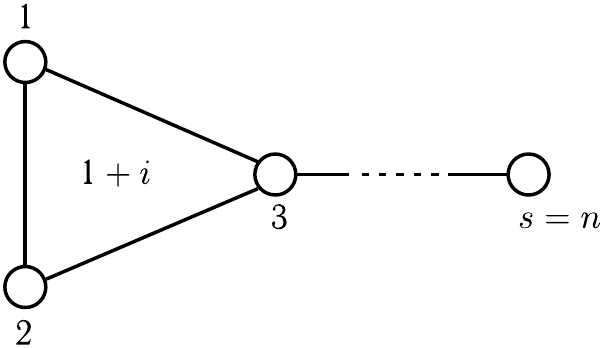}
$$


For the choosen path
connecting vertices numbered 1 and $l \geqslant 3$ (respectively, $1$ and $2$)
the sequence of numbers of its vertices
is $1,3,4,\ldots ,\ell $
(respectively, $1,3,2$).

Here $\Zee  [\Tr K] = \Zee  [i]$ and
$$
\langle e_1 \mid e_3 \rangle =\langle e_2 \mid e_3 \rangle =\langle e_3 \mid e_4 \rangle =\ldots =\langle e_{n-1}\mid e_n \rangle =-\frac{1}{2}.
$$

Hence
\begin{align*}
\Delta &=[1,i],\\
\Lambda _1&=[1,i]e_1,\\
\Lambda _l&=[1,i](1-(-1))^{l -2}\left (-\frac{1}{2} \right )^{l -2}e_l =[1,i]e_l\;\; \mbox{for  $l \geqslant 3$},\\
\Lambda _2&=[1,i](1-(-1))^2 \left (-\frac{1}{2} \right )^2\!e_2=[1,i]e_2.
\end{align*}
Therefore,
$$
\Lambda =[1,i]e_1+\cdots +[1,i]e_n.
$$
\vskip 2mm


Considering in a similar way each of the groups $K$ listed at the beginning of Section \ref{Case1}, we obtain exactly the lattices listed in column $\Tran W$ of Table 2 for the case of $\Lin W = K$.

\subsubsection{Algorithm for constructing \boldmath $K$-invariant full rank root lattices: Case 2}

We consider now the remaining irreducible finite $r$-groups $K$, i.e., the groups

\vskip 2mm

$K_2$, type $G(m,1,s)$ for $s \geqslant 2$ and $m = 2,3,4$, type $G(6,6,2)$;\vskip 1mm

$K_5$; $K_{26}$; $K_{28}$.

\vskip 2mm

We see that the graph of $K$ in these cases is \textit{a chain}. Taking a~suitable numbering, we can assume that $c_{12},c_{23},\ldots ,c_{n-1,n}$ are the only nonzero $c_{jk}$ (the numbering in Table~1
has this property).

\eject

\noindent\textbf{Algorithm}

\vskip 1.1mm

Take an arbitrary $\Zee  [\Tr K]$-invariant
lattice $\Delta_1$ of rank 2 in $\Cee $
(i.e., such that $\Zee  [\Tr K]\Delta_1\subseteq \Delta_1$). By Theorem \ref{cirg3}, we  have
$$
\Delta _1 \subseteq c_{12}^{-1}\Delta _1.
$$
Next, take an arbitrary $\Zee  [ \Tr K] $-invariant lattice $\Delta _2$ between these two lattices (such a~lattice exists, e.g., $\Delta _1$ has this property):
$$
\Delta _1 \subseteq \Delta _2 \subseteq c_{12}^{-1}\Delta _1.
$$
And so on:
\begin{gather*}
\Delta _2 \subseteq \Delta _3 \subseteq c_{23}^{-1}\Delta _2,\\
\ldots \\
\Delta _{n-1}\subseteq \Delta _n \subseteq c_{n-1,n}^{-1}\Delta _{n-1}.
\end{gather*}
For every $l=2,\ldots, n$, put
$$
\Lambda _l =\Delta _l (\ii-R_l )(\ii-R_{l -1})\cdots (\ii-R_2)e_1
=\Delta _l \bigg(\prod\limits _{j=2}^l (1-\mu _j)\langle e_{j-1}\,|\,e_j \rangle \!\bigg)e_l.
$$


\textit{We claim that}
\begin{enumerate}[\hskip 7.2mm\rm a)]\itemsep -.1ex
\item
$\Gamma =\Lambda _1+\cdots +\Lambda _n$ \textit{is a~$K$-invariant root lattice in $V$\;of
rank $2n$}.
\item
\textit{Each $K$-invariant root lattice in $V$\;of
rank $2n$ is obtained in this way}.
\end{enumerate}

\vskip 2mm
\noindent{\it Proof.}
Let us check the conditions a) and c) of Theorem~\ref{2lat}.
In fact we only need to check~c), because a) is obvious. We have
\begin{align*}
\Lambda _l &=\Delta _l (\ii-R_l )(\ii-R_{l -1})\cdots (\ii-R_2)e_1,\\
\Lambda _{l +1}&=\Delta _{l +1}(\ii-R_{l +1})(\ii-R_l )\cdots (\ii-R_2)e_1.
\end{align*}
Therefore,
\begin{align*}
\hskip 5mm (\ii-R_{l +1})\Lambda _l &=\Delta _l (\ii-R_{l +1})\cdots (\ii-R_2)e_1\\
&
\subseteq \Delta _{l +1}(\ii-R_{l +1})\cdots (\ii-R_2)e_1\\
&=\Lambda _{l +1}\subseteq c_{l , l +1}^{-1}\Delta _l (\ii-R_{l +1})\cdots (\ii-R_2)e_1\\
&=c_{l , l +1}^{-1}(\ii-R_{l +1})\Lambda _l .\hskip 14mm \text{\textsquare}
\end{align*}

The same argument is in case 1 shows that, up to similarity,
\begin{align*}
\Delta_1 =\begin{cases}
\Zee  [ \Tr K]  &\text{if  $\Zee  [\Tr K]\ne \Zee$},\\
[1,\tau]\;\; \mbox{for any $\tau\in \Omega$} &\text{if  $\Zee  [\Tr K]=\Zee$}.
\end{cases}
\end{align*}

\subsubsection{Example}
$K = K_2$, type $G(3,1,s)$, $s \geqslant 2$. The graph is
$$
\includegraphics{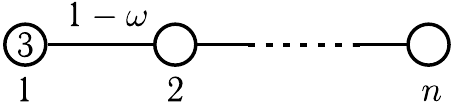}
$$
\eject
We have $\Zee  [ \Tr K] = \Zee  [\omega]$ and
\begin{gather*}
\langle e_1 \mid e_2 \rangle =1/\sqrt 2,\quad \langle e_2 \mid e_3 \rangle =\ldots =\langle e_{n-1}\mid e_n \rangle =-\frac{1}{2},\\
c_{12}=1-\omega ,\quad c_{23}=\ldots =c_{n-1,n}=1.
\end{gather*}
Therefore,
\begin{align*}
(\ii-R_l )\cdots (\ii-R_2)e_1&=(1-(-1))^{l -1}\left (\!-\frac{1}{2} \right )^{l -2}\frac{1}{\sqrt 2}\,e_l\\
& =(-1)^l \sqrt 2\,e_l \;\;\;\mbox{for $l =2,\ldots ,n.$}
\end{align*}
We have
\begin{align*}
&\Delta _1=[1,\omega ],\\
&\Delta _2=\ldots =\Delta _n,\\
[1,\omega ]\subseteq &\Delta _2 \subseteq (1-\omega )^{-1}[1,\omega ].
\end{align*}

But $|1-\omega |^2=3$. Hence
$$
\Lambda _2=[1,\omega ]\quad \text{or}\quad (1-\omega )^{-1}[1,\omega ]=\frac{i}{\sqrt 3}\,[1,\omega ].
$$
Thus, we have only two possibilities: either
$$
\Lambda =[1,\omega ]e_1+[1,\omega ]\sqrt 2\,e_2+\cdots +[1,\omega ]\sqrt 2\,e_n,
$$
or
$$
\Lambda =[1,\omega ]e_1+[1,\omega ]i \sqrt {\frac{2}{3}}\,e_2+\cdots +[1,\omega ]i \sqrt {\frac{2}{3}}\,e_n.
$$


In this way, for the groups $K$ under consideration, we obtain all the lattices listed in Table 2 in column ${\rm Tran}\,W$  for the case $\Lin W=K$, however,
if $K = K_2$, types $G(2,1,2)$, $G(6,6,2)$, $K = K_5$, and $K=K_{28}$, then apart of them, we obtain also some additional lattices. They are not listed in Table 2 because they do not give new $r$-groups, i.e., the semidirect product of $K$ and such a lattice is equivalent to one of the groups from Table 2. We skip the details.

\subsection{Invariant lattices in the case \boldmath$s = n+1$}\label{4.8}

\textit{Now we briefly consider the case $s = n+1$}. These are the groups:
$$
G(4,2,n),\quad G(6,2,n),\quad G(6,3,n),\quad K_{12},\quad \text{and}\quad K_{31}.
$$
We explain the approach by several examples.

\subsubsection{Examples}\label{refff}
\

a) $K = K_2$, type $G(6,2,n)$ or $G(6,3,n)$.


The graphs of these groups are
$$
\begin{array}{ll}
G(6,2,n)&\begin{array}{@{}l@{}}
\includegraphics[width=80mm]{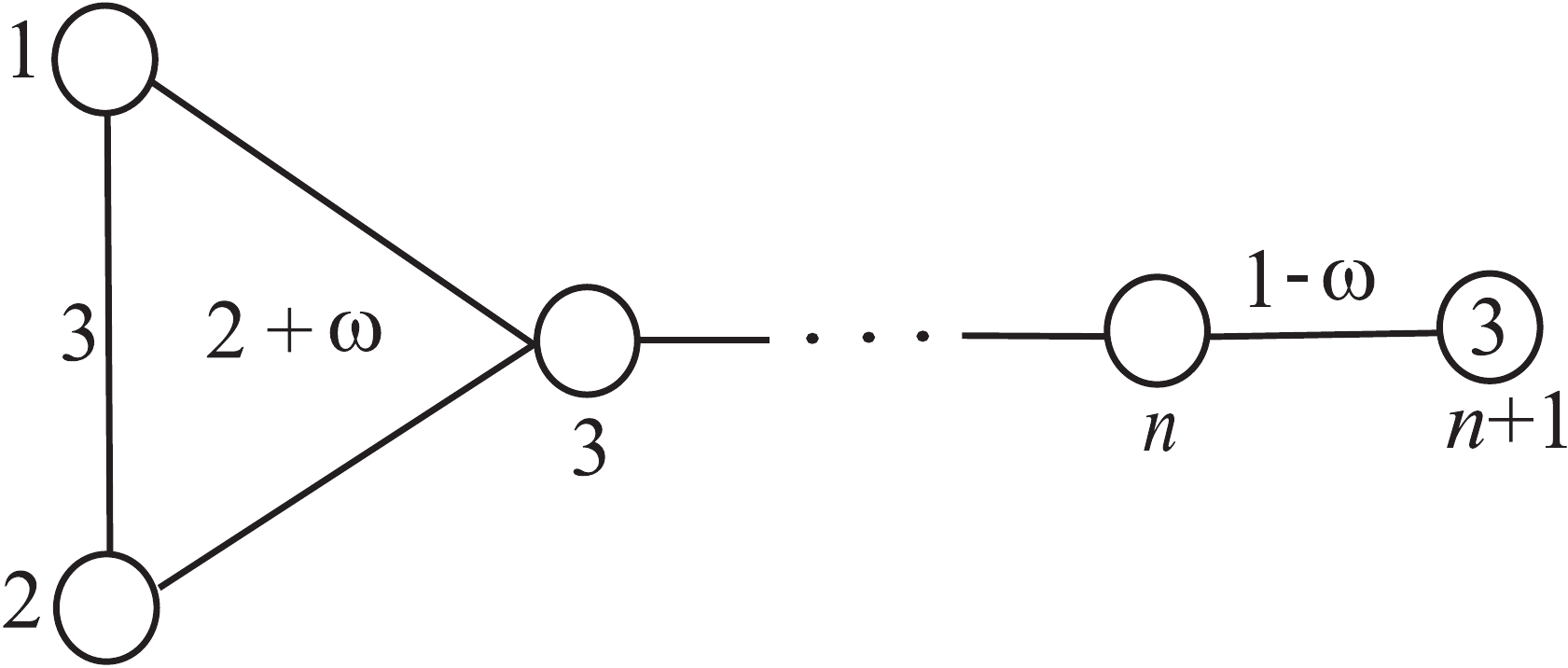}
\end{array}
\\
G(6,3,n)&\begin{array}{@{}l@{}}
\includegraphics[width=80mm]{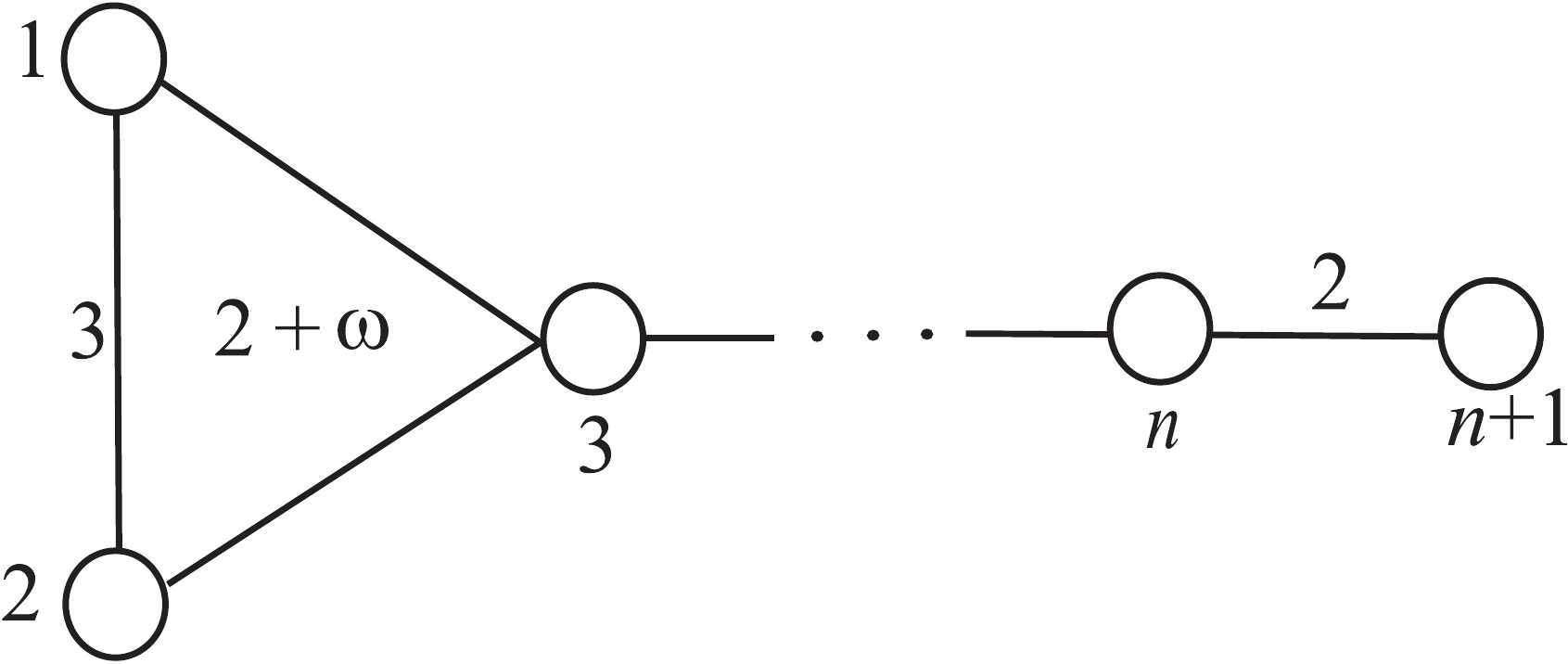}
\end{array}
\end{array}
$$

Using the notation of Section \ref{sn+1}, we see that $K'=K_2$, type $G(6,6,n)$, see Table~1.

Applying  the described algorithm (case 1),
we see
that, up to similarity, there exists only one $K'$-invariant full rank \textit{root} lattice, namely,
$$
\Gamma =[1,\omega ]e_1+\cdots +[1,\omega ]e_n.
$$
But for $K'$, we have
$$
{\mbox{\rm det}}\, S=4 \,{\mbox{\rm sin}}^2 \frac{\pi}{m} \Bigr |_{m=6}=1,
$$
see Example b) in Section \ref{exa}.
It follows from Corollary \ref{Cooo}
that, up to similarity, $\Gamma $ is a~unique $K'$-invariant full rank lattice. But $K$ \textit{must} have an invariant lattice of full rank thanks to Theorem~\ref{cirg2}.
Therefore, this lattice has to be $\Gamma $ (of course, one can also verify $K$-invariance of
$\Gamma $
by a
direct computation).

\vskip 1mm

b) $K=K_2$, type $G(4,2,n)$ for $n \geqslant 3$. The graph is
$$
\includegraphics[width=80mm]{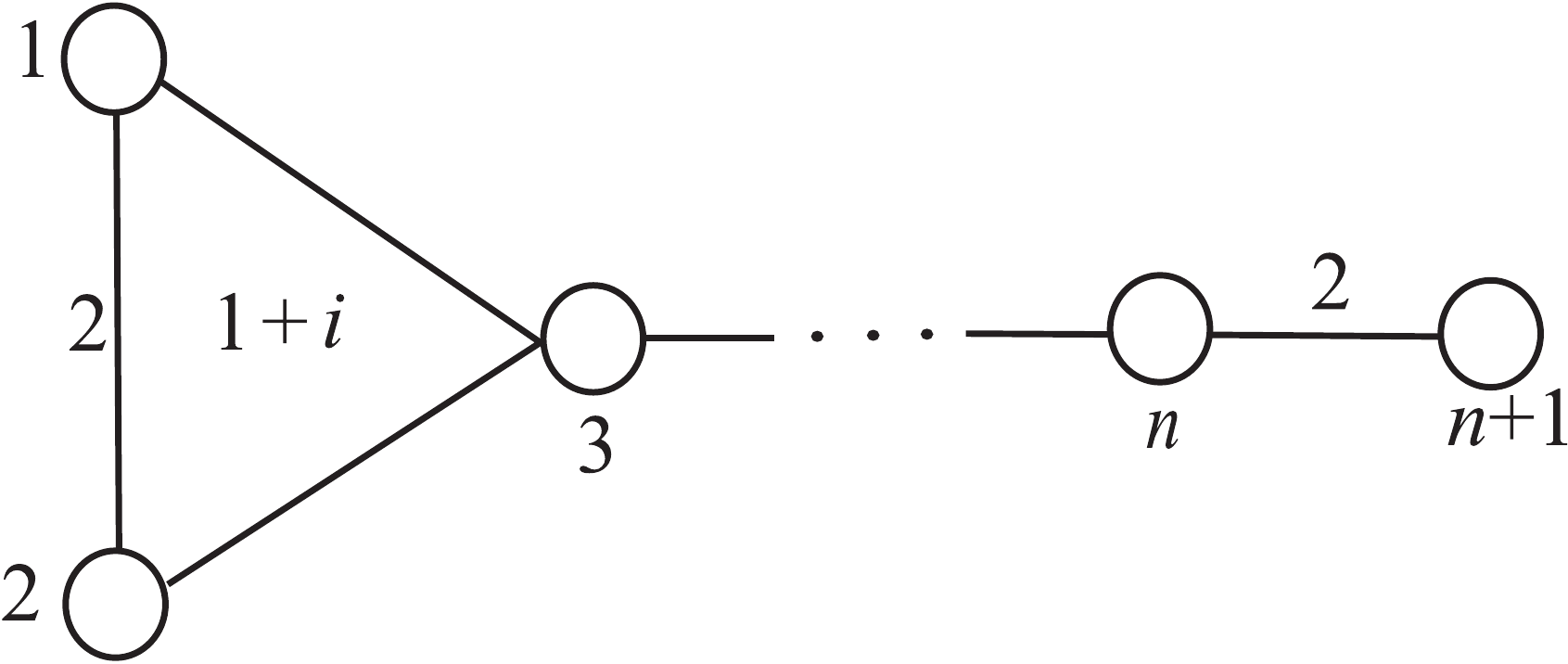}
$$
We see that $K' = K_2$, type $G(4,4,n)$.
As in the previous example we see that there exists only one (up to similarity) $K'$-invariant full rank root lattice in $V$, namely,
$$
\Lambda =[1,i]e_1+\cdots +[1,i]e_n.
$$

So, to describe all $K$-invariant full rank lattices in $V$ we need to find all lattices $\Gamma $ in $V$ which share the following properties:
\begin{enumerate}[\hskip 7.2mm\rm a)]\itemsep -.1ex
\item
$\Gamma ^0=\Lambda $ (with respect to $K'$),
\item
$\Lambda \subseteq \Gamma \subseteq \Lambda ^*=S^{-1}\Lambda $ (with respect to $K'$);
\item
$(\ii-R_{n+1})\Gamma \subseteq \Gamma $.
\end{enumerate}

\noindent
We see that
every
vertex of the graph of $K$ is endpoint of
an edge with weight 1. Hence $(\Lambda ^*)^0=\Lambda $, see Theorem~\ref{***}.
Therefore, b) $\Rightarrow $ a).

The matrix of $S$ with respect to the basis $e_1,\ldots ,e_n$ has the form
$$
\left (\begin{array}{ccrrcrr}
2&1-i&-1\\
1+i&2&-1\\
&-1&2&-1\\
&&-1&2\\
&&&&\ddots\\
&&&&&2&-1\\
&&&&&-1&2
\end{array}
\right ).
$$
We have ${\mbox{\rm det}}\, S = 4 {\mbox{\rm sin}}^2 \frac{\pi}{4}=2$, see Section  \ref{4.4}. Therefore, the coefficients of $S^{-1}$ lie in $\Zee  \left [i,\frac{1}{2} \right ]$ and $\Lambda ^*/\Lambda $ is the groups of order 4
(see\,Theorem~\ref{proS}). 
These facts imply that
\[
\Lambda ^*/\Lambda \simeq \Zee  /2 \Zee  \oplus \Zee  /2 \Zee .
\]
It is readily checked that $e_j=Sf_j$ for $j=1,2$, where
\begin{align*}
f_1&=\frac{n}{2}\,e_1+\frac{-1-(n-1)i}{2}\,e_2+\frac{(1-i)(n-2)}{2}\,e_3+
\frac{(1-i)(n-3)}{2}\,e_4+\cdots +\frac{(1-i)1}{2}\,e_n,\\
f_2&=\frac{-(n-1)-i}{2}\,e_1+\frac{ni}{2}\,e_2+\frac{(i-1)(n-2)}{2}\,e_3+
\frac{(i-1)(n-3)}{2}\,e_4+\cdots +\frac{(i-1)1}{2}\,e_n.
\end{align*}
From $f_1, f_2\in S^{-1}\Lambda=\Lambda^*$ and $f_1-f_2\notin \Lambda$ we infer that
$f_1$, $f_2$, $f_3=f_1+f_2$ are representatives of different nonzero elements of $\Lambda ^*/\Lambda $. Therefore, every $K$-invariant full rank lattice is similar to one of the folowing lattices:
\[
\left.\begin{array}{l}
\Lambda ;\\[-1mm]
\Lambda \cup (\Lambda +f_j)\quad \text{for}\quad j=1,2,3;\\[-1mm]
\Lambda \cup \bigcup\limits _{j=1}^3(\Lambda +f_j)
\end{array}\right\}
\eqno\boxed{5}
\]
Using the equalities
\begin{align*}
e_{n+1}&=\frac{-1-i}{\sqrt 2}\,e_1+\frac{-1+i}{\sqrt 2}\,e_2-\sqrt 2\,e_3+\sqrt 2\,e_4-\ldots -\sqrt 2\,e_n,\\
(\ii-R_{n+1})v&=\begin{cases}
0 &\text{if}\quad v=e_1,\ldots ,e_{n-1},f_3,\\
\sqrt 2\,e_{n+1}&\text{if}\quad v=e_n,\\
\frac{1-i}{\sqrt 2}\,e_{n+1}&\text{if}\quad v=f_1,\\
\frac{i-1}{\sqrt 2}\,e_{n+1}&\text{if}\quad v=f_2,
\end{cases}
\end{align*}
one can
directly verify that \textit{all lattices $\boxed{5}$
are $K$-invariant}.

The same considerations can be carried out for other groups $K$ from the above list and obtain the description of \textit{all} (up to similarity) $K$-invariant full rank lattices. We leave it
to the reader (the most complicated case is $K = K_2$, type $G(4,2,2)$).

\section{\bf The structure of \boldmath $r$-groups in the case $s=n+1$}\label{Sr-s}

When $s = n+1$, it is no longer true, in general, that an infinite complex irreducible crystallographic $r$-group $W$ is the semidirect product of $\Lin W$ and $\Tran W$. We explain here how one can find the corresponding extensions of $\Tran W$ by $\Lin W$ in this case. We assume that $\ka = \Cee $.

\subsection{The cocycle $c$}\label{cocyc}

We recall that the structure of extension is given by a~$V$-valued 1-cocycle $c$ of $\widetilde K$, where $K = \Lin W$; see Section \ref{cohom} the notation of which we use.

Taking a point $a\in E$ as the origin, we can identify $A(E)$ and $\GL(V) \ltimes V$; then $W$ consists of the elements $(P,c(P)+v
)$, $P\in K$, $v\in \Tran W$.

We assume that the order of $R_j$ is equal to $m(H_{R_j} )$ for all $1 \leqslant j \leqslant s$, see Section \ref{auxi}.


First of all, we note the following:

\vskip 2mm

\textit{
Replacing $c$ by a~suitable cocycle cohomologous to $c$, one can assume that}
$$
c(r_1)=\ldots =c(r_n)=0.
$$

\begin{proof} For every $j=1,\ldots n$, let
$
\gamma_j\in W$ be a
reflection with $\Lin \gamma _j=R_j$ (see Theorem \ref{rR}).
We assume, as usual, that the group $K'$ generated by
$R_1,\ldots, R_n$ is irreducible. We have $\bigcap
_{j=1}^nH_{R_j}=0$. Hence $\bigcap
_{j=1}^nH_{\gamma _j}$ is a single point of $E$, say $b$.
Then  we have
$$
\kappa _b(\gamma _j)=(R_j,0),\quad j=1,\ldots ,n,
$$
and we are done.
\end{proof}

Given this, we assume now that
$$
c(r_1)=\ldots =c(r_n)=0.
$$

Therefore, $c$ is defined by only one vector $c(r_{n+1})$. Moreover, one can take
$$
c(r_{n+1})=\lambda e_{n+1}\quad \text{for some}\quad \lambda \in \Cee ,
$$
because there exists a~reflection $(R_{n+1},v)$ in $W$ (see Theorem \ref{rR} and Proposition \ref{prope}).

Therefore, the problem can be reformulated as follows:
 given a $K$-invariant full rank lattice $\Gamma $ in $V$,
find all
$\lambda \in \Cee $ such that
the following properties hold:

\vskip 1.5mm

a) The $V$-valued cocycle $c$ of $\widetilde K$, given by the equalities
\begin{equation*}
\left.\begin{split}
c(r_1)&=0,\\[-1mm]
\ldots&\ldots\\
c(r_n)&=0,\\
c(r_{n+1})&=\lambda e_{n+1}
\end{split}
\right\}
\eqno\boxed{6}
\end{equation*}
satisfies the condition
$$
c(F) \in \Gamma\;\; \mbox{for every relation $F\in \widetilde K$ of $K$}
$$
(i.e., for every $F\in \Ker \phi $, see Section \ref{cohom}).

\vskip 1mm

b) The extension $W$ of $\Gamma$ by $K$ determined by
this cocycle $c$
is an $r$-group.

\subsubsection{Theorem}\label{T??}

{\it
\begin{enumerate}[\hskip 7.2mm\rm 1)]\itemsep -.1ex
\item
If $\Gamma =\Gamma ^0$, then {\rm a)} implies
{\rm b)}.
\item
If $c(F) \in \Gamma ^0$ for every relation $F\in \widetilde K$ of $K$, then {\rm b)} implies the equality $\Gamma =\Gamma ^0$.
\end{enumerate}
}

\begin{proof}
1) Let $\Gamma =\Gamma ^0$. We know that
$$
\Gamma ^0=\Gamma _1+\cdots +\Gamma _{n+1}.
$$
But $\Gamma '=\Gamma _1+\cdots +\Gamma _n$ is a~root lattice for $K'$ and the condition
$$
c(r_1)=\ldots =c(r_n)=0
$$
shows that the semidirect product $W'$ of $K'$ and $\Gamma '$ lies in $W$. Theorem \ref{ssp} implies that $W'$
is an $r$-group.

By construction,
$W$ contains the following set of reflections:
\[
(R_1,0),\ldots, (R_n, 0),\;\;
\mbox{and\;\; $(R_{n+1},\lambda e_{n+1}+t)$  for
for each $t\in \Gamma _{n+1}$}.
\eqno\boxed{7}
\]

Take an element $\gamma =(P,v)\!\in\! W$.
As  $P$ is a~product of some elements of the set
$R_1,\ldots, R_{n+1}$,
multiplying $\gamma $ by some elements of set $\,\boxed{7}$\,,
we can obtain an element of the form $(\ii, t)$.\;As each element $(\ii,t')$ for $t'\in \Gamma '$ lies in $W'$, it is also a product of reflections from $W$, because $W'$ is an $r$-group. This proves that, multiplying $\gamma $ by reflections, we can obtain $(\ii,t)$ for some $t\in \Gamma _{n+1}$. But as $c$ is a cocycle, the reflection $(R_{n+1}^{-1},-\lambda R_{n+1}^{-1}e_{n+1}+t)$ lies in $W$.
From this and the equality
$$
(R_{n+1}^{-1},-\lambda R_{n+1}^{-1} e_{n+1}+t)(R_{n+1},\lambda e_{n+1})=(\ii, t).
$$
we then infer that
$\gamma $ is a~product of reflections. Hence $W$ is an $r$-group.

2) Let $c(F) \in \Gamma ^0$ for every relation $F\in\widetilde K$ of $K$.\;Then the cocycle $c$ defines, in fact, a 1-cocycle of $K$ with values in $V/\Gamma ^0$. Let $W'$ be the group defined by $c$, with $\Lin W' = K$ and $\Tran W' =\Gamma ^0$. It is an $r$-group because of 1). Besides, we
have the~group $W$ defined by $c$, with $\Lin W =K$ and $\Tran W=\Gamma $.

Let $\gamma \in W$ be a~reflection. By Theorem \ref{ctpgr},
there exists $\delta \in W'$
such that $\delta \gamma \delta ^{-1}=(R_j^l ,t)$ for certain $l$, $j$, $t$. As $\delta \gamma \delta ^{-1}$ is also a~reflection, we have $t \perp H_{R_j}$. But $t=c(r_j^l )+v$ for some $v\in \Gamma $. We have
$$
c(r_j^l )=\big(\ii+R_j+R_j^2+\cdots +R_j^{l -1}\big)c(r_j).
\eqno\boxed{8}
$$
By the definition of $c$ (see $\boxed{6}$), we have $c(r_j)\perp H_{R_j}$.\;In view of  $\,\boxed{8}$\,, this yields $c(r_j^l )\perp H_{R_j}$. There\-fore, $v \perp H_{R_j}$, or, in other words, $v\in \Gamma ^0$. This means that $\delta \gamma \delta ^{-1}=W'$. As $\delta \in W'$, this yields $\gamma \in W'$. Therefore, if $W$ is an $r$-group, then $W =W'$ and $\Gamma =\Gamma ^0$.
\end{proof}

The following simple observation is very useful in practice because it gives strong restrictions on the choice of $\lambda $:

\subsubsection{Theorem}\label{lamb}
{\it
Let condition {\rm a)} of Section {\rm \ref{cocyc}}
holds and let $P\in K'$ be an element  such that
$$
R_{n+1}PR_{n+1}^{-1}\in K'.
$$
Then
$$
\lambda (\ii-R_{n+1}PR_{n+1}^{-1})e_{n+1}\in \Gamma .
$$
}

\begin{proof}
It follows from the definition of $c$ (see $\boxed{6}$) that $(P,0)$ and $(R_{n+1},\lambda e_{n+1}) \in W$.
The equality
$$
(R_{n+1},\lambda e_{n+1})(P,0)(R_{n+1},\lambda e_{n+1})^{-1}
=(R_{n+1}PR_{n+1}^{-1}, -R_{n+1}PR_{n+1}^{-1}\lambda e_{n+1}+\lambda e_{n+1}).
$$
now implies the assertion made.
\end{proof}

\subsubsection{Example}\hskip -2mm

$K=K_{31}$. The graph is
$$
\includegraphics{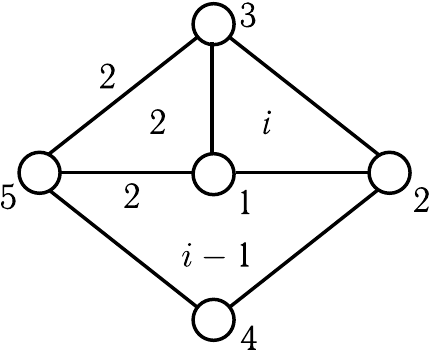}
$$
The vectors $e_1,e_2,\ldots ,e_5$ are given in Table~\ref{tab2}. Note that $e_5=ie_1+e_2+e_3$.

There exists only one (up to similarity) $K$-invariant full rank lattice $\Gamma $ in $V$, namely,
$$
\Gamma =[1,i]e_1+\cdots +[1,i]e_4.
$$

It is known (see [4]) that  for the system of generators $R_1,\ldots, R_5$, the presentation of $K$
is given by the following relations:
\begin{equation*}
\left.
\begin{split}
&r_1^2,\; r_2^2,\; r_3^2,\; (r_2r_3)^3,\; (r_3r_1)^3,\; (r_1r_2)^3,\\
&(r_2r_1r_3r_1)^4,\\
&(r_4r_5)^3,\\
&r_5^2,\; (r_5r_2)^2,\; (r_5r_1r_3r_1)^2,\; (r_5r_3)^4,\;
r_1(r_5r_3r_2r_3)r_1(r_5r_3r_2r_3)^{-1},\\
&r_4^2,\;(r_4r_1)^2,\;(r_4r_3)^2,\;(r_4r_2)^3
\end{split}
\right\}
\eqno\boxed{9}
\end{equation*}
(i.e., $\Ker \phi $ is
the normal closure in $\widetilde K$ of the set $\boxed{9}$; see Section \ref{cohom}).

We consider cocycle $\boxed{6}$.
From the equalities
$$
R_1(R_5R_3R_2R_3)R_1(R_5R_3R_2R_3)^{-1}=R_5^2=\ii
$$
it follows that $R_5R_1R_5\in K'$.
Therefore, by Theorem \ref{lamb}, if condition a) holds, then
$$
\Gamma \ni \lambda (\ii-R_5R_1R_5)e_5=\lambda ((1+i)e_1+2e_2+2e_3).
$$
Hence $\lambda =
{(a+bi)}/2$ for $a,b\in \Zee  $, and $a
\equiv b \;({\rm mod}\; 2)$.\;As $[1,i]e_5\in \Gamma $,
this reduces our con\-si\-de\-ra\-ti\-ons to checking whether $c(F)\in \Gamma$ holds or not for each relation $F$ listed in $\boxed{9}$  and\;for
$$
\lambda =\frac{1+i}{2}.
$$
This
is done
by direct
computations:
\begin{align*}
c(r_5)^2&=c(r_5)+r_5c(r_5)=\frac{1+i}{2}\,(1+r_5)e_5=0\in \Gamma ,\\
c((r_4r_5)^3)&=(1+r_4r_5+(r_4r_5)^2)(c(r_4)+r_4c(r_5))\\
&=\frac{1+i}{2}\,(1+r_4r_5+(r_4r_5)^2)(r_4e_5)=0\in \Gamma ,
\end{align*}
and so on (one only needs to consider the relations which involve $R_5$).\;This checking shows that $\lambda =
{(1+i)}/2$ indeed gives a~cocycle of $K$,
and hence defines an $r$-group $W$ with $\Lin W = K$ and $ \Tran W = \Gamma $.
It can be straightforwardly verified  that this cocycle
is not a~coboundary, i.e., $W$ is not a~semidirect product.

\subsection{}\label{fis} \hskip -2mm
The same considerations can be given for each $K$ with $s = n + 1$, and, as a~result, one obtains Table~\ref{tab2}.

\textit{A posteriori} it appears that in all the cases either $\Gamma $ is a~root lattice, or $c(F) \in \Gamma ^0$ for every relation $F$ of $K$. Therefore,
the following theorem holds:

\subsubsection{Theorem}\label{}
{\it
$\Tran W$ is
a~root lattice for every $r$-group $W$.
}


\vskip 2mm

Similar straightforward computations yield
the proofs of
Theorems \ref{th1} and \ref{th2}.

As for completing the proof of Theorem~\ref{cirg3}
(whose first part is given in Section \ref{4.6}), its second part (about minimality of $\Zee  [\Tr K]$) follows from \cite{9}, and
third part
from Table~1.

\EditInfo{April 26, 2023}{July 13, 2023}{Dimitry Leites and Sofiane Bouarroudj}
\end{document}